\documentclass[11pt,leqno, a4paper]{amsart}
\usepackage{amsmath,amsthm,amscd,amssymb,amsfonts, amsbsy}
\usepackage{latexsym}
\usepackage{txfonts}
\usepackage{exscale}
\usepackage{graphicx}
\usepackage{pgf}
\usepackage{color}

\definecolor{purple}{RGB}{128,0,128}

\calclayout%
\allowdisplaybreaks

\theoremstyle{plain}
\newtheorem{theorem}{Theorem}[section]
\newtheorem{lemma}[theorem]{Lemma}
\newtheorem{corollary}[theorem]{Corollary}
\newtheorem*{corollary*}{Corollary}
\newtheorem*{theorem*}{Theorem}
\newtheorem{proposition}[theorem]{Proposition}

\newtheorem{definition}[theorem]{Definition}

\theoremstyle{remark}
\newtheorem{remark}[theorem]{Remark}

\numberwithin{equation}{section}

\newcommand{\rn}{\mathbb{R}^n}

\newcommand{\cH}{\mathcal{H}}

\newcommand{\W}{\mathcal{W}}

\newcommand{\I}{\mathcal{I}}
\newcommand{\cA}{\mathcal{A}}
\newcommand{\cF}{\mathcal{F}}

\newcommand{\cG}{\mathcal{G}}

\newcommand{\Hn}{\mathcal{H}^{n-1}}

\newcommand{\N}{\mathbb{N}}

\newcommand{\St}{\mathbb{S}}
\newcommand{\R}{\mathbb{R}}

\newcommand{\intrn}{\int_{\rn}}

\newcommand{\G}{\mathcal{G}}
\newcommand{\eps}{\varepsilon}

\newcommand{\supp}{{\rm{spt}}}
\newcommand{\spt}{{\rm{spt}}~}

\newcommand{\hm}{\omega}
\newcommand{\A}{\mathcal{A}}

\newcommand{\RNum}[1]{\uppercase\expandafter{\romannumeral #1\relax}}
\newcommand{\B}{\mathcal{B}}

\newcommand{\diam}{\operatorname{diam}}

\newcommand{\dist}{\operatorname{dist}}

\renewcommand{\tilde}{\widetilde}
\renewcommand{\hat}{\widehat}

\newcommand{\lip}{\operatorname{Lip}}

\newcommand{\loc}{\operatorname{loc}}

\renewcommand{\div}{\operatorname{div}}

\newcommand{\re}{\mathbb{R}}

\newcommand{\defeq}{:=}

\newcommand{\om}{\Omega}
\newcommand{\pom}{{\partial\Omega}}
\DeclareMathOperator{\interior}{int}

\newcommand{\ree}{\mathbb{R}^{n+1}}

\newcommand{\dif}{d}

\def\Xint#1{\mathchoice
{\XXint\displaystyle\textstyle{#1}}%
{\XXint\textstyle\scriptstyle{#1}}%
{\XXint\scriptstyle\scriptscriptstyle{#1}}%
{\XXint\scriptscriptstyle\scriptscriptstyle{#1}}%
\!\int}
\def\XXint#1#2#3{{\setbox0=\hbox{$#1{#2#3}{\int}$ }
\vcenter{\hbox{$#2#3$ }}\kern-.6\wd0}}

\def\fint{\Xint-}

\newcommand{\restr}{\mathbin{\vrule height 1.6ex depth 0pt width
0.13ex \vrule height 0.13ex depth 0pt width 1.3ex}}

\title[Flatness and Oscillation]{Two phase free boundary problem for Poisson kernels}
\author{S. Bortz}
\address{Department of Mathematics \\
University of Alabama, Tuscaloosa, AL, 35487, USA}
	\email{sbortz@ua.edu} 
	
\author{M. Engelstein}
\address{Department of Mathematics\\University of Minnesota, Minneapolis, MN, 55455, USA }
\email{mengelst@umn.edu}

\author{M. Goering}
\address{Department of Mathematics, University of Washington, Seattle, WA
98195, USA}
\email{mgoering@uw.edu}

\author{T. Toro}
\address{Department of Mathematics, University of Washington, Seattle, WA
98195, USA}
\email{toro@uw.edu}

\author{Z. Zhao}
\address{Department of Mathematics, University of Chicago, Chicago, IL 60637, USA}
\email{zhaozh@uchicago.edu}

\date{\today}

\thanks{Some of this work was done while S.B., M.E., T.T. and Z.Z. were in residence at the MSRI Harmonic Analysis program, supported under NSF Grant No. DMS-1440140.  During the preparation of this manuscript S.B. was partially supported by the NSF INSPIRE Award DMS-1344235 and M.E. was supported by, NSF DMS-2000288, NSF MSPRF DMS-1703306 and by David Jerison's grant DMS 1500771. M.G. was partially supported by NSF grants DMS-1664867 and DMS-1500098. T.T.  was partially supported by the Craig McKibben \& Sarah Merner Professor in Mathematics, by NSF grant number DMS-1664867, and by the Simons Foundation Fellowship 614610. Z.Z. was partially supported by the Institute for Advanced Study and by NSF grants DMS-1664867 and DMS-1902756.}

\begin{document}
\maketitle

\begin{abstract}
We provide a potential theoretic characterization of vanishing chord-arc domains under mild assumptions. In particular we show that, if a domain has Ahlfors regular boundary, the oscillation of the logarithm of the interior and exterior Poisson kernels yields a great deal of geometric information about the domain. We use techniques from classical calculus of variations, potential theory and quantitative geometric measure theory to accomplish this. {One} feature of this work, compared to \cite{kenigtorotwophase} and \cite{bortz2017singular}, {is that {\it a priori} we only require that the domains in question are connected.}

\end{abstract}

\section{Introduction}\label{s:intro}

Questions concerning the connections between the geometry of a domain and the regularity of its boundary with the potential theoretic properties of the domain, the behavior of singular integrals on the boundary, and the boundary regularity to solutions of elliptic PDEs have generated a flurry of activity in the area of non-smooth analysis (see \cite{toronotices} and \cite{toronotices-2}
for a brief recent history and references). In this paper we focus on the potential theoretic properties of a domain and its complement and explore their ties to the geometry of the domain.
In particular, we show that if $\Omega: =\Omega^+ \subset\rn$ and the interior of its complement $\Omega^-$ are connected,  have a shared boundary which is Ahlfors regular (see Definition \ref{def:AR}), and the logarithm of the Poisson kernel of each domain is in $VMO_{\loc}$, then the unit normal is also in $VMO_{\loc}$ and the domain is vanishing Reifenberg flat (see Definitions \ref{def:Reifenberg} and \ref{def:bmo-vmo}). We contrast our result with those in the literature in order to emphasize the wealth of geometric information  (thus far overlooked) encoded in the assumption concerning the oscillation of the logarithm of the Poisson kernels.

In \cite{kenigtorotwophase} the authors established the following: suppose that $\Omega^\pm$ are chord-arc domains (i.e, NTA domains with Ahlfors regular boundary), and that $k^\pm$ are the Poisson kernels of $\om^\pm$ with poles $X^\pm\in\om^\pm$. If $\log k^\pm \in VMO_{\loc}(\sigma)$ then the unit normal vector $\nu\in VMO_{\loc}(\sigma)$ where $\sigma=\cH^{n-1}\restr\pom$ (see Definition \ref{def:NTA}). 
In particular the assumption that $\Omega^\pm$ are chord-arc domains ensures that $\pom^\pm$ are uniformly rectifiable (see Definition \ref{URdom}).
In \cite{bortz2017singular} the authors relax the geometric conditions: to be more precise, via a novel approach using layer potentials rather than blow ups, they prove that if $\om^\pm\subset \rn$ are domains, whose common
boundary is uniformly rectifiable then $\log k^\pm \in VMO_{\loc}(\sigma)$ implies that $\nu\in VMO_{\loc}(\sigma)$. {  We also mention the recent work Prats-Tolsa \cite{Prats-Tolsa}, where the authors studied a different but closely related problem arising in Kenig-Toro \cite{kenigtorotwophase}. They study the kernel between harmonic measures $\omega^{\pm}$ of $\Omega^{\pm}$, and show that for Reifenberg flat NTA domains, small oscillation for the logarithm of that kernel is also closely linked to small oscillation for the unit normal $\nu$.} 
 
In this paper we further loosen the \textit{a priori} assumption in \cite{kenigtorotwophase} and instead deduce as much geometric information as possible from the regularity of $\log k^\pm$. Furthermore using classical tools from the calculus of variations we establish that in this context the oscillation of the unit normal controls the flatness of the boundary. More precisely, when $\pm \nu$ are outward pointing normal vectors to $\partial \Omega^{\pm}$, $\sigma = \cH^{n-1} \restr \partial \Omega$, and $\omega_{\pm}^{X^{\pm}} = \omega^{\pm}$ is the harmonic measure for $\Omega^{\pm}$ with pole at $X^{\pm}$  we show that:

\begin{theorem}\label{Tchar2sideCAD.thrm}
Let $n \ge 3$ and suppose $\om^+ \subset \rn$ and $\om^-= \rn\setminus \overline{\om^+}$ are domains satisfying $\partial\Omega:= \pom^+ = \pom^-$,
and that $\pom$ is $(n-1)$-Ahlfors regular. Then the following are equivalent:
\begin{itemize}
\item[\emph{(i)}] $\om^\pm$ are both vanishing chord-arc domains with $\nu\in VMO_{\loc}(\sigma)$ (see Definition \ref{d:CADSC}).
\item[\emph{(ii)}] There are $X^\pm \in \om^\pm$  such that $k^\pm = \frac{d\hm_\pm^{X^\pm}}{d\sigma}$ exist and $\log k^\pm \in VMO_{\loc}(d\sigma)$.
\end{itemize}
\end{theorem}
Furthermore we obtain corresponding quantitative results, see Theorems \ref{t:osckcontrolsoscn} and \ref{FBthrm.thrm}.

 \begin{remark}
 There is some redundancy in condition (i) of Theorem \ref{Tchar2sideCAD.thrm}, that we leave in for the sake of emphasis. In particular, under the conditions of Theorem \ref{Tchar2sideCAD.thrm}, $\Omega^+$ is a vanishing chord-arc domain if and only if $\Omega^-$ is. Additionally, it is a consequence of our work in this paper (see Corollary \ref{c:VMO}) that $\nu \in VMO_{\loc}(\sigma)$ is equivalent to (under the hypothesis of Theorem \ref{Tchar2sideCAD.thrm})  $\Omega^{\pm}$ being vanishing chord-arc domains.
\end{remark} 

In this paper techniques from potential theory and geometric measure theory come together allowing one to deduce geometric properties of domains.  In Section \ref{prelim} basic definitions from both areas are presented. In Section \ref{s:excessdecay} we apply classical tools of geometric measure theory dating back to De Giorgi's original work on sets of locally finite perimeter. See \cite{maggi2012sets} for references and 
an approach motivating the one presented here. The novelty is that we extend these tools from \textit{perimeter minimizers} to sets of locally finite perimeter with Ahlfors regular boundaries\footnote{Rather a representative whose boundary agrees with the support of the Gauss-Green measure. See \eqref{e:12} and Remark \ref{rmk:chooseagoodrep}.} \footnote{ The class of perimeter minimizers is a subclass of the sets we consider in Section \ref{s:excessdecay}, as defined in \eqref{e:12}. See \cite[Section 16.2]{maggi2012sets}.}, which allows us to reduce and better understand topological hypothesis
from previous works concerning potential theory in ``rough" domains ( cf. \cite{kenigtorotwophase, bortzengelstein} and the discussion in the last paragraph of Section \ref{prelim}).
The general approach we take is analogous to \cite{maggi2012sets}, but new ideas are also implemented in various places to extend the proof to a more general class of sets.
 In particular, Corollary \ref{c:BMO}, which is analogous to a well known result that plays a fundamental role in the proof of regularity of 
perimeter minimizers, shows that control on the oscillation of the unit normal provides both local control on the flatness of the boundary as well as local separation
properties (see Definition \ref{def:Reifenberg}). In addition to the proofs of these separation properties, in Appendix \ref{s:proofofseparation} we also prove that if the unit normal has small oscillation in a ball centered on the perimeter, then a large portion of the perimeter inside a slightly smaller concentric ball is contained in the graph of a Lipschitz function. Within the smaller ball both the Lipschitz norm of this function and the symmetric difference of this graph and the boundary are quantitatively controlled by the oscillation inside the larger ball.

These results should be contrasted with those in
\cite{semmes1}, \cite{semmes2}, \cite{kenig1999free}, \cite{hofmann2010singular}, \cite{jessica1} and \cite{jessica2}. In \cite{semmes1} and \cite{semmes2}, Semmes introduced the notion 
of chord-arc surfaces with small constant. (His definition is similar to ours in Definition \ref{d:CADSC}, except that he works on $C^2$ connected embedded hypersurfaces, whereas we assume Ahlfors regularity.) He focused on characterizing such surfaces through the behavior of singular integral operators on them. One crucial tool in Semmes' study is the ``Semmes decomposition theorem" which allows one to write a large portion of the chord-arc surface with small constant as the graph of a Lipschitz function (this is in the same vein as our aforementioned Lipschitz covering in Appendix \ref{s:proofofseparation} of this paper). To obtain this decomposition, Semmes needed to assume that the surface was $C^2$ (though his estimates did not depend on the $C^2$-norm). The decomposition was later obtained in the more general context of Reifenberg flat domains by \cite{kenig1997harmonic, kenig1999free} \footnote{{We thank the referee for pointing out that no one has explicitly written the proof that a chord-arc domain with small constant in the sense of \cite{kenig1997harmonic, kenig1999free} satisfies the small $\gamma$-condition of \cite{semmes1}. Although the proof is straightforward, we include it in Remark \ref{r:fufillsemmes} to patch this gap in the literature.}} and then in the even more general context of domains with the two-sided local John condition in \cite{hofmann2010singular}. Given the Semmes decomposition one can immediately use the oscillation of the unit normal to control the Reifenberg flatness of the chord-arc surface. Our key result along these lines, Corollary \ref{BMOimpflatnessdomains.cor}, also implies that the oscillation of the unit normal controls the Reifenberg flatness of the chord-arc surface. However, our condition (which is implied by a local two-sided corkscrew condition), is weaker than two-sided local John. Furthermore, our approach does not need a Semmes-type decomposition (though, as mentioned above, it does yield such a decomposition).

In addition to his geometric study of chord-arc surfaces with small constant, Semmes expressed interest in 
obtaining potential theoretic characterizations. These characterizations were investigated by Kenig and Toro, with the {\it a priori} assumption of Reifenberg flatness in 
\cite{kenig1997harmonic}, \cite{kenig1999free} and \cite{kenig2003poisson}. As a consequence of results herein, we show that the flatness 
hypothesis is redundant\footnote{As mentioned above, to show that a domain is $\delta$-chord-arc from the oscillation of the unit normal, one can use \cite[Theorem 4.19]{hofmann2010singular} (which does not require Reifenberg flatness) instead of \cite[Theorems 4.2 and 4.4]{kenig1999free} in the presence of the two-sided local John condition. Corollary \ref{BMOimpflatnessdomains.cor} allows one to remove the two-sided John condition from \cite[Theorem 4.19]{hofmann2010singular}. Then one can state the hypotheses of some theorems in \cite{kenig1997harmonic}, \cite{kenig1999free} and \cite{kenig2003poisson}  in terms of the oscillation of the unit normal alone, that is, without assuming {\it a priori} Reifenberg flatness (or two-sided local John). See e.g. \cite[Theorem 4.21]{hofmann2010singular}. }, this in turn, allows one to remove the {\it a priori} topological assumption of Reifenberg flatness  (or, more generally, two-sided local John) from some theorems in the aforementioned works of Kenig and Toro. 
In Section \ref{s: twophaseharmonicmeasure} we focus on the local two phase free boundary problem for the Poisson kernels. 
In Section \ref{geometric-info} we show that local doubling properties of $\omega^\pm$ combined with the
Ahlfors regularity of the boundary yield the existence of corkscrew balls on both sides (locally) and therefore imply local uniform rectifiability of the boundary (see Lemma \ref{ADRandDoublingimpCS.lem} and Corollary \ref{hmdoubcorr.lem}). In Section \ref{localization-tt} we show that  in our setting, the assumption $\log k^\pm \in VMO_{\loc}(d\sigma)$ yields information about the doubling properties of $\omega^\pm$ and the local optimal behavior of $k^\pm$ (see Lemma \ref{lm:doubling}). Combining the results in 
Sections \ref{geometric-info} and \ref{localization-tt} we {\it almost} recover the hypothesis in \cite{bortz2017singular}. The proof of Theorem \ref{t:osckcontrolsoscn} follows the general scheme 
of the proof in \cite{bortz2017singular} with an additional domain approximation scheme (see Appendix \ref{URapprox.sect}) and special attention given to the constants in order to prove a quantitative result.

\section*{Acknowledgement}
We are grateful for many helpful discussions with S. Hofmann. In particular, we are indebted to him for the  proof of Lemma \ref{ADRandDoublingimpCS.lem} and the ideas in Appendix \ref{URapprox.sect}, which are adapted from \cite{HofMarMay}. Some of these results first appear in the manuscript \cite{bortzengelstein} by the first two authors. Through the preparation of a mini-course to be taught at ICMAT, the other three authors
recognized that some of the hypothesis in the original paper could be improved. We also like to thank the referee(s) for their careful reading and helpful comments.

\section{Preliminaries}\label{prelim}

In the sequel, $n$ is a natural number with $n \ge 3$. We typically use $E$ to denote a set of locally finite perimeter in $\R^{n}$,  that is, a Lebesgue measurable set such that for every compact set $K \subset \R^{n}$
$$
\sup \left\{ \int_{E} \div T(x) ~ dx : T \in C^{1}_{c}(\R^{n} ; \R^{n}), \spt T \subset K, \sup_{\R^{n}} |T| \le 1 \right\} < \infty.
$$ 
 And we use $\Omega$ to denote a domain, i.e. an open and connected set, in $\mathbb{R}^n$. Oftentimes the domain $\Omega$ will also be a set of locally finite perimeter, for example if we assume $\partial\Omega$ is upper Ahlfors regular (see \cite[Section 5.11]{evans1992measure}). We recall a few results.

\begin{proposition} \label{p:2.1} (\cite[Proposition 12.1]{maggi2012sets})
If $E$ is a Lebesgue measurable set in $\mathbb{R}^{n}$, then $E$ is a set of locally finite perimeter if and only if there exists an $\rn$-valued Radon measure $\mu_{E}$ on $\rn$ such that
$$
\int_{E} \div T(x) dx = \int_{\rn} T  \cdot d \mu_{E}, \qquad \forall T \in C^{1}_{c}(\rn;\rn).
$$
\end{proposition}
The measure $\mu_{E}$ is is called the Gauss-Green measure of $E$.

For a vector-valued Radon measure $\mu$ on $\rn$, the total variation of $\mu$ is denoted by $|\mu|$.  We recall (see \cite[Chapter 4]{maggi2012sets}) that $|\mu|$ is a non-negative Radon measure that has the following characterization on open sets $V \subset \rn$
\begin{equation} \label{e:10}
|\mu|(V) = \sup \left\{ \int_{{\rn}} T \cdot d \mu : T \in C^{1}_{c}(V; \R^{n}), | T | \le 1 \right\}.
\end{equation}

 If $E$ is a set of locally finite perimeter, and $\mu_{E}$ the associated Gauss-Green measure, recall (see \cite[Chapter 15]{maggi2012sets})
the reduced boundary of $E$, denoted $\partial^{*}E$,  is defined by
\begin{equation} \label{e:redbdry}
\partial^{*}E = \left\{ x \in \spt \mu_{E} : \lim_{r \downarrow 0} \frac{\mu_{E}(B(x,r))}{|\mu_{E}|(B(x,r))} = \nu_{E}(x) \in \St^{n-1} \right\}
\end{equation}
In fact, $\nu_{E} : \partial^{*}E \to \St^{n-1}$ defined by the limit in \eqref{e:redbdry} is a Borel function called the measure-theoretic outward pointing unit normal. Moreover, the following is a version of De Giorgi's structure theorem.
\begin{theorem}[De Giorgi's structure theorem, \cite{maggi2012sets} Theorem 15.9]
If $E \subset \rn$ is a set of locally finite perimeter, then
$$
\mu_{E} = \nu_{E} \cH^{n-1} \restr \partial^{*}E \qquad \text{and} \qquad |\mu_{E}| = \cH^{n-1} \restr \partial^{*}E.
$$
\end{theorem}

\begin{remark}\label{l:boundarycont}
For a set of locally finite perimeter $E \subset \R^{n}$ there are several notions of boundary: the reduced boundary $\partial^{*}E$, the measure-theoretic boundary $\partial_{*}E$, the support of the Gauss-Green measure, and the topological boundary (see \cite{evans1992measure} or \cite{maggi2012sets} for relevant definitions). The following relationships between different notions of the boundary hold
\begin{equation} \label{e:boundarycont}
\partial^{*}E \subset \partial_{*}E 
\subset \spt \mu_{E} \subset \partial E \quad \text{and} \quad  \mathcal{H}^{n-1}(\partial_*E \setminus \partial^* E) = 0.
\end{equation}
In particular, $\partial^{*}E = \partial E$ implies $\partial^{*}E =\partial_{*}E = \spt \mu_{E} = \partial E$.
\end{remark}

The next two propositions can be found in \cite[Propositions 4.29,4.30]{maggi2012sets}.

\begin{proposition}[Lower semi-continuity of weak$*$ convergence] \label{p:lsc}
If $\mu_{k}$ and $\mu$ are vector-valued Radon measures with $\mu_{k} \rightharpoonup \mu$, i.e. for every $\phi\in C_c(\R^n,\R^n)$
$$
\int\phi\cdot d\mu_k\rightarrow \int\phi\cdot d\mu,
$$
then for every open set $A \subset \R^{n}$ we have
\begin{equation} \label{e:lsc}
|\mu|(A) \le \liminf_{k \to \infty} |\mu_{k}|(A).
\end{equation}
\end{proposition}

\begin{proposition}
Let $\mu_{k}$ be vector-valued Radon measures on $\R^{n}$ such that $\mu_k \rightharpoonup \mu$ for some $\mu$, a vector-valued Radon measure on $\R^n$.
\begin{enumerate}
\item If additionally $|\mu_{k}| \rightharpoonup \nu$ for some $\nu$ a non-negative Radon measure on $\R^n$. Then, for every Borel set $F \subset \R^{n}$,
\begin{equation} \label{e:00}
|\mu|(F) \le \nu(F).
\end{equation}
Furthermore, if $F  \subset \rn$ is a bounded Borel set with $\nu(\partial F) = 0$, then
\begin{equation} \label{e:69}
\mu(F) = \lim_{k \to \infty} \mu_{k}(F).
\end{equation}

\item If $|\mu_{k}|(\R^{n}) \to |\mu|(\R^{n})$, and $|\mu|(\R^{n}) < \infty$, then $|\mu_{k}| \rightharpoonup |\mu|$.

\end{enumerate}
\end{proposition}

\begin{definition}[Ahlfors regularity]\label{def:AR}
A Borel measure $\mu$ on $ \R^{n}$ is said to be \emph{$d$-Ahlfors regular} if there exists a positive finite constant $C_A$ such that 
\begin{equation} \label{e:11}
C_A^{-1} r^{d} \le \mu( B(x,r)) \le C_A r^{d}
\end{equation}
  for all $x \in \spt\mu$ and all $0 < r < \diam \spt\mu$. More generally, we say that a measure $\mu$ is \emph{$d$-Ahlfors regular up to scale $r_{0}$} if \eqref{e:11} holds for all $0 < r < r_{0}$. In either case, the constant $C_{A}$ is called the \emph{Ahlfors regularity constant for $\mu$}.

Let $F\subset\R^n$ be a closed set. If \eqref{e:11} holds for the measure $\mu=\mathcal H^d\restr F$ and some $0 < d \leq n$, then $F$ is said to be ($d$-)Ahlfors regular up to scale $r_{0}$.  When $d$ is understood from context, we simply say $F$ is Ahlfors regular up to scale $r_{0}$.
\end{definition}

\begin{definition}[Uniformly Rectifiable (UR) sets]\label{defur}
Let $A \subset \R^{n}$ be a closed set that is $d$-Ahlfors regular. It is said to be uniformly rectifiable (UR) if it contains ``Big Pieces of Lipschitz Images". This means there exist a pair of constants $\theta, \Lambda > 0$ such that for all $x \in A$ and all $0 < r \le \diam(A)$ there is a Lipschitz mapping $g : B(0,r) \subset \R^{d} \to \R^{n}$ with $\lip(g) \le \Lambda$ such that $\cH^{d} \left( E \cap g(B(0,r)) \right) \ge \theta r^{d}$.
\end{definition}

One reason uniformly rectifiable sets are ubiquitous is that they are spaces on which one can develop a rich Calder\'on-Zygmund theory. An example of this, to be used (implicitly) later, is the following characterization of uniformly rectifiable sets in co-dimension $1$.

\begin{theorem}[\cite{david1991singular}, \cite{mattila1996cauchy}, and \cite{nazarov2014uniform}]
Let $F\subset \rn$ be a closed $(n-1)$-Ahlfors regular set with {the associated} measure $\sigma:= \cH^{n-1}\restr F$. Then $F$ is uniformly rectifiable if and only if the Riesz transform operator (see Definiton \ref{riesztransdef}), $\mathcal R$ is $L^2$ bounded with respect to $\sigma$, in the sense that its truncation $\mathcal{R}_\epsilon $ 
satisfies
 \begin{equation}\label{eqrtbound}
 \sup_{\eps>0} \|\mathcal{R}_\eps f\|_{L^2(F,\sigma)} \le C\rVert f \rVert_{L^2(F,\sigma)}\,\, {\forall f\in L^2(F, \sigma)},
 \end{equation}
 {with a $C > 0$ uniform in $f\in L^2(F,\sigma)$}. 
\end{theorem}

\begin{definition}[UR domain, see \cite{hofmann2010singular}]\label{URdom}
We say that an open set $\om$ is a UR domain if $\pom$ is UR,
and the measure-theoretic boundary
$\partial_*\om$ (see \cite[Chapter 5]{evans1992measure})
satisfies $\Hn(\pom \setminus  \partial_*\om)=0$.
\end{definition}

We remark that in the above definition, $\Omega$ is not required to be connected; we use the term ``UR domain'' nonetheless following the convention set by \cite[Definition 3.7]{hofmann2010singular} and also to distinguish them from UR sets (of Definition \ref{defur}).

\begin{definition}[BMO and VMO]\label{def:bmo-vmo}

Let $F \subset \R^{n}$ {be $(n-1)$-Ahlfors regular up to scale $r_{0}$}\footnote{{Of course, this notion can be defined for $d$-Ahlfors regular subsets of $\R^n$ but we are only concerned with the case $d = n-1$}}. Then, for all $0 < r < r_{0}$, {$x \in F$}, and $f \in L_{\loc}^{2}(\Hn\restr F)$, define 
\begin{equation} \label{e:tbmo}
\|f\|_{*}(x,r) = \sup_{0 < s < r} \left(\fint_{B(x,s)\cap F} \left| f(y) -\fint_{B(x,s)\cap F} f(z) \dif \Hn(z) \right|^{2} \dif \Hn(y)\right)^{\frac{1}{2}}.
\end{equation}
We say that:
\begin{enumerate}
\item $f \in$ BMO$_{\loc}(\Hn\restr {F})$ if for every compact set $K \subset \R^{n}$, there exist $R_K>0$ and $C_K>0$ such that 
\begin{equation} \label{e:bmoloc}
\sup_{0<r<R_K} \sup_{x \in {F} \cap K} \|f\|_{*}(x,r) \le C_K.
\end{equation}
\item $f \in$ BMO$_{\loc}(\Hn\restr {F})$ with constant $\kappa>0$ if for every compact set $K \subset \R^{n}$, there exists $R_K>0$ such that 
\begin{equation} \label{e:bmoloc-kappa}
\sup_{0<r<R_K} \sup_{x \in {F} \cap K} \|f\|_{*}(x,r) \le \kappa.
\end{equation}

\item $f \in$ VMO$_{\loc}(\Hn\restr {F})$ if for every compact set $K \subset \R^{n}$,
\begin{equation} \label{e:vmoloc}
\lim_{r \to 0} \sup_{x \in \partial E \cap K} \|f\|_{*}(x,r) = 0.
\end{equation}
\end{enumerate}
\end{definition}

\begin{remark}\label{comactvsballs.rmk}
It is clear that the local conditions in the definition above are equivalent to replacing arbitrary compact sets by balls centered on the boundary with radius less than, say, $(1/4){\diam(F)}$. This is obvious if ${F}$ is unbounded and if ${F}$ is bounded we can cover ${F}$ by a finite collection of such balls.
\end{remark}

\begin{definition}[Corkscrew Condition]\label{cs.def-tt}
We say an open set $E \subset \R^{n}$ satisfies the $(M,R_{0})$ interior corkscrew condition if for every $x \in \partial E$ and $r \in (0, R_{0})$ there exists a  point $x_{1}$ called the interior corkscrew point so that 
$B(x_{1}, r/M) \subset E \cap B(x,r)$. 
\end{definition}

\begin{definition}[Two-sided Corkscrew Condition]\label{tscs.def}
We say an open set $E \subset \R^{n}$ satisfies the $(M,R_{0})$ two-sided corkscrew condition if for every $x \in \partial E$ and $r \in (0, R_{0})$ there exist two  points $x_{1} \in E$ and $x_{2} \in \rn \setminus E$ such that $B(x_{1}, r/M) \subset E$ and $B(x_{2}, r/M) \subset \rn \setminus E$. We call $x_{1}$ and $x_{2}$ the interior and exterior corkscrew points respectively. 
\end{definition}

\begin{definition}[Harnack Chain Condition]
Following \cite{jerison1982boundary}, we say that a domain $\Omega$ satisfies the $(C,R)$-Harnack Chain condition if for every $0 < \rho \le R, \Lambda \ge 1$, and every pair of points $X,X^{\prime} \in \Omega$ with $\delta(X), \delta(X^{\prime}) \ge \rho$ and $|X-X^{\prime}| < \Lambda \rho$, there is a chain of balls $B_{1}, \dots, B_{N} \subset E$ with $N \le C \log_{2} \Lambda + 1$, and $X \in B_{1}, X^{\prime} \in B_{N}$, $B_{k} \cap B_{k+1} \neq \emptyset$ for all $k = 1, \dots, N-1$ and $C^{-1} \diam(B_{k}) \le \dist(B_{k}, \partial \Omega) \le C \diam (B_{k})$ for all $k = 1, \dots, N$. The chain of balls is called a ``Harnack Chain''.
\end{definition}

\begin{definition}[NTA and Chord-Arc Domain]\label{def:NTA} 
We say that $\Omega \subset \R^{n}$ is a Non-Tangentially Accessible Domain (NTA) with constants $(M,R_{0})$, if it satisfies the $(M,R_{0})$-Harnack chain condition and the $(M,R_{0})$ two-sided corkscrew condition. If $\Omega$ is unbounded, we require that $\R^{n} \setminus \partial \Omega$ consists of two, non-empty, connected components. Note that if $\Omega$ is unbounded, then $R_{0} = \infty$ is allowed.

Finally, if $\Omega$ is an NTA domain whose boundary is Ahlfors regular we say that $\Omega$ is a chord-arc domain.
\end{definition}

\begin{remark}
	Sometimes in the definition of unbounded NTA domains, it is required that $R_{0} = \infty$ (see, e.g. \cite{kenig1997harmonic}, \cite{kenigtorotwophase}). In particular, this allows one to obtain estimates on harmonic measure/functions at arbitrarily large scales. Since we are only interested in local geometric properties of $\Omega$, we allow $R_{0} < \infty$ even for unbounded domains 
$\Omega$. 

Also note that if $\om$ is an open set with an Ahlfors regular boundary and satisfies the two-sided corkscrew condition with $R_0 \approx \diam(\pom)$, then it is a UR domain (see \cite[Theorem 1]{david1990lipschitz} and also Badger \cite{badger2012null}\footnote{In fact, Badger shows that upper Ahlfors regularity is not necessary for the quantitative interior approximation by Lipschitz domains shown in \cite{david1990lipschitz}.}). In addition, {having interior and exterior corkscrews} at arbitrarily small scales forces $\partial_*\Omega = \pom$.
\end{remark}

 Let $\Sigma \subset \R^{n}$ be a closed set. For any $x \in \Sigma$ and $r>0$, we define
\begin{equation} \label{e:flatness}
\Theta(x,r) = \inf_L \left\{\frac{1}{r}D[ \Sigma \cap \overline{B(x, r)},L\cap \overline{B(x, r)}]  \right\}
\end{equation}
where the infimum is taken over all $(n-1)-$planes containing $x$. Here $D$ denotes the Hausdorff distance, that is, for non-empty sets $A,B \subset \rn$,
$D[A,B] \defeq \sup\{d(a,B): a \in A\}  + \sup\{d(b,A): b \in B\}$. With this in hand, we can define flatness as in Reifenberg \cite{reifenberg1960solution};

\begin{definition}[Reifenberg Flat and Vanishing Reifenberg Flat sets]\label{def:Reifenberg-set}
 We say a closed set $\Sigma \subset \R^{n}$ is 
$\delta-$Reifenberg
 flat for some $\delta > 0$ if for each compact set $K \subset \R^{n}$ there exists $R_K>0$ such that
\begin{equation} \label{e:rfs}
\sup_{r \in (0,R_K]}\sup_{x \in K \cap \Sigma} \Theta(x,r) < \delta.
\end{equation}
We say $\Sigma$ is a vanishing Reifenberg flat set if for every compact set $K \subset \R^{n}$
$$\lim_{r \to 0} \sup_{x \in \Sigma \cap K} \Theta(x,r) = 0.$$

\end{definition}

\begin{definition}[Reifenberg Flat and Vanishing Reifenberg Flat domains]\label{def:Reifenberg}
 Let $\delta \in (0, \delta_{n})$ where $\delta_{n}$ is chosen appropriately (see Remark  \ref{rmk:deltan}) and depends only on the dimension $n$. We say that a domain $\Omega \subset \mathbb R^{n}$ is an $\delta$-Reifenberg flat domain (or vanishing Reifenberg flat domain), if $\partial \Omega$ is $\delta$-Reifenberg flat (resp. vanishing Reifenberg flat) and $\Omega$ satisfies the {\bf separation property}{\normalfont:}  for every compact set $K \subset \rn$ there exists $R_{K} > 0$ such that for any $y \in \pom \cap K$ and $0 < r < R_K$ there exists a $\nu\in \mathbb S^{n-1}$ so that if $x \in B(y,r)$ and $\langle x - y, \nu \rangle > \delta r$ then $x \in \Omega^{c}$, and if $\langle x-y, \nu \rangle < - \delta r$ then $x \in \Omega$.

Additionally, if $\om$ is unbounded it is further required that $\R^{n} \setminus \pom$ consists of two connected components,
and that $\delta\leq \delta_n$. \footnote{Note that the definition above is slightly different from the one in \cite[Definition 1.6]{kenig2003poisson} as we do not require flatness at large scales.}
\end{definition}

\begin{definition}[Chord-arc domains with small constants and vanishing chord-arc domains] \label{d:CADSC}
Let $\delta\in (0,\delta_n)$  (where $\delta_n$ is from Definition \ref{def:Reifenberg}, see the remark below). A set of locally finite perimeter $\Omega \subset \R^{n}$ is said to be a $\delta$-chord-arc domain (or chord-arc domain with small constant) if $\Omega$ is a $\delta$-Reifenberg flat domain, $\partial \Omega$ is Ahlfors regular and for each compact set $K \subset \R^{n}$ there exists some $R > 0$ such that
\begin{equation}\label{def:oscns}
	\sup_{x \in \partial \Omega \cap K } \| \nu_{\Omega} \|_{*}(x,R) < \delta.
\end{equation}
We say a domain $\Omega$ is a chord-arc domain with vanishing constant if it is a chord-arc domain with small constant and for each compact set $K \subset \R^{n}$
\begin{equation} \label{e:vanconst}
\lim_{r \to 0} \sup_{x \in \partial \Omega \cap K} \| \nu_{\Omega} \|_{*} (x,r) = 0,
\end{equation}
that is if $\nu_{\Omega} \in VMO_{\loc}(\cH^{n-1} \restr \partial \Omega)$.
\end{definition}

\begin{remark}\label{rmk:deltan}
 We recall from \cite[Theorem 3.1]{kenig1997harmonic} that there exists a $\delta_{n}>0$ such that if $\Omega \subset \rn$ is a $\delta$-Reifenberg flat domain for some $\delta <\delta_{n}$, then $\Omega$ is (locally) an NTA domain. If $\partial \Omega$ is also assumed to be Ahlfors regular, then $\Omega$ is a chord-arc domain (as in Definition \ref{def:NTA}). This justifies the name $\delta$-chord-arc domain (or chord-arc domain with vanishing constant).
\end{remark}

The reader may wonder whether the smallness in \eqref{def:oscns} implies the smallness in \eqref{e:rfs}, for example when $\partial\Omega$ is smooth. In the planar case ($n=2$) one can show that $\sup_{x,r} \Theta(x,r) \lesssim \|\nu_{\Omega}\|_*$; but in higher dimensions this estimate holds only if we know the smallness of both parameters \emph{a priori}, otherwise $\partial\Omega$ might have small \textit{handles}. See the discussions and main theorem in \cite{semmes3}. However, when $\partial\Omega$ is assumed to be Ahlfors regular (plus some weak topological assumptions), we will show in Section \ref{s:excessdecay} how to bound $\Theta(x,r)$ by $\|\nu_{\Omega}\|_*$.

{
\begin{remark}\label{r:fufillsemmes} 
Here we record a straightforward argument that chord arc domains with small constant in the sense of \cite{kenig1997harmonic, kenig1999free} satisfy the quantitative conditions in the definition of a chord-arc surface from \cite{semmes1, semmes2}. 

By Definitions \ref{def:Reifenberg-set} and \ref{def:Reifenberg}, an $\eta$-Reifenberg flat domain $E$ satisfies the following flatness condition: for any $q \in \partial E, r > 0$ there exists some unit vector $n_{q,r}$ so that 
\begin{equation} \label{e:gammacondition}
| \langle n_{q,r}, y-q \rangle| \le  \eta r \qquad \forall y \in \overline{B(q,r)} \cap \partial E.
\end{equation}

In \cite{semmes1,semmes2} it was assumed not only that chord-arc surfaces with small constant had small BMO norm, but also that they satisfied a flatness condition like \eqref{e:gammacondition} where $n_{q,r}$ is replaced with the specific vector $\nu_{q,r} = \fint_{B(q,r) \cap \partial^{*}E} \nu_{E} d \cH^{n-1}$. This height bound is an unsurprising consequence of being both Reifenberg flat and having a small BMO norm\footnote{Being Reifenberg flat is not necessary, \emph{a priori}, as seen by Corollary \ref{c:BMO}.}. More precisely,

{\bf Claim: }
If $E$ has Ahlfors regular boundary, is an $\eta$-Reifenberg flat domain, and satisfies $\| \nu_{E}\|_{*}(q,r) \le \delta$ for some $\delta \le 1/2$, then
\begin{equation} \label{e:11gammacond}
| \langle y-q, \nu_{q,r} \rangle | \le C r \sqrt{\eta + \delta} \qquad \forall y \in B(q,r) \cap \partial E.
\end{equation}

\begin{proof}
Let $\sigma = \cH^{n-1} \restr \partial^{*} E$ and $n_{q,r}$ be the direction from the $\eta$-Reifenberg flat condition. We claim it suffices to show \eqref{e:.1} - \eqref{e:.35},
\begin{equation} \label{e:.1}
\left| \int_{B(q,r)} n_{q,r} \cdot \nu_{E} d \sigma - \omega_{n-1} r^{n-1} \right| \le C r^{n-1} \eta,
\end{equation}
\begin{equation} \label{e:.2}
\left| |\nu_{q,r}| -1 \right| \leq \delta
\end{equation}
\begin{equation} \label{e:.35}
(1 - n \eta^{2}) \omega_{n-1} r^{n-1} \le \sigma(B(q,r)) \le (1 + 2 \delta) \omega_{n-1} r^{n-1}.
\end{equation}

Indeed, \eqref{e:.1} and \eqref{e:.35} together ensure
$$
| 1 - \nu_{q,r} \cdot n_{q,r} | = \left| 1- \fint_{B(q,r)} n_{q,r} \cdot \nu_E d\sigma \right| \leq C(\eta+ \delta).
$$
Combining \eqref{e:.2} with the preceding inequality we deduce
\begin{align*}
\left| n_{q,r} \cdot \frac{\nu_{q,r}}{|\nu_{q,r}|} - 1 \right| \le \left| n_{q,r} \cdot \frac{\nu_{q,r}}{|\nu_{q,r}|} - \frac{1}{|\nu_{q,r}|}\right| + \left| \frac{1}{|\nu_{q,r}|} - 1 \right| \le \frac{C}{1-\delta} \left[\eta + \delta \right] = C(\eta + \delta).
\end{align*} 
This in turn implies
\begin{equation} \label{e:11close}
\left|n_{q,r} - \frac{\nu_{q,r}}{|\nu_{q,r}|} \right|^{2} \le C ( \eta + \delta).
\end{equation}

Consequently, for $y \in B(q,r) \cap \partial E$,
\begin{align*}
| \langle y-q, \nu_{q,r} \rangle | &\le (1 + \delta) \left| \left \langle y-q, \frac{\nu_{q,r}}{|\nu_{q,r}|} \right \rangle \right| \\
& \le (1 + \delta) \left\{  \left| \left \langle y-q, n_{q,r} \right \rangle \right| +  \left| \left \langle y-q, \frac{\nu_{q,r}}{|\nu_{q,r}| }- n_{q,r}\right \rangle \right|  \right\} \\
& \le (1 + \delta) \left\{ \eta r + C r \sqrt{\eta + \delta} \right\},
\end{align*}
where the first inequality used \eqref{e:.2} and the final inequality follows from \eqref{e:11close} and the fact that $\partial E$ is $\eta$-Reifenberg flat. Since $\delta$ is small, this verifies \eqref{e:11gammacond}. Hence, it remains to check \eqref{e:.1} - \eqref{e:.35}.

We compare $E \cap B(q,r)$ to $B(q,r) \cap \{ \langle y-q, n_{q,r}\rangle \le 0 \}$ to verify \eqref{e:.1}. Indeed, for any constant vector $e$ it follows
\begin{equation} \label{e:11div0}
0 = \int_{B(q,r)\cap E} \div e = \int_{B(q,r)\cap \partial E} e \cdot \nu_{E} d \sigma + \int_{\partial B(q,r) \cap E} e \cdot \frac{y-q}{|y-q|} d \cH^{n-1}.
\end{equation}
Plugging in $e = n_{q,r}$ we get 
\begin{align}
\label{e:11.1} & \int_{\partial B(q,r) \cap E} n_{q,r}  \cdot \frac{y-q}{|y-q|} d \cH^{n-1}= -\int_{B(q,r) \cap \partial E} n_{q,r} \cdot \nu_E d\sigma
\end{align}
Since $n_{q,r}$ comes from the $\eta$-Reifenberg flat condition,
$$
\begin{cases}
E \cap \overline{B(q,r)} \subset \{ \langle y - q, n_{q,r} \rangle \leq \eta r \} \cap \overline{B(q,r)}  \\
\{ \langle y  - q, \nu_{q,r} \rangle \leq - \eta r \} \cap \overline{B(q,r)}  \subset E \cap \overline{B(q,r)}.
\end{cases}
$$
Using the divergence theorem as it was used in \eqref{e:11div0} it follows
\begin{equation} \label{e:11.2}
\int_{\partial B(q,r) \cap \{ \langle y-q, n_{q,r} \rangle \le 0 \}} n_{q,r} \cdot \frac{y-q}{|y-q|} d\mathcal{H}^{n-1} = - \int_{B(q,r) \cap \{ \langle y-q, n_{q,r} \rangle = 0 \}} n_{q,r} \cdot n_{q,r} ~d\sigma = - \omega_{n-1} r^{n-1}.
\end{equation} 
and a very generous estimate ensures
\begin{equation} \label{e:11.3}
\left| \int_{\partial B(q,r) \cap E} n_{q,r} \cdot \frac{y-q}{|y-q|} d \cH^{n-1} - \int_{\partial B(q,r) \cap \{ \langle y-q, n_{q,r} \rangle \le 0 \}} n_{q,r} \cdot \frac{y-q}{|y-q|} d \cH^{n-1} \right| \le C r^{n-1} \eta.
\end{equation}
Combining \eqref{e:11.1} - \eqref{e:11.3} confirms \eqref{e:.1}. Equation \eqref{e:.2} follows from $\|\nu\|_{*}(B(q,r)) \le \delta$. Details are included when the same statement is verified in \eqref{e:44}. 

It only remains to show \eqref{e:.35}. The lower bound follows immediately from $\partial E$ being $\eta$-Reifenberg flat and the separation property since then
$$
\sigma(B(q,r) \ge \omega_{n-1} \left(r \sqrt{1 - \eta^{2}} \right)^{n-1} \ge \omega_{n-1}\left(1 - \frac{n-1}{2} \eta^{2} \right) r^{n-1}.
$$
For the upper bound, we use \eqref{e:11div0} with $e = \nu_{q,r}$ to obtain the estimate
\begin{align*} \label{e:11.4}
\sigma(B(q,r)) |\nu_{q,r}|^{2} = \left| \int_{B(q,r)} \nu_{q,r} \cdot \nu_{E} d \sigma \right| & = \left| \int_{\partial B(q,r) \cap E} \nu_{q,r} \cdot \frac{y-q}{|y-q|} d \cH^{n-1} \right| \\
	& \leq |\nu_{q,r}| \omega_{n-1} r^{n-1}.
\end{align*}
Therefore by \eqref{e:.2}
\[ \sigma(B(q,r)) \leq \frac{1}{|\nu_{q,r}|} \omega_{n-1} r^{n-1} \leq (1+2\delta) \omega_{n-1} r^{n-1}.  \]
\end{proof}

\end{remark}
}

\section{Flatness from Control on Oscillation}\label{s:excessdecay}

In this section we introduce a class of well-behaved sets $\cA(C_{A},r_0)$, and prove our key geometric result, Corollary \ref{c:BMO}.
Namely in the class, $\cA(C_{A},r_0)$, the oscillation of the unit normal controls the flatness (in the sense of Reifenberg) of the boundary. One key tool is the ``excess" of a set of locally finite perimeter, first introduced by De Giorgi in \cite{DG61} and ubiquitous in the calculus of variations. Due to Lemma \ref{l:excessvsoscillation}, all of our arguments could also be written in terms of the mean oscillation of the unit normal. Given $r_0 \in (0, \infty)$ and $C_A \in [1,\infty)$, we define a class of sets

\begin{equation} \label{e:12}
\cA(C_{A}, r_{0}) = \left\{ E \subset \R^{n} \bigg|  \substack{ E \text{ is a set of locally finite perimeter satisfying $\partial E = \spt \mu_{E}$ and its perimeter} \\ \text{ measure } |\mu_{E}| \text{ is $(n-1)$-Ahlfors regular up to scale $r_{0}$ with constant $C_{A}$} } \right\}.
\end{equation}

Uniformly rectifiable domains (up to choosing a representative from the equivalence class, see Remark \ref{rmk:chooseagoodrep}) with Ahlfors regularity constant $C_{A}$ form a subset of $\cA(C_{A},r_{0})$. {The complement of a quasiminimal surface of codimension $1$ is the disjoint union of two open domains of $\R^n$ (see \cite{david1998quasiminimal}), and each of these domains would fall within the class $\cA(C_{A},r_{0})$.}

\begin{remark}\label{rmk:chooseagoodrep}
The condition that $\partial E = \spt \mu_{E}$ corresponds to choosing a representative for our set amongst the \textit{equivalence class} of sets of locally finite perimeter (see \cite[Proposition 12.19, Remark 16.11]{maggi2012sets}): for any set of finite perimeter $E$, we can find a Borel set $F$ such that
\[ |E\Delta F|=0, \quad \partial F = \spt\mu_F = \spt\mu_E. \] 
This choice is necessary since we want to deduce information on the topological boundary from information on the measure-theoretic unit outer normal, which is merely defined on the reduced boundary $\partial^* E$, see for example Lemma \ref{l:SL} and Theorem \ref{t:SL}.
\end{remark}

A particularly useful property of $\cA(C_{A},r_{0})$ is that if $E \in \cA(C_{A},r_{0})$ then $\R^{n} \setminus E \in \cA(C_{A},r_{0})$. This follows since $\mu_{E} = - \mu_{\R^{n} \setminus E}$ and $\partial E = \partial (\R^{n} \setminus E)$.

\begin{remark}\label{r:variationAR}
If $E \in \cA(C_{A},r_{0})$ then $\partial E$ is $(n-1)$-Alhfors regular since $\partial E = \spt \mu_{E}$ and $\cH^{n-1}(\partial E \setminus \partial^{*} E) = 0$ (see Theorem 6.9 \cite{Mattila}). Thus
$$
|\mu_{E}| =  \cH^{n-1} \restr \partial^{\ast} E=\cH^{n-1} \restr \partial E.
$$

\end{remark}

\begin{definition}[Cylinders and excess: {c.f. \cite[Chapter 22]{maggi2012sets}}]
For $r > 0, x \in \R^{n}$, and some $\nu \in \St^{n-1}$, we let
\begin{equation} \label{e:cylinder}
C(x,r,\nu) = \{ y : |\langle x - y , \nu \rangle | < r, | x-y - \langle x-y, \nu \rangle \nu | < r \}.
\end{equation}
Note that $C(x,r,\nu)$ is a cylinder with center $x$, radius and height $r$, and axial direction $\nu$. For a set of locally finite perimeter $E$, $x \in \partial E$, $r > 0$, and $\nu \in \St^{n-1}$ we define the \emph{cylindrical excess}
\begin{equation} \label{e:cylex}
e(E, x , r , \nu) = \frac{1}{r^{n-1}} \int_{C(x,r,\nu) \cap \partial^{*} E} \frac{ |\nu_{E} - \nu|^{2}}{2} \dif \cH^{n-1} 
\end{equation}
\end{definition}

The following lemma elucidates the relationship between oscillation of the unit normal and excess. 

\begin{lemma}\label{l:excessvsoscillation}
Let $E \in \cA(C_{A},r_0)$ and let $Q \in \partial E$ and $0 < r < r_0$. There exists some constant $0 < C < \infty $ (which depends only on $C_A$ and the dimension) such that \begin{equation}\label{e:excessvsoscillation1}\fint_{B(Q,r)\cap \partial^* E}\left|\nu_E -(\nu_{E})_{Q,r}\right|^2 \dif \cH^{n-1} \leq Ce(E, Q, r, \nu) \end{equation} for any $\nu \in \St^{n-1}$, where $(\nu_E)_{Q,r}$ represents the integral average of $\nu_E$ with respect to $\cH^{n-1}$ on $B(Q,r)\cap \partial^* E$. 
Furthermore, as long as $|(\nu_E)_{Q,r}| \neq 0$, then

\begin{equation}\label{e:excessvsoscillation2} e\left(E, Q, \frac{r}{\sqrt{2}}, \frac{ (\nu_{E})_{Q,r}}{|(\nu_{E})_{Q,r}|}\right) \leq C\fint_{B(Q,r)\cap \partial^* E}\left|\nu_E -(\nu_E)_{Q,r}\right|^2 \dif \cH^{n-1}. \end{equation}
  \end{lemma}

\begin{proof}
We first prove \eqref{e:excessvsoscillation1}. Note that $B(Q, r)\cap \partial^* E \subset C(Q, r, \nu)\cap \partial^* E$ for any $Q \in \partial E$ and $r > 0$. Thus 
$$e(E, Q, r, \nu) \geq c\fint_{B(Q,r)\cap \partial^* E} \frac{|\nu_E - \nu|^2}{2}  \dif \cH^{n-1} ,$$ 
where $c$ is a constant that depends only on the Ahlfors regularity of $E$. We can compute 
\begin{equation}\label{e:excessoverosc} \begin{aligned}\fint_{B(Q,r)\cap \partial^* E} & |\nu_E(x) - (\nu_E)_{Q,r}|^2 \dif  \cH^{n-1} \\
& \leq 2\fint_{B(Q,r)\cap \partial^* E} |\nu_E(x) - \nu|^2 \dif  \cH^{n-1} + 2\fint_{B(Q,r)\cap \partial^* E} |\nu - (\nu_E)_{Q,r}|^2\dif  \cH^{n-1}\\
& \leq 4\fint_{B(Q,r)\cap \partial^* E} |\nu_E(x) - \nu|^2 \dif  \cH^{n-1} \leq Ce(E,Q,r,\nu),\end{aligned}
\end{equation} 
where the second inequality above follows from the triangle inequality and Jensen's inequality. This is exactly \eqref{e:excessvsoscillation1}.

To prove \eqref{e:excessvsoscillation2} it suffices to consider
$$
\fint_{B(Q,r)\cap \partial^* E} | \nu_{E} - (\nu_{E})_{Q,r}|^{2} \dif \cH^{n-1} = \epsilon < 1.
$$
We first estimate $|(\nu_E)_{Q,r}|$; note,
\begin{align} \label{e:44}
(|1 - |(\nu_{E})_{Q,r}|)^{2} &=\fint_{B(Q,r)\cap \partial^* E} (|\nu_{E}| - | (\nu_{E})_{Q,r}|)^{2} \dif\cH^{n-1} \\
\nonumber & \le\fint_{B(Q,r)\cap \partial^* E} |\nu_{E} - (\nu_{E})_{Q,r}|^{2} \dif\cH^{n-1} = \epsilon
\end{align}
and
\begin{equation} \label{e:45}
| (\nu_{E})_{Q,r}| = \left|\fint_{B(Q,r)\cap \partial^* E} \nu_{E} \dif\cH^{n-1}  \right| \le\fint_{B(Q,r)\cap \partial^* E} | \nu_{E}| \dif\cH^{n-1} = 1.
\end{equation}
Combining \eqref{e:44} with \eqref{e:45} ensures that $1 - \sqrt{\epsilon} \le |(\nu_{E})_{Q,r}| \le 1$. Let $\nu_0 \equiv \frac{(\nu_{E})_{Q,r}}{|(\nu_{E})_{Q,r}|}$ and compute,
\begin{align*}
| \nu_{E} - \nu_{0}| &\le |\nu_{E} - (\nu_{E})_{Q,r}| + |(\nu_{E})_{Q,r}| \left| 1 - \frac{ 1}{|(\nu_{E})_{Q,r}|} \right| \\
& \le | \nu_{E} - (\nu_{E})_{Q,r}| + |1 - (\nu_{E})_{Q,r}| \\
& \le | \nu_{E} - (\nu_{E})_{Q,r}| + \epsilon^{1/2},
\end{align*}
so that 
\begin{equation} \label{e:46}
|\nu_{E} - \nu_{0}|^{2} \le 2|\nu_{E} - (\nu_{E})_{Q,r}|^{2} + 2 \epsilon.
\end{equation}
Notably, \eqref{e:46} and $C(Q,\frac{r}{\sqrt{2}},\nu_{0}) \subset B(Q,r)$ imply
\begin{align*}
&e\left(E,Q,\frac{r}{\sqrt{2}},\nu_{0}\right)  = \frac{2^{(n-1)/2}}{r^{n-1}} \int_{C(Q,\frac{r}{\sqrt{2}},\nu_{0})\cap \partial^* E} \frac{|\nu_{E}-\nu_{0}|^{2}}{2} \dif\cH^{n-1} \\
 & \quad \le \frac{2^{(n-1)/2}}{r^{n-1}}\int_{C(Q,\frac{r}{\sqrt{2}},\nu_{0})\cap \partial^* E} |\nu_{E} - (\nu_{E})_{Q,r}|^{2} \dif\cH^{n-1} \\
 &\qquad\qquad + \frac{2^{(n-1)/2}\cH^{n-1} (C(Q,\frac{r}{\sqrt{2}},\nu_{0})\cap\partial^\ast E)}{r^{n-1}} \epsilon \\
& \quad \le 2^{(n-1)/2}\frac{\cH^{n-1} (B(Q,r)\cap\partial^\ast E)}{r^{n-1}} \left(\fint_{B(Q,r)\cap \partial^* E} |\nu_{E} - (\nu_{E})_{Q,r}|^{2} \dif\cH^{n-1}  +  \epsilon\right) \\
& \quad = 2^{(n+1)/2} \frac{\cH^{n-1} (B(x,r)\cap \partial^\ast E)}{r^{n-1}} \epsilon \le C_n\cdot C_{A} \epsilon.
\end{align*}
\end{proof}

\begin{remark}\label{r:excessunderoperations}
We recall some basic properties of the cylindrical excess (see for instance \cite[Chapter 22]{maggi2012sets} for more details). The cylindrical excess is invariant under translation and scaling in the sense that if $E_{x,r} = \frac{ E - x}{r}$, then
\begin{equation} \label{e:scalex}
e(E_{x,r},0,1, \nu) = e(E,x,r, \nu).
\end{equation}

Furthermore, if $r < s$, the non-negativity of the integrand ensures
\begin{equation*}
\frac{1}{r^{n-1}} \int_{C(x,r,\nu) \cap \partial^{*} E} \frac{ |\nu_{E} - \nu|^{2}}{2} \dif \cH^{n-1} \le \left( \frac{s}{r} \right)^{n-1} \frac{1}{s^{n-1}} \int_{C(x,s,\nu) \cap \partial^{*} E} \frac{ |\nu_{E} - \nu|^{2}}{2} \dif \cH^{n-1},
\end{equation*}
that is,
\begin{equation} \label{e:compscales}
e(E, x , r , \nu) \le \left( \frac{ s}{r} \right)^{n-1} e(E,x,s,\nu).
\end{equation}

Finally, since $\nu, \nu_{E}$ are each of unit length, $
\frac{|\nu_{E} - \nu|^{2}}{2} = 1 - \langle \nu_{E}, \nu \rangle
$
so that
\begin{equation} \label{e:cylex2}
e(E, x , r , \nu) = \frac{1}{r^{n-1}} \int_{C(x,r,\nu) \cap \partial^{*} E} \left( 1 - \langle \nu_{E} , \nu \rangle \right) \dif \cH^{n-1}.
\end{equation}
\end{remark}

Given a sequence of sets of locally finite perimeter $\{E_{k}\}_{k \in \N}$ in $\rn$, we say that $\{E_{k}\}$ converges to $E$ in $L^{1}_{\loc}(\rn)$ and write $E_{k} \xrightarrow{L^{1}_{\loc}(\rn)} E$ if $\lim_{k \to \infty} \cH^{n} \left( E \Delta E_{k} \right) = 0$.
The following compactness theorem is the key tool used in proving the flatness result. 
\begin{theorem} \label{t:cacomp}
If $\{E_{k}\}_{k \in \N} \subset \cA(C_{A}, r_{0})$ with $0 \in \partial E_{k}$ for all $k \ge 1$, there exist a subsequence $\{E_{k_{j}}\}_{j \in \N}$, a set $E$ of locally finite perimeter, and a non-negative Radon measure, $\mu$ such that as $j$ approaches infinity,
\begin{equation} \label{e:13}
E_{k_{j}} \xrightarrow{L^{1}_{\loc}(\R^{n})} E, \quad \mu_{E_{k_{j}}} \rightharpoonup \mu_{E} \quad  \text{and} \quad |\mu_{E_{k_{j}}}| \rightharpoonup \mu.
\end{equation}
Additionally, $\partial E = \spt \mu_E$ and $\mu$ is $(n-1)$-Ahlfors regular up to scale $r_{0}$ with constant $C_{A}$. Furthermore, $|\mu_{E}| \le \mu$ and
\begin{enumerate}
\item If $x \in \partial E$, then for all $j \in \N$ there exist $x_{k_{j}} \in \partial E_{k_{j}}$ such that $\lim_{j \to \infty} x_{k_{j}} = x$.
\item If $x \in \spt \mu$, then for all $j \in \N$ there exist $x_{k_{j}} \in \partial E_{k_j}$ so that $\lim_{j \to \infty} x_{k_{j}} = x$. 
\item If for all $j \in \N$, $x_{k_{j}} \in \partial E_{k_{j}}$ and $\lim_{j \to \infty} x_{k_{j}} =  x$ then $x \in \spt \mu$.
\end{enumerate}
\end{theorem}
\begin{remark}
\begin{itemize}
	\item We note that (2) and (3) in Theorem \ref{t:cacomp} combine to say that $x \in \spt \mu$ if and only if there exists $x_{k_{j}} \in \partial E_{k_{j}}$ such that $x_{k_{j}} \to x$. However, without additional hypotheses, all that is known is that
$$
\spt \mu_{E} \subseteq \spt \mu.
$$
	\item  Unlike in the analogous theorem \cite[Theorem 21.14]{maggi2012sets} for perimeter minimizers, here in general we do not have $\mu=|\mu_E|$ because of possible cancellations for sets of finite perimeter. However, with further information on the excess, we will be able to conclude $\mu = |\mu_E|$, see for example Lemma \ref{l:SL}.
\end{itemize}
\end{remark}
\begin{proof}
 Standard techniques and a diagonalization argument (see for instance \cite[Sections 12.4, 21.5]{maggi2012sets}) verify that sets whose boundary are uniformly Ahlfors-regular (i.e. Ahlfors regular with constants independent of the element in the sequence) are pre-compact in the space of sets of locally finite perimeter. That is to say, there exists some set of locally finite perimeter $E \subset \rn$ so that $\chi_{E_{k_j}} \rightarrow \chi_E$ in $L^1_{\loc}$ and $\mu_{E_{k_j}} \rightharpoonup \mu_E$ in a weak star sense. Without loss of generality (see Remark \ref{rmk:chooseagoodrep}) we may assume that $\mathrm{spt}\, \mu_E = \partial E$. Finally, note the  $|\mu_{E_{k_j}}|$ are uniformly Ahlfors regular  (see Remark \ref{r:variationAR}) and hence precompact. Without explicitly relabeling the new subsequence, there exists a Radon measure $\mu$ on $\rn$ so that  $|\mu_{E_{k_j}}| \rightharpoonup \mu$ in the weak star sense. Thus \eqref{e:13} holds. 

The fact that $|\mu_E| \leq \mu$ follows from \eqref{e:00}. This ensures that $\spt \mu_{E} \subset \spt \mu$, so (2) which is a standard fact implies (1). Moreover (2) and the uniform upper regularity of $\{|\mu_{E_{k_j}}|\}$ imply the upper Ahlfors regularity of $\mu$. 

We show (3) and lower Ahlfors regularity of $\mu$ simultaneously. For each $j \in \N$ suppose $x_{k_{j}} \in \partial E_{k_{j}} = \spt |\mu_{E_{k_j}}|$ such that $x_{k_{j}} \to x$. 

 Fix $0 < s < r_{0}$ and fix $\epsilon \in (0,1)$. Note that for $k_{j}$ large enough, 
$B(x_{k_{j}}, s(1 - \epsilon)) \subset B(x, s(1 - \epsilon/2))$.  Since $E_{k_{j}} \in \cA(C_{A},r_{0})$ it follows that 
$$
C_{A}^{-1} (s(1 - \epsilon))^{n-1} \le |\mu_{E_{k_{j}}}|(B(x_{k_{j}}, s (1- \epsilon))) \le | \mu_{E_{k_{j}}}| \left( \overline{B(x,s (1- \epsilon/2)) }\right) 
$$
so that by weak$*$ convergence of $|\mu_{E_{k_{j}}}|$ to $\mu$
$$
C_{A}^{-1} (s (1- \epsilon))^{n-1} \le \limsup_{j} | \mu_{E_{k_{j}}}|  \left( \overline{B(x,s(1 - \epsilon/2)) }\right)  \le \mu \left( \overline{B(x, s(1 - \epsilon/2))} \right),
$$
taking $\epsilon \to 0$ results in $C_{A}^{-1} s^{n-1} \le \mu(B(x,s))$ for all $s \in (0, r_{0})$; in particular $x \in \spt \mu$, verifying (3). On the other hand, since (2) and (3) combine to show that $x \in \spt \mu$ if and only if there exists $x_{k_{j}} \in \partial E_{k_{j}}$ such that $x_{k_{j}} \to x$, this demonstrates that $\mu$ is $(n-1)$-lower Ahlfors regular up to scale $r_{0}$ with constant $C_{A}$. 
\end{proof}

We now prove that small excess implies local measure theoretic separation.  To simplify notation, define $e_n(E, x, r) = e(E, x, r, e_n)$.

\begin{lemma}[Separation Lemma ({compare with \cite[Lemma 22.10]{maggi2012sets})}] \label{l:SL}
Given $C_{A} \ge 1$, $t_{0} \in (0,1)$, there exists $\omega(n,t_{0},C_{A}) \in (0, \infty)$ such that if $E \in \cA(C_{A}, 2r)$ for some $r > 0$ and if there exist $x_{0} \in \partial E$ and $\nu \in \St^{n-1}$ with 
$$
e(E, x_{0},2 r, \nu) \le \omega(n, t_{0}, C_{A}),
$$
then
\begin{equation} \label{e:15}
| \langle x - x_{0}, \nu \rangle| < t_{0} r \qquad \forall x \in C(x_{0}, r, \nu) \cap \partial E,
\end{equation}
\begin{equation} \label{e:16}
\left| \{ x \in C(x_{0}, r, \nu) \cap E \mid \langle x - x_{0}, \nu \rangle > t_{0} r \} \right| = 0,
\end{equation}
and
\begin{equation} \label{e:17}
\left| \{ x \in C(x_{0}, r, \nu) \cap E^{c} \mid \langle x - x_{0} , \nu \rangle < - t_{0} r \} \right| = 0.
\end{equation}
(Note, here and below for any  Lebesgue measurable set $\mathcal O \subset \mathbb R^n$ we write $|\mathcal O|$ to denote the Lesbesgue measure of $\mathcal O$). 
\end{lemma}

\begin{proof}
The proof follows by a compactness-contradiction argument. If Lemma \ref{l:SL} does not hold, there exist $C_A>1$, $t_{0} \in (0,1)$, a sequence of sets $\{F_{k} \}_{k \in \N}$ and radii, $r_k > 0$, such that $F_k \in \cA(C_{A},2r_{k})$, a sequence of points $x_{k} \in \partial F_{k}$, and a sequence of directions $\nu_{k} \in \St^{n-1}$, with
$$
e(F_{k}, x_{k}, 2r_{k}, \nu_{k}) \le 2^{-k},
$$
such that at least one of the following conditions holds for infinitely many $k$:

\begin{equation} \label{e15}
\{x \in C(x_{k}, r_{k}, \nu_{k}) \cap \partial F_{k} \mid  |q_{k}(x)|  > t_{0} r_{k} \} \neq \emptyset,
\end{equation}
\begin{equation} \label{e16}
\left| \{ x \in C(x_{k}, r_{k}, \nu_{k}) \cap F_{k} \mid  q_{k}(x) > t_{0}r_{k} \} \right| > 0,
\end{equation}
or
\begin{equation} \label{e17}
\left| \{ x \in C(x_{k}, r_{k}, \nu_{k}) \cap  F_{k}^{c} \mid  q_{k}(x) < - t_{0} r_{k} \} \right| > 0,
\end{equation}
where $q_{k}(x) = \langle x - x_{k}, \nu_{k} \rangle$. 

By rescaling, recentering, and rotating (see Remark \ref{r:excessunderoperations}) we may assume that $\nu_k \equiv e_n, x_k \equiv 0$ and $r_k \equiv 1$. Note that the transformed domains are now in $\cA(C_A,2)$. Abusing notation we call these new sets $F_k$. Note that, 
\begin{equation} \label{e:fksmall}
e_{n}(F_{k},0,2) \le 2^{-k} \qquad \forall k \ge 1.
\end{equation}

Writing $C_{r} = C(0,r,e_{n})$ and $q(x) = \langle x, e_{n} \rangle$ we rewrite \eqref{e15} - \eqref{e17} as,
\begin{equation} \label{15}
\{ x \in C_{1} \cap \partial F_{k} \mid t_{0} \le |q(x)| \} \neq \emptyset,
\end{equation}
\begin{equation} \label{16}
| \{  x\in C_{1} \cap F_{k} \mid q(x) > t_{0} \}| >  0,
\end{equation}
or
\begin{equation} \label{17}
| \{ x \in C_{1} \setminus F_{k} \mid q(x) < - t_{0} \} | > 0.
\end{equation}

By Theorem \ref{t:cacomp}, there exists a set of finite perimeter $F \subset C_{5/3}$ with $0 \in \partial F = \spt | \mu_{F}|$ and a Radon measure $\mu$ such that, by passing to a subsequence which we do not explicitly relabel, $F_{k} \cap C_{5/3} \to F$ in $L^{1}(\R^{n})$, $\mu_{F_{k}\cap C_{5/3} }\rightharpoonup \mu_{F}$, and $|\mu_{F_{k}\cap C_{5/3}}| \rightharpoonup \mu$ with $|\mu_{F}| \le \mu$.

Consider an open set $U$ such that $\overline U \subset C_{5/3}$. Then \eqref{e:cylex2} implies
\begin{align}\label{tt-1}
\left(\frac{5}{3}\right)^{n-1}e_{n}(F_{k}, 0, 5/3) & \ge \int_{U \cap \partial^{*}F_{k}} (1 - e_{n} \cdot \nu_{F_{k}} ) \dif \cH^{n-1}  = |\mu_{F_{k}}|(U) - e_{n} \cdot \mu_{F_{k}}(U)  \ge 0,
\end{align}
where the final inequality follows since  
\begin{equation} \label{e:unit}
d \mu_{F_{k}} = \nu_{F_{k}} |d \mu_{F_{k}}| \quad \text{and} \quad |\nu_{F_{k}}| = 1 \quad |\mu_{F_{k}}| -a.e. 
\end{equation}

Then  \eqref{e:compscales} and \eqref{e:fksmall} ensure that as $k$ tends to infinity 
$$
e_{n}(F_{k}, 0 ,5/3) \le \left( \frac{6}{5}\right)^{n-1} e_{n}(F_{k}, 0, 2) \to 0.
$$ 
This combined with \eqref{tt-1} yields
\begin{equation}\label{tt-2}
0\le \lim_{k\to\infty} \left\{ |\mu_{F_{k}}|(U) - e_{n} \cdot \mu_{F_{k}}(U) \right\} \le C_n\lim_{k\to\infty}e_{n}(F_{k},0, 5/3)=0.
\end{equation}
 Since \eqref{e:00} says $|\mu_{F}| \le \mu$ we can apply \eqref{e:69} to both $|\mu_{F_{k}}|$ and $\mu_{F_{k}}$ to learn 
\begin{equation} \label{e:100}
\mu(U) = e_{n} \cdot \mu_{F}(U) \qquad \text{ for any open set } U\subset\subset C_{5/3}, \text{ with } \mu(\partial U) = 0.
\end{equation}

Note that by Theorem \ref{t:cacomp}, $\mu$ is Ahlfors regular with constant $C_{A}$ up to scale $2$ in the cylinder $C_{5/3}$.
Hence in particular for any $x\in C_{4/3} \cap \spt \mu$ and a.e. $r \in (0, 1/3)$,    $\ \mu(\partial B(x,r))=0$ and by \eqref{e:100} $\mu(B(x,r))=e_{n} \cdot \mu_{F}(B(x,r))$.  Consequently, for all $x \in \spt |\mu_F| \cap C_{4/3}$,
\begin{equation}\label{tt-3}
\limsup_{r\to 0}\frac{\mu(B(x,r))}{|\mu_F|(B(x,r))}= e_{n} \cdot \limsup_{r \to 0} \frac{ \mu_{F}(B(x,r))}{|\mu_{F}|(B(x,r))}\le 1,
\end{equation} 
where the final inequality uses the property \eqref{e:unit} for the set $F$. Thus in $C_{4/3}$ we have shown $\mu\le |\mu_F| \le \mu$, which implies $\mu=|\mu_{F}|=\cH^{n-1}\restr\partial^\ast F$. But then, \eqref{e:100} says $|\mu_{F}|= e_{n} \cdot \mu_{F}$ so that $\nu_F(x)=e_n$ at $\cH^{n-1}$-a.e. $x\in\partial^\ast F$. In particular, $e_n(F,0, 4/3)=0$, at which point \cite[Proposition 22.2]{maggi2012sets} asserts that $F \cap C_{4/3}$ is equivalent (in the sense of sets of locally finite perimeter) to $C_{4/3} \cap \{q(x) < 0 \}$ or $C_{4/3} \cap \{q(x) > 0\}$. Since $|\mu_{F}| = e_{n} \cdot \mu_{F}$ it follows that $F \cap C_{4/3}$ is equivalent to $C_{4/3} \cap \{q(x) < 0\}$. We write this as

\begin{equation} \label{e:19}
C_{4/3} \cap F \sim \{q(x) < 0 \} \cap C_{4/3}.
\end{equation}
We assumed, that one of \eqref{15} - \eqref{17} holds for infinitely many $k$. First suppose that \eqref{15} holds for infinitely many $k$. By passing to a subsequence, we may assume that \eqref{15} holds for all $k \in \N$. Then, for all $k \in \N$, there exists $x_{k} \in \partial F_{k} \cap C_{1}$ such that $t_{0}\le  |q(x_{k})|$. By passing to a subsequence, $x_{k} \to x_{\infty}$ for some $x_{\infty} \in \overline{C_{1}}$ and $|q(x_\infty)|\ge t_0$. By Theorem \ref{t:cacomp} (3), $x_{\infty} \in \spt \mu = \spt \mu_F=\partial F$. Hence (see \cite[Proposition 12.19]{maggi2012sets})
\begin{equation} \label{e:20}
0 < | B(x_\infty, s) \cap F| < \omega_{n} s^{n} \qquad \forall s > 0.
\end{equation}
However, due to $|q(x_\infty)| \ge t_0$, \eqref{e:19} implies that for any $s\le  \min \{ 1/8, |q(x_\infty)|/2\}$
\begin{equation} \label{e:21}
|B(x_\infty, s) \cap F| = \begin{cases} \omega_{n} s^{n} &\hbox{ if } q(x_\infty) < 0 \\ 0 & \hbox{ if } q(x_\infty) > 0 \end{cases}
\end{equation}
which contradicts \eqref{e:20}. This shows that \eqref{15} cannot hold for infinitely many $k$.

Arguing as above and invoking Theorem \ref{t:cacomp} (3) we conclude that there exists $k_{0} \in \N$ such that for all $k \ge k_{0}$
\begin{equation} \label{e:22}
\{ x \in C_{5/4} \cap \partial F_{k} \mid t_{0} < |q(x)| \le 1 \} = \emptyset.
\end{equation}
However, by \cite[Equation 16.7]{maggi2012sets} for all $r \in (1,5/4)$
$$
|\mu_{F_{k} \cap C_{r}}| = |\mu_{C_{r}}| \restr F_{k}^{(1)} + |\mu_{F_{k}}| \restr \left(C_{r} \cup \{ \nu_{F_{k}} = \nu_{C_{r}} \} \right).
$$
For almost every $r \in (1,5/4)$ we know $|\mu_{F_{k}}|(\partial C_{r}) = 0$ for all $k$. Then, for any such $r$ \eqref{e:22} demonstrates 
\begin{equation} \label{e:22.1}
|\mu_{F_{k} \cap C_{r}}| \left( \{ x \in C_{r} \mid t_{0} < |q(x)| < 1 \} \right) = 0  \qquad \forall k \ge k_{0}.
\end{equation}

We claim \eqref{e:22.1} implies that for almost every $r \in (1,5/4)$, $\chi_{C_{r} \cap F_{k}}$ is constant on each connected component of $\{ t_{0} < |q(x)| < 1\} \cap C_{r}$ which implies $\chi_{C_{1} \cap F_{k}}$ is constant on connected components of $\{t_{0} < |q(x)| < 1 \} \cap C_{1}$ . Indeed, choose $r \in (1,5/4)$ so that $|\mu_{F_{k}}|(\partial C_{r}) = 0$ for all $k$.  Consider the sets $U_{\pm} \defeq \{ t_0 < \pm q(x) < 1\} \cap C_{r}$, which are both open and connected. The definition \eqref{e:10} and \eqref{e:22.1} guarantee for all $k \ge k_{0}$,
$$
\int_{\rn} T\cdot d\mu_{F_k} = 0 \qquad \text{for all }T \in C^{1}_{c}(U_{\pm} ; \rn).
$$
(Because if the integral is nonzero, we can flip the sign of $T$ and get a countradition with \eqref{e:22.1}.)
Thus by Proposition \ref{p:2.1}
$$
\intrn \chi_{F_{k}} \div T dx = \int_{\rn} T\cdot d\mu_{F_k} = 0 \qquad \text{for all }T \in C^{1}_{c}(U_{\pm} ; \rn),
$$ 
that is, in the weak sense, that $\nabla \chi_{F_{k}} = 0$ on $U_{+}$ and $U_{-}$. This implies $\chi_{F_{k}}$ is almost everywhere constant on each $U_{\pm}$ (for instance, see \cite[Lemma 7.5]{maggi2012sets}).
Combining \eqref{e:19} with $\chi_{F_{k}}$ constant on each $U_{\pm}$ and  $F_{k}\cap C_{5/3} \xrightarrow{L^{1}(\R^{n})} F$, it follows that for $k \ge k_{0}$
$$
\chi_{F_{k} \cap C_{1}} = \begin{cases} 0 & \text{ for almost every } x \in C_{1} \cap \{t_{0} < q(x) < 1 \} \\ 1 & \text{ for almost every }  x \in C_{1} \cap \{-1 < q(x) < t_{0}\} \end{cases}
$$ 
This shows that \eqref{16} and \eqref{17} cannot happen for infinitely many $k$.
\end{proof}

The (qualitative) separation lemma above can be further improved to a quantative ``height bound'' of $\partial E$.  Since the proof is by fairly standard techniques in the theory of sets of locally finite perimeter, we include it Appendix \ref{s:proofofseparation} (see Theorem \ref{t:HB}). Topological considerations then imply the following theorem. 

\begin{theorem} \label{t:SL}
Given $C_{A} \ge 1$ and $n \ge 2$, there exist positive constants $C_{1} = C_{1}(n, C_{A}) < \infty$ and $\epsilon_{1} = \epsilon_1(n, C_A)$ small such that if $E \in \cA(C_{A}, 4r_{0})$ for some $ r_{0} > 0$, and $x_{0} \in \partial E$ satisfies 
\begin{equation} \label{e:CSL-2}
e(E,x_{0}, 2r, \nu) \le \epsilon_{1} 
\end{equation}
for some $\nu \in \St^{n}$ and $0 < r < 2r_{0}$, then 
\begin{equation}
	\left |\langle x-x_0,\nu\rangle \right| \leq C_{1}r e(E,x_{0}, 2r, \nu)^{\frac{1}{2(n-1)}} \quad \forall x\in C(x_0,r,\nu)\cap \partial E,
\end{equation}
\begin{equation} \label{e:ssep1}
\left\{ x \in C(x_{0}, r, \nu) \cap E \mid \langle x - x_{0}, \nu \rangle >  C_{1} r e(E,x_{0}, 2r, \nu)^{\frac{1}{2(n-1)}} \right \} = \emptyset,
\end{equation}
and
\begin{equation} \label{e:ssep2}
\left\{ x \in C(x_{0}, r, \nu) \cap E^{c} \mid \langle x - x_{0}, \nu \rangle < - C_{1} r e(E,x_{0}, 2r, \nu)^{\frac{1}{2(n-1)}} \right\} = \emptyset.
\end{equation}
\end{theorem}

An immediate quantitative consequence of Lemma \ref{l:excessvsoscillation} and Theorem \ref{t:SL} is 

\begin{corollary} \label{c:BMO}
Given $n \ge 2$, and $C_{A} \ge 1$ there exist constants $\epsilon_{2} = \epsilon_2(n,C_{A})$ and $C_{2} = C(n,C_{A})$ (both positive and finite) such that if $E \in \cA(C_{A},r_0)$ (for some $ r_{0} > 0$) satisfies
\begin{equation} \label{e:local}
\sup_{r < r_{0}} \left(\fint_{B(x,r)\cap \partial^* E} |\nu_{E} - (\nu_{E})_{x,r}|^{2} \dif \cH^{n-1} \right)^{\frac{1}{2}} \le \epsilon_{2},
\end{equation}
for some $x\in\partial E$, then
\begin{equation} \label{e:locflat}
\sup_{\rho < r_{0}/8} \Theta(x, \rho) \le C_{2} \epsilon_{2}^{\frac{1}{n-1}}.
\end{equation}
In particular, if $\Omega \subset \R^{n}$ is a domain such that $\partial_{*} \Omega = \partial \Omega$, $\partial \Omega$ is $(n-1)$-Ahlfors regular, and $\nu_{\Omega}$ satisfies
\begin{equation} \label{e:sbmo}
\sup_{r < r_{0}} \sup_{x \in \partial \Omega} \left(\fint_{B(x,r)\cap\partial\Omega} |\nu_{\Omega} - (\nu_{\Omega})_{x,r}|^{2} \dif \cH^{n-1}\right)^{\frac{1}{2}} \le \epsilon_{2},
\end{equation} 
then 
$\Omega$ is a $C_{2} \epsilon_{2}^{\frac{1}{n-1}}$-Reifenberg flat domain.
\end{corollary}

\begin{proof}
As in Remark \ref{l:boundarycont}, $\partial \Omega = \partial_* \Omega$ and $\partial \Omega$ is Ahlfors regular imply 
	\[ \partial \Omega = \spt \mu_\Omega, \quad |\mu_{\Omega}| \text{ is Ahlfors regular}. \]
That is, $\Omega \in \mathcal{A}(C_A, r_0)$ for some constants $C_A$, and all $r_0$. Therefore the corollary is a consequence of Theorem \ref{t:SL}. 
\end{proof}

An immediate qualitative consequence of Lemma \ref{c:BMO} and Theorem \ref{t:SL} is 

\begin{corollary} \label{c:VMO}
If $\Omega \subset \R^{n}$ is a domain such that $\partial_{*} \Omega = \partial \Omega$, $\partial \Omega$ is $(n-1)$-Ahlfors regular, and $\nu_{E} \in VMO_{\loc}(\cH^{n-1} \restr \partial\Omega)$
then $\partial \Omega$ is a vanishing Reifenberg flat set.
\end{corollary}

Corollary \ref{c:BMO} also has the following quantitative consequence for $\delta$-CADs, see Definition \ref{d:CADSC}.

\begin{corollary}\label{BMOimpflatnessdomains.cor} 
Let $\Omega \subset \mathbb R^n$ be a domain with $\partial_* \Omega = \partial \Omega$ and with $(n-1)$-Ahlfors regular boundary with constant $C_A$. Further assume, if $\Omega$ is unbounded, that $\mathbb R^n \backslash \partial \Omega$ consists of two nonempty connected components. Then, there exists a $\delta_n > 0$ such that for $\delta \in (0, \delta_n]$, there exists $\epsilon_\delta < \epsilon_2$ (where $\epsilon_2 > 0$ is as in Corollary \ref{c:BMO}) such that if for every compact set $K \Subset \mathbb R^n$ there exists an $R_K > 0$ such that $\sup_{x\in \partial \Omega\cap K} \|\nu\|(x, R_K) < \epsilon_\delta$, then $\Omega$ is a $\delta$-chord-arc domain.
\end{corollary}

\section{An Application to a Two-Phase Problem For Harmonic Measure}\label{s: twophaseharmonicmeasure}
In this section, we consider a two-phase free boundary problem for harmonic measure, originally studied by Kenig-Toro in \cite{kenigtorotwophase} and later by \cite{bortz2017singular}. In particular, we complete the proof of Theorem \ref{Tchar2sideCAD.thrm}, and prove a quantitative version of it (Theorem \ref{FBthrm.thrm}).  

\subsection{The Existence of Corkscrews} \label{geometric-info}

The goal of this subsection is to show that  
the doubling of harmonic measure implies interior corkscrews (Lemma \ref{ADRandDoublingimpCS.lem}). Later, we will show that control on the oscillation of the logarithm of the Poisson kernel implies doubling. This is an important step in proving Theorem \ref{t:osckcontrolsoscn} as it will allow us use the theory of UR domains (by way of Appendix \ref{URapprox.sect}). First we recall what it means for harmonic measure to be doubling.

\begin{definition}
Let $\om \subset \rn$ be a domain with harmonic measure $\hm$.
We say that $\omega$ is locally doubling with constant $C$, if for every compact set $ K$ there exists $r_K>0$ such that
\begin{equation} \label{e:doubling}
\hm(B(x,2r)) < C\hm(B(x,r)).
\end{equation}
for all $x\in\partial\Omega\cap K$ and all $r\in (0, r_K)$. We also refer to $r_K$ as the (local) doubling condition radius.
\end{definition}

\begin{remark}
	We often assume that $r_K$ is sufficiently small compared to the distance from the pole of $\omega$ to the boundary $\pom$. This allows us to focus on local regions away from the pole, so that we can use preliminary estimates on the harmonic measure with ease.
\end{remark}
To prove estimates that are uniform on compacta, it is important to keep track what the value of each constant depends on, and in particular, whether or not it depends on the choice of compact set. For simplicity we may say the value depends on allowable constants, if it only depends on the dimension $n$ and the Ahlfors regularity constant, and does not depend on the compact set.
The following Lemma \ref{ADRandDoublingimpCS.lem}, which might be considered folklore, shows the existence of interior corkscrews given the doubling of harmonic measure. This is an essential step, as it allows us to gain topological information on $\Omega$ from the regularity of the Poisson kernel. We sketch the proof here, which is a small modification of the proof of \cite[Lemma 3.14]{hofmann2015uniform} (see also \cite[Lemma 4.24]{hofmann2017weak}). 

\begin{lemma}\label{ADRandDoublingimpCS.lem}
Let $\Omega \subset \R^{n}$ be a domain whose boundary is Ahlfors regular with constant $C_{A}$. Fix $X_{0} \in \Omega$. Suppose $\omega^{X_{0}}$ is locally doubling with constant $C_{0}$. There exists an $\eta = \eta(n, C_A) > 0$ such that for every closed ball $K$, if $r_K << \delta(X_0)$ is the doubling radius of $\omega^{X_0}$ in $K$, then  $\Omega$ admits an interior corkscrew ball at every $x \in \partial \Omega \cap K$ up to radius $s_{K} := \eta r_{K}$ with constant $C_{1} = C(n, C_{A}, C_{0})$. 
\end{lemma}

\begin{proof}
Fix the closed ball $K$ and recall that $r_K$ is the local doubling radius. The proof of this lemma requires a slight modification of the argument in \cite[Lemma 3.14]{hofmann2015uniform}.
Recall the following relationship between the Green function and the harmonic measure. For $\Phi \in C_c^\infty(\ree)$ 
\begin{equation}\label{representationofGhm.eq}
\int_{\pom} \Phi(y) \, d\hm^X(y)  - \Phi(X) = -\iint_{\om} \nabla \G(X,Y) \nabla \Phi(Y) \, dY, \quad \text{ a.e. } X \in \om,
\end{equation}
where $\hm := \hm^X$ and $\G(Y):= \G(X,Y)$ are the harmonic measure and Green's function for $\om$ with pole at $X$. 

It was proven in \cite[Lemma 2.40]{hofmann2015uniform} that there exists $\kappa_0 > 2$ depending only on dimension and the Ahlfors regularity constant such that for all $x\in \pom$ and $0< r <  \min\{ \delta(X)/\kappa_0, \diam(\pom)\}$, for $B = B(x,r)$
\begin{equation}\label{Carltype1.eq}
\sup_{\tfrac{1}{2}B} \G(Y) \lesssim \frac{1}{|B|}\iint_{B} \G(Y) \, dY \lesssim r \frac{\hm(C B)}{\sigma(C B)},
\end{equation}
where all implicit (and explicit) constants depend only on dimension and the Ahlfors regularity constant.

Now let $x \in \pom$ and $0< r < \min\{\delta(X_0)/\kappa_0, 10^{-3}\diam(\pom), 10^{-3}r_K/C\}$, where $r_K$ is the doubling condition radius for $\hm$ and $C$ is as in \eqref{Carltype1.eq}. Without loss of generality we may assume $r_K \ll \min\{ \delta(X),  \diam(\pom) \}$, so that the above minimum equals $10^{-3} r_K/C$. 
Set $B := B(x,r)$ and $\Phi \in C_c^\infty(\tfrac{1}{2}B)$ be such that $0 \le \Phi \le 1$, $\Phi \equiv 1$ on $\tfrac{1}{100}B$ and $|\nabla \Phi| \lesssim 8/r$. Using \eqref{representationofGhm.eq} with $X = X_0$\footnote{We may move $X_0$ slightly using the Harnack inequality.}
we obtain
\begin{equation}\label{initGhmrepsplit.eq}
\begin{split}
r \hm(\tfrac{1}{100}B) &\le r \int_{\pom \cap \tfrac{1}{100}B} \Phi(y) \, d\hm(y) =  -r\iint_{\om} \nabla \G(Y) \nabla \Phi(Y) \, dY
\\ &\le 8 \iint_{\om \cap \tfrac{1}{2}B} |\nabla \G(Y)|\, dY
\\ & \le 8  \iint_{\left(\tfrac{1}{2}B \cap \om\right) \setminus \Sigma_\rho(r)} |\nabla \G(Y)|\, dY 
+  8 \iint_{\tfrac{1}{2}B \cap \Sigma_\rho(r)} |\nabla \G(Y)|\, dY
\\ & = \A + \B,
\end{split}
\end{equation}
where $\Sigma_\rho(r)$ is the `boundary strip', $\Sigma_\rho(r) := \{Y \in \om: \delta(Y) \le \rho r\}$ and $\rho > 0$ is a small number to be chosen momentarily. Let $\W = \{I\}$ be a Whitney decomposition of $\om$ and let $\I := \{I \in \W: I \cap \tfrac{1}{2}B \cap \Sigma_\rho(r) \neq \emptyset\}$. Then using standard interior estimates (the Caccioppoli inequality and the Moser estimate)
\begin{equation}\label{Bwhitneybreak.eq}
\B \le 8 \sum_{I \in \I} \iint_{I} |\nabla \G(Y)|\, dY \le C' \sum_{I \in \I} \ell(I)^{n-1} |\G(Y_I)|,
\end{equation}
where $Y_I$ is the center of the Whitney cube $I$ and $\ell(I)$ is the side length of $I$. 
For each $I \in \I$ we use the H\"older continuity at the boundary of the Green function (which only depends on dimension and the Ahlfors regularity constant), in conjunction with \eqref{Carltype1.eq}, to obtain the estimate
\[\G(Y_I) \lesssim \left(\frac{\ell(I)}{r}\right)^{\alpha} \frac{1}{|2B|} \iint_{2B \cap \om} \G(Y) \, dY
\lesssim \left(\frac{\ell(I)}{r}\right)^{\alpha}r\frac{\hm(CB)}{\sigma(C B)}.\]
Summing over $I \in \I$, and using an elementary geometric argument, whose proof we temporarily postpone,
we have that
\begin{equation}\label{tt-5}
\B \lesssim \rho^{\alpha} r \hm(CB) \lesssim \rho^{\alpha} r \hm(\tfrac{1}{100}B),
\end{equation}
where we used that the harmonic measure is doubling up to $r_K$. 

{  Then there exists $\rho > 0$ depending on $C_{0}, n,$ and $C$ (which depended additionally on $C_{A}$), small enough so that the upper bound in \eqref{tt-5} can be absorbed in the left hand side of \eqref{initGhmrepsplit.eq} at which point we have} 
$$\A = 8\iint_{(\tfrac{1}{2}B \cap \om) \setminus \Sigma_\rho(r)} |\nabla \G(Y)|\, dY  \ge \frac12 r \hm(\tfrac{1}{100}B) > 0.$$
Since $\A > 0$ there exists a point $Y_B \in \tfrac{1}{2}B \cap \om$ such that $\delta(Y_B) > \rho r$, which shows that $\om$ satisfies the $(\tfrac{1}{\rho}, R_0)$-interior corkscrew condition, where \newline
\[ R_0 = \min\{\delta(X)/\kappa_0, 10^{-3}\diam(\pom), 10^{-3}r_K/C\} = 10^{-3} r_K/C =:s_K. \]
Hence we finish the proof of the lemma with constant $\eta:= 10^{-3}/C$.

Now we sketch the `elementary geometric argument', that is, how we used the estimate on $\G(Y_I)$ and \eqref{Bwhitneybreak.eq} to obtain \eqref{tt-5}. If we let 
\[\widetilde{\I} := \{I \in \W: I \cap \tfrac{1}{2}B  \neq \emptyset\}\]
then we observe that the Whitney property of each $I \in \widetilde{\I}$ ensures that $\ell(I) \lesssim r$ and for each $I \in \widetilde{\I}$ there exists $\hat{x}_I$ in $B(x,Cr) \cap \pom$ such that 
\[\ell(I) \approx \dist(I, \pom) \approx \dist(\hat{x}_I, Y), \quad \forall Y \in I.\]
Now fix $k$ such that $2^{-k} \lesssim \rho r$, denote $\widetilde{\I}_k:= \{I \in \widetilde{\I}: \ell(I) = 2^{-k}\}$ and  cover $B(x,Cr) \cap \pom$ by balls $\{B_{k,j}\}_{j} = \{B(x_{k,j}, 2^{-k})\}$ with $x_{k,j} \in \pom$ such that  $\{\tfrac{1}{5}B_{k,j}\}_j$ are disjoint. Using Ahlfors regularity to compare surface areas we see that for each fixed $k$, 
\[\#\{B_{k,j}\}_{j} \approx r^{n-1}2^{k(n-1)}.\]
Now for each $I \in \widetilde{\I}_k$ associate an index $j$ such that $x_{I} \in B_{k,j}$ and notice  we have $\dist(Y, x_{k,j}) \lesssim 2^{-k}$ for all $Y \in I$. Since the $I \in \widetilde{\I}_k$ are disjoint, comparing volumes we have that for fixed $j$  
\[\#\{I \in \widetilde{\I}_k: I \text{ is associated to } j\} \le C,\] where $C$ depends on dimension. It follows from our bound on $\#\{B_{k,j}\}_{j}$ that 
\[\#\widetilde{\I}_k \lesssim r^{n-1}2^{k(n-1)}.\]
Now breaking the sum over $k$ in \eqref{Bwhitneybreak.eq} and using our bound for $\G(Y_I)$ we obtain
\begin{align*}
\B &\lesssim \hm(CB)r^{2-n -\alpha} \sum_{k \gtrsim -\log_2(\rho r)} \sum_{I \in \widetilde{\I}_k} 2^{-k(n-1 + \alpha)} 
\\& \lesssim  \hm(CB)r^{2-n -\alpha} \sum_{k \gtrsim -\log_2(\rho r)}r^{n-1}2^{k(n-1)} 2^{-k(n-1 + \alpha)}
\\ &\lesssim  \rho^\alpha r\hm(CB)
\end{align*}
as desired, where we used $\sigma(C B)\approx r^{n-1}$ in the first line.
\end{proof}

One immediate corollary is that domains with Ahlfors regular boundaries have uniformly rectifiable boundaries whenever the interior and exterior harmonic measures are doubling.

\begin{corollary}\label{hmdoubcorr.lem}
Suppose $\om^+ \subset \rn$ and $\om^-= \rn\setminus \overline{\om^+}$ are domains with common topological boundary $\pom := \pom^+ = \pom^-$ and $\diam(\pom^+) < \infty$, which has the additional property that $\pom$ is $(n-1)$-Ahlfors regular. Suppose further that there exists $X^+ \in \om^+$ and $X^- \in \om^-$ such that the harmonic measures $\hm_\pm^{X^\pm}$ are doubling. Then $\pom$ is uniformly rectifiable and $\pom = \partial_*\om$. In particular, $\om^\pm$ are UR domains.
\end{corollary}

\subsection{A Localization Result} \label{localization-tt}

The major technical result of this section is Theorem \ref{t:osckcontrolsoscn}, which, roughly, states that the local oscillation of the Poisson kernel controls the local oscillation of the unit normal. 
Perhaps contrary to the spirit of a ``localized result" the {\it scale} at which we get control of the oscillation of the unit normal depends on the compact set; however the quantitative control does not, see \eqref{tt-10} and \eqref{tt-11}.

Our main tool in the proof of Theorem \ref{t:osckcontrolsoscn} is the single layer potential, we recall its definition now:

\begin{definition}[Riesz transform and the single layer potential]\label{riesztransdef}
Let $F\subset \rn$ be a closed $(n-1)$-Ahlfors regular set with surface measure $\sigma=\cH^{n-1}\restr F$.
We define the (vector-valued) Riesz kernel as
\begin{equation}\label{RieszKern}
\mathcal{K}(X) = \tilde{c}_n\frac{X}{|X|^{n}}
\end{equation}
where $\tilde{c}_n$ is chosen so that $\mathcal{K}(X) = \nabla \frac{c_n}{|X|^{n-2}}$ and $c_n$ is such that  $-\Delta \frac{c_n}{|X|^{n-2}} = \delta_0$ (here $\delta_0$ is the Dirac mass at the origin).
 
Let $f \in L^p(d\sigma)$ for some $p \in [1,n-1)$.  We define the Riesz transform of $f$ (relative to $F$)
\begin{equation}\label{RieszXform}
\mathcal{R}f(X) := \mathcal{K} \ast (f\sigma)(X) = \int_{F} \mathcal{K}(X-y)f(y)\, d\sigma(y) \quad X \in \rn\setminus F,
\end{equation} 
as well as the truncated Riesz transforms for $X \in F$
$$\mathcal{R}_\eps f(X):= \int_{F\,\cap\,\{y:|X-y|>\eps\}} \mathcal{K}(X-y)f(y)\, d\sigma(y)\,,\qquad  \eps>0\,.$$
We define $\mathcal{S}$ the (harmonic) single layer potential of $f$ relative to $F$ to be
\begin{equation}
\mathcal{S}f(X): = \int_{F} \mathcal{E}(X -y) f(y) \, d\sigma(y),
\end{equation}
where $\mathcal{E}(X) = c_n|X|^{2-n}$. 
\end{definition}

\begin{remark} For $f$ as above we have that $\mathcal{S}f(X)$ makes sense as an absolutely convergent integral for $X \not \in F$. To see this we may use the upper Ahlfors regularity to break the boundary up into dyadic annuli centered at $x_0 \in F$ such that $\dist(X, F) = |X - x_0|$ and see that $\mathcal{E}(X -y)$ is in $L^{p'}(d\sigma)$ for all  $p \in [1,n-1)$, where $p'$ is the H\"older conjugate exponent to $p$ (albeit with bounds depending on $X$). Notice also that for such $f$, $\nabla \mathcal{S}f(X) = \mathcal{R}f(X)$  for $X \not\in F$ and $\mathcal{R}f(X)$ makes sense as an absolutely convergent integral for $X \not\in F$ (here we use the same argument as for $\mathcal{E}$ to show that $\mathcal{K}(X-y) \in L^{p'}(d\sigma)$ for $p \in [1,\infty)$).  To see $\nabla \mathcal{S}f(X) = \mathcal{R}f(X)$ for $X \not\in F$ we form the difference quotients for $\mathcal{S}f$ and use the dominated convergence theorem. Every function we apply the single layer potential to in the proof of Theorem \ref{t:osckcontrolsoscn} is in the space $L^1(d\sigma)$.
\end{remark}

The singular layer potential is useful in that it generates solutions to the Neumann problem (see e.g. \cite[Section 5.5]{hofmann2010singular}). However, in order to make sense of boundary data in a rough domain we need to introduce the concept of non-tangential regions: 

\begin{definition}[Nontangential approach region and maximal function]\label{NTapp} Fix $\alpha > 0$ and let $\om$ be a domain,
then for $x \in \pom$ we define the nontangential approach region (or ``cone")
\begin{equation}\label{NTapp1}
\Gamma(x) = \Gamma_\alpha(x) = \{Y \in \om : |Y - x| < (1 + \alpha)\delta(Y)\}. 
\end{equation}
We also define the nontangential maximal function for $u: \om \to \re$
\begin{equation}\label{NTapp2}
\mathcal{N}u(x) = \mathcal{N}_\alpha u(x) = \sup_{Y \in \Gamma_\alpha(x)}|u(Y)|, \quad x \in \pom.
\end{equation}
We make the convention that $\mathcal{N}u(x) = 0$ when $\Gamma_\alpha(x)=\emptyset$\footnote{In the settings treated here, this is always a set of $\cH^{n-1}$ measure zero \cite[Proposition 2.9]{hofmann2010singular}.}  and that $\alpha = 1$ when no subscript appears in $\Gamma$.
\end{definition}

The relationship between the two definitions above is made clear in the following two lemmas:

\begin{lemma}[\cite{hofmann2010singular}]\label{NTmaxbound.lem}
Suppose that $\om$ is a UR domain (recall Definition \ref{URdom}) and $f \in L^{q}(d \sigma)$ for some $q \in [1,n-1)$. For all $p \in (1,\infty)$ we have
\begin{equation}\label{eq15}
\lVert \mathcal{N}(\nabla \mathcal{S} f)\rVert_{L^p(d\sigma)}  \le C \lVert f \rVert_{L^p(d\sigma)},
\end{equation}
where $C$ depends on the UR character of $\pom$, dimension, $p$, and the aperture of the cones defining
$\mathcal{N}$.  
\end{lemma}
The bound for the non-tangential maximal function of $\nabla \mathcal{S}f$ follows from uniform bounds for the truncated singular integrals  \cite{david1991singular}, plus a  
Cotlar Lemma argument;  the details may be found in \cite[Proposition 3.20]{hofmann2010singular}.

In addition, we have the following result proved in \cite{hofmann2010singular}.
\begin{lemma}[\cite{hofmann2010singular} Proposition 3.30]\label{hmtlemma} 
If $\om$ is a UR domain, whose measure theoretic and topological boundary agree up to a set of $\cH^{n-1}$ measure zero, 
then for  a.e. $x \in \pom$, and for all   $f \in L^p(d\sigma)$, $1\le p<n-1$,
\begin{equation}\label{eq16}
\lim_{\substack{Z \to x \\ Z \in \Gamma^-(x)}}\nabla \mathcal{S} f (Z) = -\frac{1}{2}\nu(x)f(x) + \mathcal{T}f(x)\,,
\end{equation}
and 
\begin{equation}\label{eq30}
\lim_{\substack{Z \to x \\ Z \in \Gamma^+(x)}}\nabla \mathcal{S} f (Z) = \frac{1}{2}\nu(x)f(x) + \mathcal{T}f(x)\,.
\end{equation}
where $\Gamma^+(x)$ is the cone at 
$x$ relative to $\om$, $\Gamma^-(x)$ is the cone at $x$ relative to ${\Omega_{\rm ext}}$, $\nu$ is the unit outer normal to $\om$, and $\mathcal{T}$ is a (vector-valued) principal value singular integral operator:
\[ \mathcal{T} f(x) = \lim_{\epsilon \to 0+} \int_{y\in \pom \setminus B(x,\epsilon)} \nabla \mathcal{E}(x-y) f(y) d\sigma(y). \]
\end{lemma}

\begin{remark} As in \cite{bortz2017singular}, we have taken our fundamental solution to be positive, so 
for that reason there are some changes in sign in both \eqref{eq16} and \eqref{eq30} 
as compared to the formulation in \cite{hofmann2010singular}. \end{remark}

Next we show that if $\log k$ has \textit{small} BMO norm, the measure $\omega = k \, d\sigma$ is doubling. The proof uses the fact that $\sigma$ is doubling. We remark that in general, 
the fact that $\|\log k\|_{BMO} <\infty$ or that $k$ satisfies a reverse H\"older inequality, does not ensure that $\omega = k \, d\sigma$ is doubling, see the discussions and example in \cite[Chapter I]{ST89}.

\begin{lemma}\label{lm:doubling}
	Let $\sigma$ be a doubling measure on $\rn$ and $\omega = k\, d\sigma$ be another Radon measure with $0\leq k \in L^1_{\loc}(d\sigma)$. There exists $\tau_0 > 0$ depending on the doubling constant of $\sigma$, such that if 
	\begin{equation}\label{tt-20}
	\|\log k\|_{*}(B(x_0,4r_0)) < \tau\le\tau_0 \hbox{  for some   }x_0\in \spt \sigma\hbox{  and  }r_0>0, 
	\end{equation}
	then the following holds for $B\subset B(x_0,2r_0)$ with $B$ a ball centered on $\spt\sigma$:
	\begin{enumerate}
		\item There is a constant $C$ depending on $n$ such that
			\begin{equation}\label{eq:avglog}
				 \frac{1}{1+C\tau} \, \fint_B k \, d\sigma \leq  e^{\fint_B \log k \, d\sigma} \leq \fint_B k \, d\sigma = \frac{\omega(B)}{\sigma(B)}.
				\end{equation} 
		\item Given $p>1$, there exists $\tau(p)\le \tau_0$  such that if \eqref{tt-20} holds with $\tau\le\tau(p)$ then
		for any Borel set $E\subset B$, where $B$ is as before
	\begin{equation}\label{eq:Apratio}
		\frac{\omega(E)}{\omega(B)} \geq c(p,\tau) \, \left( \frac{\sigma(E)}{\sigma(B)} \right)^p.
	\end{equation} 
	Here the constant $c(p,\tau) \to 1$ as $\tau \to 0$.
		\item In particular, for $x\in\spt\sigma$ such that $B(x,2r)\subset B(x_0,2r_0)$
			\begin{equation}\label{eq:doubling}
			\omega(B(x,2r)) \leq C \omega ( B(x, r) ),
			\end{equation} 
			where the constant $C$ depends on $n$ and the doubling constant of $\sigma$.
		\item Given $r>1$ there exists $\tilde\tau(r)\le \tau_0$  such that if \eqref{tt-20} holds with $\tau\le\tilde\tau(r)$ then
		the weight $k$ satisfies the reverse H\"older inequality for $r$, i.e.
			\begin{equation}\label{eq:RHr}
				\left( \fint_B k^r \, d\sigma \right)^{1/r} \leq C(r, \tau) \, \fint_B k \, d\sigma.
			\end{equation}
			Here the constant $C(r, \tau) \to 1$ as $\tau \to 0$.
	\end{enumerate}
\end{lemma}

\begin{proof}
	By the local version of John-Nirenberg inequality for doubling measures (see \cite[Theorem 5.2]{ABKY}) we have 
	\[ \sigma\left( \left\{x\in B: |\log k(x) - (\log k)_B | > \lambda \right\} \right) \leq C_1 e^{-C_2\frac{\lambda}{\tau}} \sigma(B) \]
	for all $\lambda>0$, where the constants $C_1$ and $C_2$ depend on the doubling constant for $\sigma$. Therefore
	\begin{align}
		& \fint_B e^{|\log k - (\log k)_B|} d \sigma \nonumber \\
		 &\hspace{-.1in} \qquad = \frac{1}{\sigma(B)} \int_0^\infty \sigma \left( \left\{x\in B: e^{|\log k(x) - (\log k)_B |} > s \right\} \right) ds \nonumber \\
		& \hspace{-.1in} \qquad   \le \frac{1}{\sigma(B)} \int_0^1 \sigma(B) ds +\frac{1}{\sigma(B)} \int_0^\infty \sigma \left( \left\{x\in B: |\log k(x) - (\log k)_B | > \lambda \right\} \right) e^\lambda d\lambda \nonumber \\
		&\hspace{-.1in} \qquad \leq  1+ C_1 \int_0^\infty e^{-\frac{C_2}{\tau} \lambda + \lambda} d\lambda \nonumber \\
		&\hspace{-.1in} \qquad \leq 1+ C\tau, \label{eq:explogavg}
	\end{align}
	if $\tau$ is sufficiently small (depending on the constant $C_2$). \eqref{eq:avglog} follows immediately.
	
	Similarly, provided $\tau$ is small enough depending on $p$ we also have
	\begin{equation}\label{eq:pfAp}
		\fint_B e^{\frac{1}{p-1}|\log k - (\log k)_B|} d \sigma \leq 1+ C_p \tau. 
	\end{equation} 
	Henceforth, $\tau_0 >0$ is chosen so that \eqref{eq:pfAp} holds with $p = 2$ and $\tau \le \tau_0$.
	Let $q=p/(p-1)$ be the H\"older conjugate of $p$.
	It follows that
	\begin{align*}
	\fint_B k \, d\sigma \cdot \left( \fint_B k^{-\frac{q}{p}} d\sigma \right)^{\frac{p}{q}} & =  \fint_B e^{\log k} d\sigma \cdot \left(\fint_B e^{-\frac{1}{p-1} \log k} d\sigma \right)^{p-1}\\
	& =  \fint_B e^{\log k - (\log k)_B} d\sigma \cdot \left( \fint_B e^{-\frac{1}{p-1} \left( \log k - (\log k)_B \right)} d\sigma \right)^{p-1} \\
	& \leq \fint_B e^{|\log k - (\log k)_B|} d\sigma \cdot \left( \fint_B e^{\frac{1}{p-1} |\log k - (\log k)_B|} d\sigma \right)^{p-1 } \\
	& \leq (1+C_p \tau)^p, 
\end{align*} 
i.e. $k\in A_p(\sigma)$, where $A_p$ is the Muckenhaupt class with power $p>1$.

Let $g\geq 0$ be an arbitrary measurable function on $B$. We have
\begin{align*}
	\int_B g \, d \sigma & \leq \left( \int_B g^p k \, d\sigma \right)^{\frac{1}{p}} \left( \int_B k^{-\frac{q}{p} } \, d\sigma \right)^{\frac{1}{q}} \\
	& \leq  (1+C_p \tau) \sigma(B) \left( \int_B k\, d\sigma \right)^{-\frac{1}{p}} \left( \int_B g^p k \, d\sigma \right)^{\frac{1}{p}}. \\
\end{align*}
In particular for any Borel set $E\subset B$, by plugging in the above inequality $g= \chi_E$, we get
\[ \frac{\sigma(E)}{\sigma(B)} \leq (1+C_p \tau) \left( \frac{\omega(E)}{\omega(B)} \right)^{\frac{1}{p}}, \]
or equivalently
\[ \frac{\omega(E)}{\omega(B)} \geq c(p,\tau) \left( \frac{\sigma(E)}{\sigma(B)} \right)^p \]
with $c(p,\tau) =1/(1+C_p\tau)^p$. The doubling property \eqref{eq:doubling} follows by taking $E= (1/2)B$, $p = 2$ and $\tau = \tau_0$.

Let $r>1$, then \eqref{eq:pfAp} applied to $p=1+1/r$ implies that for $\tau$ small enough depending on $r$ we have
\[ \fint_B k^r \, d\sigma \leq (1+C_r \tau) e^{r (\log k)_B }. \]
Taking $r$-th root on both sides of the inequality and using \eqref{eq:avglog}, we get
\[ \left( \fint_B k^r \, d\sigma \right)^{1/r} \leq (1+C_r \tau)^{1/r} \, e^{(\log k)_B} \leq (1+ C_r \tau)^{1/r} \, \fint_B k \, d\sigma, \]
i.e. $k\in RH_r(\sigma)$, where $RH_r$ denotes weight that satisfies the reverse H\"older inequality with power $r>1$.
\end{proof}

After we establish the reverse H\"older inequality \eqref{eq:RHr}, one can show 
\begin{equation}\label{1-k/aest.eq}
	\left( \fint_B \left| 1- \frac{k}{a} \right|^2  \, d\sigma \right)^{1/2} \leq C\left( \|\log k \|_*(4B) \right)^{1/8} \leq C \tau^{1/8},
\end{equation}
where $a= e^{\fint_B \log k \, d\sigma}$. For details of the proof we refer interested readers to \cite[Lemma 1.33]{bortz2017singular}.

\bigskip

The following result states that control on the oscillation of the logarithm of the interior and exterior Poisson kernels provides control on the oscillation of the unit normal.

\begin{theorem}\label{t:osckcontrolsoscn}
Let $\Omega^+ \subset \R^n, \Omega^- = \R^n\backslash \overline{\Omega^+}$ be domains with common (topological) boundary, $\pom^+ = \pom^- \equiv \pom$. Assume that $\pom$ is $(n-1)$-Ahlfors regular and let $X^\pm \in \Omega^{\pm}$ be such that $k^{\pm} = \frac{d\omega^{\pm}}{d\sigma}$ exist.  Here $\hm^{\pm} = \hm_{\pm}^{X^\pm}$, where $\hm_{\pm}^{X^\pm}$ is the harmonic measure for $\om^\pm$ with pole at $X^\pm$. Given $\epsilon>0$ there exists 
$\kappa_1>0$ depending on $\delta(X^\pm)$, $\epsilon$, $n$ and the Ahlfors regularity constant $C_A$ such that if $\log k^\pm\in BMO_{\loc}( \sigma)$ with constant $0<\kappa\le \kappa_1$, then $\nu\in BMO_{\loc}( \sigma)$ with 
constant at most $\epsilon$. In particular, if $\log k^\pm\in VMO_{\loc}( \sigma)$, then $\nu\in VMO_{\loc}( \sigma)$.
 \end{theorem}
 \begin{remark}
 	The proof of the above theorem yields a quantitative estimate, see \eqref{tt-10} and \eqref{tt-11}.
 \end{remark}
\begin{proof}
Let $A > 2$ be a constant depending on dimension and the Ahlfors regularity constant\footnote{We use $A$ to simplify notation. In fact, we take $A=C_{3}$ as in Lemma \ref{Dcubes.lem} and used in Lemma \ref{URdomApprox.lem} } such that if $x_0 \in \pom$ and $r_0 \in (0, \diam \pom)$ then there exists\footnote{See Remark \ref{goodchoicecube.rmk}.} a dyadic cube $Q$ as in Lemma \ref{Dcubes.lem} such that 
\[\Delta(x_0, r_0/A) \subset Q \subset \Delta(x_0,r_0).\]

Let $\tau(p)$ be as in Lemma \ref{lm:doubling} such that \eqref{eq:Apratio} holds with power $p= 1+1/(2(n-1))$.
Suppose that $\log k^\pm \in$ BMO$_{\loc}( \sigma)$ with BMO$_{\loc}$ semi-norm $\kappa$ satisfying $\kappa \in (0, \kappa_1)$, where $\kappa_1 \le \tau(p)$ will be determined after \eqref{tt-8}. Notice that in the case when  $\log k^\pm\in VMO_{\loc}( \sigma)$ this holds for every $\kappa > 0$. Fix $B^* = B(y_0, 4R)$ for some $y_0 \in \pom$ and  $R \in (0, \diam(\pom)/4)$ and set $\widetilde{B} = \tfrac{1}{4}\overline{B^*}$.  Since $\log k^\pm \in$ BMO$_{\loc}( \sigma)$ with constant $\kappa$, there exists a radius $r_0 = r_0(\tau(p), B^*) < c\min\{R, \delta(X^\pm)\}$ (with $c > 0$ depending on dimension and Ahlfors regularity) such that  
\[\|\log k\|_*(B(z_0,2 r_0)) < \kappa, \quad \forall z_0 \in B^* \cap \pom\] The proof of 
Lemma \ref{lm:doubling} establishes that $\omega^\pm$ are doubling\footnote{Here we have uniform control on the doubling constant by Lemma \ref{lm:doubling} and the choice of $\kappa_1$.} up to radius $r_0$ on balls centered on $B^* \cap \pom$, with a doubling constant depending on $n$ and $C_A$. 
Moreover by choice of $c$ and Lemma 
\ref{ADRandDoublingimpCS.lem}, the domains $\Omega^\pm$ both admit an interior corkscrew ball for every $x\in B^* \cap \partial\Omega$ up to radius $r_0$. Thus, we record for later use that, in the language of Appendix \ref{URapprox.sect}, $\om$ satisfies the $(x_0,M_0, r_0)$-DLTSCS\footnote{This is a local two-sided corkscrew condition, see Definition \ref{DLTSCS.def}.} for all $x_0 \in \widetilde{B}$.

From this point forth, $x_0$ will denote an arbitrary point in $\widetilde{B} \cap \pom$. Let $1<M<\infty$ and $\theta\in (0,1)$ be determined later. For $x\in B(x_0,r_0/(20A))\cap\pom$, let $r\in (0,\theta r_0)$ be such that $\Delta := \Delta(x,r) \subset\Delta^* := \Delta(x, Mr) \subset  B(x_0, r_0/(5A))$. 

 For any $y, z \in \Delta$, we let $y^*$ and $z^*$ denote arbitrary points in the non-tangential approach regions in $\om^-$, $\Gamma^-(y) \cap B(y, r/2)$ and $\Gamma^-(z) \cap B(z, r/2)$, respectively. Following \cite[Theorem 1.1]{bortz2017singular} we first show
\begin{eqnarray}\label{2.1BH.eq}
&&\left(\fint_\Delta \left| \nabla \mathcal{S} 1_{\Delta^*}(z^\ast) -\fint_\Delta \nabla \mathcal{S} 1_{\Delta^*}(y^*) \, d\sigma(y) \right|^2 \, d \sigma(z) \right)^\frac{1}{2} \\
&&\qquad\qquad\qquad\le \frac{C_1}{\omega(B(x_0, r_0/(5A)))} \cdot\left(\frac{r}{r_0}\right)^{1/2} \cdot\frac{1}{\sqrt{M}}+  C_2M^{\frac{n-1}{2}}\kappa^{\frac{1}{8}}+ \frac{C_3}{M}, \nonumber
\end{eqnarray}
where $\omega$ is the harmonic measure of $\Omega^+$ with pole $X^+$, and the constants $C_1, C_2, C_3 > 0$ only depend on $n$, the Ahlfors regularity constant $C_A$ and $\delta(X^\pm)$. In particular, 
$\omega = k^+d\sigma$.
We decompose $1_{\Delta^*}$ as
\begin{equation}\label{inddecomp.eq}
1_{\Delta^*} = \left[\left(1 - \frac{k^+}{a}\right)1_{\Delta^*}\right] + \left[\frac{k^+}{a}\right] - \left[\left(\frac{k^+}{a}\right)1_{(\Delta^*)^c}\right],
\end{equation}
where $a = a_{x,Mr} = e^{\fint_{\Delta_*} \log k^\pm \, d\sigma}$.
We want to estimate the left hand side of \eqref{2.1BH.eq} by using this decomposition and the triangle inequality. This gives three terms, which we denote as $\RNum{1}, \RNum{2}$ and $\RNum{3}$: 
\begin{equation}\label{eq19}
\RNum{1} = \left(\fint_\Delta \left| \nabla \mathcal{S} \left[\left(1 - \frac{k}{a}\right)1_{\Delta^*}\right] (z^\ast) -\fint_\Delta \nabla \mathcal{S} \left[\left(1 - \frac{k}{a}\right)1_{\Delta^*}\right](y^*) \, d\sigma(y) \right|^2 \, d \sigma(z) \right)^\frac{1}{2},
\end{equation}
\begin{equation}\label{eq20}
\RNum{2} = \left(\fint_\Delta \left| \nabla \mathcal{S} \left[\frac{k}{a}\right] (z^\ast) -\fint_\Delta \nabla \mathcal{S} \left[\frac{k}{a}\right](y^*) \, d\sigma(y) \right|^2 \, d \sigma(z) \right)^\frac{1}{2},
\end{equation}
and
\begin{equation}\label{eq21}
\RNum{3} = \left(\fint_\Delta \left| \nabla \mathcal{S} \left[\left(\frac{k}{a}\right)1_{(\Delta^*)^c}\right] (z^\ast) -\fint_\Delta \nabla \mathcal{S} \left[\left(\frac{k}{a}\right)1_{(\Delta^*)^c}\right](y^*) \, d\sigma(y) \right|^2 \, d \sigma(z) \right)^\frac{1}{2}.
\end{equation}
For simplicity we drop the super-index and write $k=k^+$.
We will leave the estimate of $\RNum{1}$ for last as it requires the use of the localization Lemma \ref{URdomApprox.lem}.
 
For $\RNum{2}$, we recall that $k = k^{+}$ is the Poisson kernel for $\om$ with pole at $X^+$.  Moreover,
$\mathcal{E}(\cdot - z^*)$ and $\mathcal{E}(\cdot - y^*)$ are 
harmonic in $\om$ since $z^*, y^* \in \Omega^{-}$, and decay to 0 at infinity, and are therefore equal 
to their respective Poisson integrals in $\om$.   Consequently,
\begin{equation}\label{eq24}
\RNum{2} 
\le \frac{1}{a}  \left(\fint_\Delta\fint_\Delta \left| \nabla \mathcal{E}({X^+} - z^*)  - \nabla \mathcal{E}({X^+} - y^*) \, d\sigma(y) \right|^2 \, d \sigma(z) 
\right)^\frac{1}{2}\,.
\end{equation}
Note that, since $y^*,z^* \in B(x,2r)$ and $|X^+ - x| > r_0$ 
\begin{equation*}\label{eq37}
\left| \nabla \mathcal{E}({X^+} - z^*)  - \nabla \mathcal{E}({X^+} - y^*)   \right| 
\lesssim \frac{r}{r_0^{n}}.
\end{equation*}
Then continuing \eqref{eq24}, we have, using \eqref{eq:Apratio} with power $p=1+1/(2(n-1))$,
\begin{multline}\label{eq38}
\begin{split}
\RNum{2}  &\lesssim \frac{1}{a r_0^n} r 
\approx \frac{\sigma(\Delta^*)}{r_0^n\hm(\Delta^*)} r  =\frac{\sigma(\Delta^*)}{\hm(B(x_0,r_0/(5A)))}\frac{\hm(B(x_0,r_0/(5A)))}{r_0^n\hm(\Delta^*)} r 
\\ & \leq \frac{C}{\hm(B(x_0, r_0/(5A)))} \,\left(\frac{Mr}{r_0}\right)^{n-1}\left(\frac{r_0}{Mr} \right)^{n-\frac{1}{2}}\frac{r}{r_0} \\
&  \leq \frac{C}{\hm(B(x_0, r_0/(5A)))}
\left(\frac{r}{r_0}\right)^{\frac{1}{2}}\cdot\frac{1}{\sqrt{M}},
\end{split}
\end{multline}
where $C > 0$ depends on $n$ and the Ahlfors regularity constant.

For $\RNum{3}$, we use basic Calder\'{o}n-Zygmund type estimates as follows.  Let 
$$\Delta_j := \Delta(x, 2^jr)\,,\qquad A_j := \Delta_j \setminus \Delta_{j-1}\,,$$ 
so that
\begin{multline}\label{eq25}
\RNum{3} = \\ 
 \left(\fint_\Delta \left|\fint_\Delta 
 \left(\nabla \mathcal{S} \left[\left(\frac{k}{a}\right)
 1_{(\Delta^*)^c}\right] (z^\ast) -  \nabla \mathcal{S} \left[\left(\frac{k}{a}\right)
 1_{(\Delta^*)^c}\right](y^*)\right) \, d\sigma(y) \right|^2 \, d \sigma(z) \right)^\frac{1}{2} \\
 = \left(\fint_\Delta \left|\fint_\Delta \int_{\pom\setminus \Delta^*} \right[\nabla 
\mathcal{E}(z^* - w) -\nabla \mathcal{E}(y^*-w)\left]\frac{k(w)}{a} \, d\sigma(w)  
\, d\sigma(y) \right|^2 \, d \sigma(z) \right)^\frac{1}{2}  \\
\le \sum_{\{j\mid 2^j \ge M\}} \left(\fint_\Delta \left[\fint_\Delta   
 \int_{A_j} \left|\nabla \mathcal{E}(z^* - w) -\nabla \mathcal{E}(y^*-w)\right|\frac{k(w)}{a} \, d\sigma(w) 
  \, d\sigma(y) \right]^2 \, d \sigma(z) \right)^\frac{1}{2} \\
 \lesssim  \sum_{\{j \mid 2^j \ge M \}} \left(\fint_\Delta \left[\fint_\Delta  \int_{A_j} 
\frac{r}{(2^jr)^{n}}\frac{k(w)}{a} \, d\sigma(w)  \, d\sigma(y) \right]^2 \, d \sigma(z) \right)^\frac{1}{2},
\end{multline}
where we understand that, if $\mathrm{diam}(\partial \Omega) < \infty$, the sums are finite and terminate for $2^j r  \ge \mathrm{diam}(\partial \Omega)$.

 \begin{equation}\label{e:3estimate2} \begin{aligned}
\RNum{3} \leq &\sum_{\{j\,\mid\, 2^j \ge M\}} \left(\fint_\Delta \left[\fint_\Delta  \int_{A_j} 
\frac{r}{(2^jr)^{n}}\frac{k(w)}{a} \, d\sigma(w)  \, d\sigma(y) \right]^2 \, d \sigma(z) \right)^\frac{1}{2}\\
\lesssim& \sum_{\{j\,\mid\, M\le 2^j\le \frac{r_0}{2r}\}} \frac{r \omega(A_j)}{(2^j r)^n a} +\sum_{\{j\,\mid \,2^j\ge \frac{r_0}{2r}\}} \frac{r \omega(A_j)}{(2^j r)^n a} =\RNum{3}_a+ \RNum{3}_b.
\end{aligned}
\end{equation}
To estimate $\RNum{3}_a$ and $\RNum{3}_b$ we use \eqref{eq:avglog}, the fact that $A_j\subset\Delta_j$ (in $\RNum{3}_a$), that $\omega$ is a probability measure (in $\RNum{3}_b$) and \eqref{eq:Apratio} again with $p=1+1/(2(n-1))$. 
\begin{equation}\label{e:3a-estimate} \begin{aligned}
\RNum{3}_a =&\sum_{\{j\,\mid \,M\le 2^j\le \frac{r_0}{2r}\}} \frac{r \omega(A_j)}{(2^j r)^n a} 
\lesssim \sum_{\{j\,\mid \,M\le 2^j\le \frac{r_0}{2r}\}} \frac{r \omega(A_j)}{(2^j r)^n }\cdot\frac{\sigma(\Delta^\ast)}{\omega(\Delta^\ast)}\\
\lesssim & \sum_{\{j\,\mid \, M\le 2^j\le \frac{r_0}{2r}\}} \frac{r \sigma(\Delta^\ast)}{(2^j r)^n }\cdot\frac{\omega(\Delta_j)}{\omega(\Delta^\ast)}
\lesssim  \sum_{\{j\,\mid \, M\le 2^j\le \frac{r_0}{2r}\}} \frac{1}{M}\cdot\frac{(Mr)^n}{(2^j r)^n }\cdot\left(\frac{2^jr}{Mr}\right)^{n-1/2}\\
\lesssim &\frac{1}{\sqrt{M}}\sum_{\{j\,\mid \, M\le 2^j\le \frac{r_0}{2r}\}}2^{-j/2}=\frac{C}{M}
\end{aligned}
\end{equation}

\begin{equation}\label{e:3b-estimate} \begin{aligned}
\RNum{3}_b =&\sum_{\{j\,\mid \,2^j\ge \frac{r_0}{2r}\}} \frac{r \omega(A_j)}{(2^j r)^n a} 
\lesssim \sum_{\{j\,\mid \,\ 2^j\ge \frac{r_0}{2r}\}} \frac{r \omega(A_j)}{(2^j r)^n }\cdot\frac{\sigma(\Delta^\ast)}{\omega(\Delta^\ast)}\\
\lesssim & \sum_{\{j\,\mid \,2^j\ge \frac{r_0}{2r}\}} \frac{r }{(2^j r)^n }\cdot\frac{\sigma(\Delta^\ast)}{\omega(\Delta^\ast)}\lesssim 
\frac{r}{r_0^n}\cdot\frac{\sigma(\Delta^\ast)}{\omega(B(x_0,r_0/(5A)))}\cdot\frac{\omega(B(x_0,r_0/(5A)))}{\omega(\Delta^\ast)}\\
\lesssim &\frac{1}{M}\cdot\left(\frac{Mr}{r_0}\right)^n\cdot\frac{1}{\omega(B(x_0,r_0/(5A)))}\left(\frac{r_0}{Mr}\right)^{n-1/2}\\
\le &  \frac{C}{\omega(B(x_0, r_0/(5A)))} \cdot\left(\frac{r}{r_0}\right)^{1/2} \cdot\frac{1}{\sqrt{M}}
\end{aligned}
\end{equation}
As before the constant $C > 0$ in \eqref{e:3a-estimate} and \eqref{e:3b-estimate} depends only on $n$ and the Ahlfors regularity constant. Combining \eqref{eq25}, \eqref{e:3estimate2},
\eqref{e:3a-estimate} and \eqref{e:3b-estimate} we conclude that 
 \begin{equation}\label{e:3estimate2-tt} 
 \RNum{3} \leq  \frac{C(n,C_A)}{M} +  \frac{C(n,C_A)}{\omega(B(x_0, r_0/(5A)))} \cdot\left(\frac{r}{r_0}\right)^{1/2} \cdot\frac{1}{\sqrt{M}}.
\end{equation}

The idea to estimate $\RNum{1}$ is to approximate $\Omega$, locally, by UR domains, so that we may exploit Lemmas \ref{NTmaxbound.lem} and \ref{hmtlemma} on those approximate domains. Using the fact that { the $(x_0,M_0, r_0)$-DLTSCS holds , we may invoke Lemma \ref{URdomApprox.lem} to construct two UR `domains'  $T_Q^\pm \subseteq \om^\pm$, where $Q$ is a dyadic cube such that $\Delta(x_0, r_0/(4A)) \subset Q \subset \Delta(x_0,r_0/4)$, where the definition of $A$ above allows us to find such a cube.}   In particular,
\[ \partial T_Q^\pm \cap \Delta(x_0, r_0/(4A)) = \Delta(x_0, r_0/(4A))\]
and for $\cH^{n-1}$ a.e. $x \in \Delta(x_0, r_0/(4A))$ the unit outer normals $\nu_{T_Q^\pm}(x)$ exist and satisfy
\begin{equation}\label{TQnormagree.eq}
\nu_{T_Q^\pm}(x) = \pm \nu_{\om^+}(x).
\end{equation}
For any open set $U$ with Ahlfors regular boundary define
\[\mathcal{S}_U f(X): = \int_{\partial U} \mathcal{E}(X -y) f(y) \, d\sigma(y).\]
 In our context $U$ is either $\Omega^{\pm}$ or $ T_Q^{\pm}$.
The coincidence of $\partial T_Q^\pm \cap \Delta(x_0, r_0/(4A))$ and $\Delta(x_0, r_0/(4A))$ allows us to conclude for $f \in L^2(\Delta(x_0, r_0/(4A)))$ with $\supp f \subseteq \Delta(x_0, r_0/(4A))$,
\begin{equation}\label{SLcoincide.eq}
\mathcal{S}_{\Omega^{+}} f(X) = \mathcal{S}_{\Omega^{-}} f(X) = \mathcal{S}_{T_Q^\pm} f(X),
\end{equation}
for all $X \not \in \Delta(x_0, r_0/(4A))$.

Recall \[\RNum{1} = \left(\fint_\Delta \left| \nabla \mathcal{S} \left[\left(1 - \frac{k}{a}\right)1_{\Delta^*}\right] (z^\ast) -\fint_\Delta \nabla \mathcal{S} \left[\left(1 - \frac{k}{a}\right)1_{\Delta^*}\right](y^*) \, d\sigma(y) \right|^2 \, d \sigma(z) \right)^\frac{1}{2},\] where $z^\ast$ and $y^{\ast}$ are in non-tangential regions in $\Omega^-$ over $y,z\in \pom$. We want to dominate $\nabla \mathcal{S} \left[\left(1 - \frac{k}{a}\right)1_{\Delta^*}\right] (z^\ast)$ by a non-tangential maximal function in $T_Q^{-}$. To this end, we make the observation that if $r/r_0$ is sufficiently small (which we may ensure by adjusting the value of $\theta$) then for any $y\in \Delta$, the non-tangential cone $\Gamma^-(y) \cap B(y,r/2) \subset T_Q^-$ provided we take the constant $K$ in the definition of $T_Q^{\pm}$ large enough depending on dimension and the Ahlfors regularity of $\pom$ \footnote{This does not affect the validity of Lemma \ref{URdomApprox.lem}.}. To see this, one needs to inspect the definition of $\W_Q$ (see Appendix \ref{URapprox.sect}) and note that if $Z \in \Gamma^-(y) \cap B(y, 2r)$ then $\delta(Z) \sim |Z - y| < 2r$ and therefore $Z$ is inside a Whitney cube $I$ for $\om^-$ with 
\[ \dist(I, y) \sim \ell(I) \sim \delta(Z)<2r \lesssim \ell(Q). \]
By choosing $K$ sufficiently large, depending on allowable parameters, we can guarantee the existence of a cube $Q' \subset Q$ containing $y \in Q'$ with length $\ell(Q') \approx_K \ell(I)$. Hence $Z \in U_{Q'}^- \subset T_Q^-$. Moreover, in the construction of the Whitney region $U_{Q'}$, $\interior I^* \subset U_{Q'}$ where $I^* = (1+ \tau) I$ for some (small) parameter $\tau > 0$ (see Appendix \ref{URapprox.sect}, and note this $\tau$ is unrelated to $\tau(p)$ above). This forces $\dist(Z, \partial T_Q^-) \gtrsim_\tau \ell(I) \sim |Z- y|$ and therefore
\[ Z\in \Gamma_{\beta, T_Q^-}(y)  := \{Y \in T_Q^- : |Y - y| < (1 + \beta)\dist(Y, \partial T_Q^-)\},\]
where $\beta = \beta(n, C_A, \theta) \gg_\tau 1$. We conclude that
\begin{equation}\label{eq:contTQNTc}
	\Gamma^-(y) \cap B(y,r/2) \subset \Gamma_{\beta, T_Q^-}(y) \cap B(y,r/2).
\end{equation} 

With these observations in hand, we can estimate $\RNum{1}$.  By \eqref{eq15} and \eqref{1-k/aest.eq}, we have
\begin{multline}\label{eq23}
\RNum{1} \le  2 \left(\fint_\Delta \left| \widetilde{\mathcal{N} }
\left( \nabla \mathcal{S}_{T_Q^-} \left[\left(1 - \frac{k}{a}\right)
1_{\Delta^*}\right]\right) \right|^2 \, d \sigma \right)^{\frac{1}{2}} \\*
\leq C \, \left(\frac{\sigma(\Delta^*)}{\sigma(\Delta)}\right)^{1/2}
\left(\fint_{\Delta^\ast} 
\left| 
1 - \frac{k}{a} 
\right|^2 \, d \sigma \right)^{\frac{1}{2}}\\*
\leq \, CM^{\frac{n-1}{2}}\left(\lVert \log k  \rVert_*(B(x_0,r_0))\right)^{1/8}\le CM^{\frac{n-1}{2}} \kappa^{1/8},
\end{multline}
where $\widetilde{\mathcal{N}}$ is the non-tangential maximal function in $T^-_Q$ with aperture $\beta$ (which dominates $\mathcal{S}_{T_Q^-} \left[\left(1 - \frac{k}{a}\right)
1_{\Delta^*}\right](y^*)$ by the arguments in the preceding paragraph). {  Note that $C > 0$ above depends only on $\beta>0$, $n$, $C_A$ and the UR constants of $\pom$, which in turn only depend on $n$, $C_A$ and $\delta(X^\pm)$ .} 

Putting \eqref{eq38}, \eqref{e:3estimate2-tt} and \eqref{eq23} together we finally obtain \eqref{2.1BH.eq}.
The estimate analogous to \eqref{2.1BH.eq} when $y^*$ and $z^*$ are in $\Gamma^+(y) \cap B(y, r/2)$ and $\Gamma^+(z) \cap B(z, r/2)$ is also true by symmetry. 
It remains to use the jump relations to get an estimate on the oscillation of unit outer normal. Here we again use the approximations $T_Q^\pm$.
Applying the jump relation in Lemma \ref{hmtlemma} to $T_Q^\pm$, and using \eqref{SLcoincide.eq}, \eqref{TQnormagree.eq} and the containment  $\Gamma^\pm(y) \cap B(y,r/2) \subset \Gamma_{\beta, T_Q^\pm}(y) \cap B(y,r/2)$, we obtain for $\cH^{n-1}$ a.e. {  $y \in \Delta(x_0, r_0/(4A))$}
\begin{equation}\label{subjumploc.eq}
\nu_{\Omega^+}(y) 1_{\Delta^*}(y) = \lim_{\substack{Z \to y \\ Z \in \Gamma^+(y)}}\nabla\mathcal{S}1_{\Delta^*}(Z) - \lim_{\substack{Z \to y \\ Z \in \Gamma^-(y)}}\nabla\mathcal{S}1_{\Delta^*}(Z).
\end{equation}
Here, we need to make the further observation that the principal value singular integral operators $\mathcal{T}_{T_Q^\pm}$\footnote{The operator $\mathcal{T}_U$ is defined in the same way as $\mathcal{S}_U$.} in \eqref{eq16} and \eqref{eq30} have the property that
\[\mathcal{T}_{T_Q^+}f = \mathcal{T}_{T_Q^-}f\]
{  whenever $f \in L^2( \Delta(x_0, r_0/(4A)))$ with $\supp f \subseteq  \Delta(x_0, r_0/(4A))$.
This is a consequence of the definition of $\mathcal{T}$ and that 
\[ \partial T_Q^+ \cap B(x_0,r_0/(4A)) = \partial T_Q^- \cap B(x_0,r_0/(4A)). \] }
 Taking nontangential limits\footnote{This is justified by Lemma \ref{NTmaxbound.lem} and the dominated convergence theorem.} in \eqref{2.1BH.eq} and using \eqref{subjumploc.eq}, we obtain
\begin{eqnarray}\label{e:noKuptonow}
&&\left(\fint_{B(x,r)} \left| \nu_{\om^+}(y) -\fint_{B(x,r)} \nu_{\om^+}(z) \, d\sigma (z) \right|^{2} \, d\sigma(y) \right)^{\frac{1}{2}} \\
&&\qquad\qquad \le 
 \frac{C_1}{\hm(B(x_0, r_0/(4A)))} \cdot\left(\frac{r}{r_0}\right)^{1/2} \cdot\frac{1}{\sqrt{M}}+  C_2M^{\frac{n-1}{2}}\kappa^{\frac{1}{8}}+ \frac{C_3}{M},\nonumber
\end{eqnarray}
for $x \in \pom \cap B(x_0, r_0/(20A))$ and $0 < r \le \theta r_0$. Here, as above, {  the constants $C_1,  C_3>0$ depend on $n$ and $C_A$ and $C_2$ depends on $n$, $C_A$ and $\delta(X^\pm)$.
Notice that we may apply the same argument to $\Omega^-$ and $\log k^-$ to get an analogous estimate to \eqref{e:noKuptonow}. 

 We define a constant
\begin{equation}\label{def:C4}
	C_4 = \dfrac{C_1}{\inf_{x_0 \in \widetilde{B}\cap \pom} \omega^{\pm} (B(x_0, r_0/(5A)))}.
\end{equation}  
In fact, for each $x_0\in \widetilde{B} \cap \pom$, the harmonic measure $\omega^\pm(B(x_0, r_0/(5A))) >0$ since $\sigma \ll \omega^\pm$. Consider an arbitrary pair $x_0, x'_0 \in \widetilde{B} \cap \pom$ such that $|x_0- x'_0| < r_0/(5A)$. By the doubling property of $\omega^\pm$ (up to radius $r_0$), we have
\[ \omega^\pm (B(x_0,r_0/(5A))) \leq \omega^\pm (B(x'_0, r_0)) \leq C \omega^\pm (B(x'_0, r_0/(5A))). \]
Since $\widetilde{B}\cap \pom$ is compact, it can be covered by finitely many balls centered on $\widetilde{B} \cap\pom$ with radii $r_0/(5A)$. In particular the denominator in \eqref{def:C4} is a strictly positively constant depending on the domains $\Omega^\pm$ and $\widetilde{B}$, and thus the constant $C_4$ is well-defined. Notice that the same argument applied to $\log k^-$ combined with \eqref{e:noKuptonow} and \eqref{def:C4} yields:
}
\begin{eqnarray}\label{tt-7}
&&\left(\fint_{B(x,r)} \left| \nu_{\om^\pm}(y) -\fint_{B(x,r)} \nu_{\om^\pm}(z) \, d\sigma (z) \right|^{2} \, d\sigma(y) \right)^{\frac{1}{2}} \\
&&\qquad\qquad\le 
 {C_4}\left(\frac{r}{r_0}\right)^{1/2} \cdot\frac{1}{\sqrt{M}}+  C_2M^{\frac{n-1}{2}}\kappa^{\frac{1}{8}}+ \frac{C_3}{M},\nonumber
\end{eqnarray}
where $C_4=C_4(n, C_A, \widetilde{B}, \Omega^\pm)$. {  For $\epsilon>0$ sufficiently small (satisfying $C_3 \epsilon  \leq 4$), we choose the constant $M$ such that $\frac{1}{\sqrt{M}}  = \frac{\epsilon}{4}$ and $\frac{C_3}{\sqrt{M}}\le 1$; we also choose the constant $\theta$ such that $M\theta < 1/(10A)$ and $C_4 \theta^{1/2} \leq 1$.} Then \eqref{tt-7} becomes
\begin{equation}\label{tt-8}
\left(\fint_{B(x,r)} \left| \nu_{\om^\pm}(y) -\fint_{B(x,r)} \nu_{\om^\pm}(z) \, d\sigma (z) \right|^{2} \, d\sigma(y) \right)^{\frac{1}{2}} 
\le 
\frac{\epsilon}{2} +  C_5{\epsilon^{-(n-1)}}\kappa^{\frac{1}{8}},
\end{equation}
where $C_5$ depends on $n$ and $C_A$. Note that in the above estimate, only $\theta$ depends on $\widetilde{B}$.
Thus, perhaps further shrinking $\kappa_1$ (depending on $\epsilon$, $n$, $C_A$ and $\delta(X^\pm)$ and independent of $\widetilde{B}$), \eqref{tt-8} becomes
\begin{eqnarray}\label{tt-9}
&&\left(\fint_{B(x,r)} \left| \nu_{\om^\pm}(y) -\fint_{B(x,r)} \nu_{\om^\pm}(z) \, d\sigma (z) \right|^{2} \, d\sigma(y) \right)^{\frac{1}{2}} \\
&&\qquad\qquad\qquad\le 
\frac{\epsilon}{2} +  C_5(n, C_A) \epsilon^{-(n-1)}\kappa_1^{\frac{1}{8}}\le \epsilon.\nonumber
\end{eqnarray}

To sum up, we have shown that given $\epsilon>0$ there exists a small constant $\kappa_1$ depending on $\epsilon, n, C_A$ and $\delta(X^\pm)$ such that the following holds: For every ball $B^*$ centered on the boundary with radius less than $(1/4)\diam(\pom)$, if there is a radius $r_0 = r_0(B^*)$ such that 
\begin{equation}\label{tt-10}
\sup_{x_0\in B^* \cap \pom} \|\log k^\pm\|_*(B(x_0, r_0)) \leq \kappa \leq \kappa_1,
\end{equation}
then we can find $\theta\in (0,1)$ depending on $n, C_A$, the domains $\Omega^\pm$ and $\widetilde{B} := \tfrac{1}{4}\overline{B^*}$ so that
\begin{equation}\label{tt-11}
\sup_{x_0 \in \widetilde{B} \cap  \in\pom} \|\nu\|_*(B(x_0, \theta r_0)) \leq \epsilon.
\end{equation}
Thus $\nu\in BMO_{\loc}( \sigma)$ with 
constant at most $\epsilon$ (see Remark \ref{comactvsballs.rmk}). This concludes the proof of Theorem \ref{t:osckcontrolsoscn}. 
\end{proof}

 \subsection{Free Boundary Results}

In this section we combine Theorem \ref{t:osckcontrolsoscn} with Corollaries \ref{c:BMO} and \ref{c:VMO} to obtain information about the local geometry 
of a domain (with minimal hypothesis) from the local oscillation of the logarithm of the interior and exterior Poisson kernels.

\begin{theorem*}[Theorem \ref{Tchar2sideCAD.thrm}]
Let $n \ge 3$ and suppose $\om^+ \subset \rn$ and $\om^-= \rn\setminus \overline{\om^+}$ are domains satisfying $\partial\Omega:= \pom^+ = \pom^-$,
and that $\pom$ is $(n-1)$-Ahlfors regular. Then the following are equivalent:
\begin{itemize}
\item[\emph{(i)}] $\om^\pm$ are both vanishing chord-arc domains (see Definition \ref{d:CADSC})
\item[\emph{(ii)}] There exists $X^+ \in \om^+$ and $X^- \in \om^-$ such that $k^+ = \frac{d\hm_+^{X^+}}{d\sigma}$ and $k^- = \frac{d\hm_-^{X^-}}{d\sigma}$ exist and $\log k^\pm \in VMO_{\loc}(d\sigma)$.
\end{itemize}
\end{theorem*}

\begin{proof}[Proof of Theorem \ref{Tchar2sideCAD.thrm}]
(i) implies (ii) is the main theorem in \cite{kenig2003poisson}. That (ii) implies (i) follows from Theorem \ref{t:osckcontrolsoscn}. Indeed, by Corollary \ref{c:VMO}, to show that $\om^{\pm}$ are vanishing chord-arc domains it suffices to prove that $\nu \in VMO_{\loc}(d\sigma)$. Theorem \ref{t:osckcontrolsoscn} asserts that this is the case when $\log k^\pm \in VMO_{\loc}(d\sigma)$.
\end{proof}

The following is a quantified version of Theorem \ref{Tchar2sideCAD.thrm} which results from the remark at the end of the proof of Theorem \ref{t:osckcontrolsoscn}.

\begin{theorem}[Quantified version of Theorem \ref{Tchar2sideCAD.thrm}]\label{FBthrm.thrm}
Let $\om^+ \subset \rn$ and $\om^-= \rn\setminus \overline{\om^+}$ be domains with common (topological) boundary $\pom = \pom^+ = \pom^-$. Assume that $\pom$ is $(n-1)$-Ahlfors regular and let $X^\pm \in \om^\pm$  be such that $k^\pm = \frac{d\hm_\pm^{X^\pm}}{d\sigma}$ exist. Given $\delta > 0$ there exists 
$\kappa= 
\kappa (\delta,n, C_A,\delta(X^\pm)) > 0$ such that if $\log k^\pm\in BMO_{\loc}(\sigma)$ with constant less than $\kappa$, then
 $\om^+$ and $\om^-$ are $\delta$-chord-arc domains. 

Conversely, for every $\kappa> 0$  there exists $\delta = \delta(\eta, n, C_A) > 0$ if $\nu\in BMO_{\loc}(\sigma)$ with constant less than $\delta$, then
$\log k^\pm\in BMO_{\loc}(\sigma)$ with  constant less than $\kappa$.
\end{theorem}

\begin{proof}
This is a combination of Theorem \ref{t:osckcontrolsoscn}, Corollary \ref{c:BMO} and the work in \cite{kenig1999free}.
\end{proof}

\appendix

\section{Proof of Theorem \ref{t:SL}}\label{s:proofofseparation}

In this section we prove Theorem \ref{t:SL}; that small excess implies flatness in the sense of Reifenberg. This is a corollary of the height bound, Theorem \ref{t:HB}. Many of the techniques, included for completeness, are standard. Another consequence of Theorem \ref{t:HB} is a Lipschitz Approximation theorem, Theorem \ref{t:LA}, which is proven at the end of this section. It is of independent interest and is not used in this paper.

The next lemma is contained in \cite[Lemma 22.11]{maggi2012sets}. We recall some notation introduced in other sections. We define $q(x) = \langle x , e_{n} \rangle$, $p(x) = x - q(x) e_{n}$, $C_{r} =\{|q(x)| < r \} \cap \{|p(x)| < r \}$, $D_{r} =p(C_{r})$ and $D = p(C_{1})$. We consider $D, D_{r}$ to be subsets of $\R^{n-1}$. Finally, when the set $E$ is clear from context, recall $e_{n}(x,r) = e(E,x, r, e_{n})$ and if $x = 0$, $e_{n}(r) = e(E,0,r,e_{n})$.

\begin{lemma}[Excess Measure] \label{l:EM}
If $E \subset \R^{n}$ is a set of locally finite perimeter in $\R^{n}$ with $0 \in \partial E$, such that for some $t_{0} \in (0,1)$
\eqref{e:15}, \eqref{e:16}, and \eqref{e:17} are each satisfied with $r = 1$ and $\nu = e_{n}$, then writing $M  = C_{1} \cap \partial^{*} E$ it follows that for any Borel $G \subset D$,

\begin{equation} \label{e:24}
\cH^{n-1}(G) = \int_{M \cap p^{-1}(G)} \langle \nu_{E}, e_{n} \rangle d \cH^{n-1}.
\end{equation}
Moreover, for every $\varphi \in C^{0}_{c}(D)$ and $t \in (-1,1)$
\begin{equation} \label{e:25}
\int_{D} \varphi dx = \int_{M } \varphi(p(x)) \langle \nu_{E}(x), e_{n} \rangle \dif \cH^{n-1}
\end{equation}
and
\begin{equation} \label{e:26}
\int_{E_{t} \cap D} \varphi dx = \int_{M \cap \{q(x) > t \}} \varphi(p(x)) \langle \nu_{E}(x), e_{n} \rangle \dif \cH^{n-1} \qquad \forall t \in (-1,1)
\end{equation}
where $E_{t} = \{ z \in \R^{n-1} \mid (z,t) \in E \}$. In fact, the set function
\begin{equation} \label{e:27}
\zeta(G) = \cH^{n-1}(M \cap p^{-1}(G)) - \cH^{n-1}(G)
\end{equation}
defines a Radon measure in $D$, and is called the excess measure of $E$ over $D$ since $\zeta(D) = e(E,0,1,e_{n})$.
\end{lemma}

\begin{theorem}[Height bound: {compare with \cite[Theorem 22.8]{maggi2012sets}}] \label{t:HB}
Given $C_{A} \ge 1, r_{0} >0$, and $n \ge 2$, there exist constants $\epsilon_1= \epsilon(n,C_A)>0$ and $C_{1} = C(n,C_{A}) \ge 1$ such that if $E \subset \R^{n}$ is Ahlfors regular with constant $C_{A}$ up to scale $4r_{0}$ and $x_{0} \in \partial E$ satisfies
\begin{equation}
	e_n(x_0,4r_0)\leq \epsilon_1,
\end{equation}
then
\begin{equation} \label{e:29}
\frac{1}{r_{0}} \sup \{ | q(x_{0}) - q(y)| : y \in C(x_{0}, r_{0}, e_{n}) \cap \partial E \} \le C_{1} e_{n}(x_{0},4 r_{0})^{\frac{1}{2(n-1)}}.
\end{equation}
\end{theorem}

\begin{proof}  By Remark \ref{r:excessunderoperations} we let $x_0 = 0$ and $2r_0 = 1$. We then want to show that $|q(x)| \le c_{0}(n) e_{n}(2)^{\frac{1}{2(n-1)}}$ whenever $x \in C_{1/2} \cap \partial E$. 

We first assume that $\epsilon_{1} \le \min \left\{ \omega(n, \frac{1}{4}, C_{A}), 2^{-n} \cH^{n-1}(D)\right\}$, with $\omega(n, \frac{1}{4}, C_{A})$ from Lemma \ref{l:SL}. Then, by Lemma \ref{l:SL}, $|q(x)| \le \frac{1}{4}$ whenever $x \in C_{1} \cap \partial^* E=:M$, and moreover $E$ satisfies the hypotheses of Lemma \ref{l:EM} with $t_0 = \frac14$. Therefore
\begin{equation} \label{e:30}
0 \le \cH^{n-1}( M) - \cH^{n-1} (D) \le e_{n}(1) \le 2^{n-1} e_{n}(2)
\end{equation}
and
\begin{equation} \label{e:31}
0 \le \cH^{n-1}(M \cap \{ q(x) > t \}) - \cH^{n-1}(E_{t} \cap D) \le 2^{n-1} e_{n}(2) \quad \forall t \in (-1,1).
\end{equation}

Now, we consider $f : (-1,1) \to [0, \cH^{n-1}(M)]$ defined by
\begin{equation} \label{e:32}
f(t) = \cH^{n-1}(M \cap \{ q(x) > t \} ). 
\end{equation} 
By Lemma \ref{l:SL} 
\begin{equation} \label{e:33}
f(t) = \begin{cases}
\cH^{n-1}(M) & -1 < t < -1/4 \\
0 & 1/4 < t < 1.
\end{cases}
\end{equation}

Since $f$ is decreasing and right-continuous there exists $|t_{0}| < \frac{1}{4}$ such that
\begin{equation} \label{e:34}
\begin{cases}
f(t) \le \frac{ \cH^{n-1}(M)}{2} & t \ge t_{0} \\ f(t) > \frac{\cH^{n-1}(M)}{2} & t < t_{0} .
\end{cases}
\end{equation}

\textbf{Claim:} If $x \in C_{1/2} \cap \partial E$ then $|q(x) - t_{0}| \le c(n) e_{n}(2)^{\frac{1}{2(n-1)}}$. In particular, since $0 \in \partial E$, this ensures $|t_{0}| \le c(n)e_{n}(2)^{\frac{1}{2(n-1)}}$.

The claim will be verified by showing that $q(x) - t_{0} \le c(n)e_{n}(2)^{\frac{1}{2(n-1)}}$, then considering $\R^{n}\setminus E$ to get $|q(x) - t_{0}| \le c(n)e_{n}(2)^{\frac{1}{2(n-1)}}$. Since $\partial E = \spt \mu_E = \overline{\partial^* E}$ and the projection function $q$ is continuous, it suffices to prove the estimate for $x\in C_{1/2}\cap \partial^* E$. To bound $q(x) - t_{0}$, we first show there exists $t_{1}$ with $q(x) - t_{1} \le c(n) e_{n}(2)^{\frac{1}{2(n-1)}}$ and then that $t_{1} - t_{0}$ satisfies a similar upper-bound.

By choice of $\epsilon_{1}$, $\sqrt{e_{n}(2)} < \frac{1}{2C_{A}} \le \frac{\cH^{n-1}(M)}{2}$. So, we choose $t_{1} \in (t_{0}, \frac{1}{4})$ such that
\begin{equation} \label{e:35}
\begin{cases}
f(t) \le \sqrt{e_{n}(2)} & \forall t \ge t_{1} \\
f(t) > \sqrt{e_{n}(2)} & \forall t < t_{1}.
\end{cases}
\end{equation}

To see $q(x) - t_{1} \le c(n) e_{n}(2)^{\frac{1}{2(n-1)}}$ for all $x \in C_{1/2} \cap \partial^* E$, note if $y \in C_{1/2} \cap \partial^* E$ and $q(y) > t_{1}$, then $q(y) - t_{1} < \frac{1}{2}$ since $t_{1} \in (t_{0}, 1/4)$ and $|q(y)| < \frac{1}{4}$. In particular, $(q(y)-t_{1})$ is a small enough scale for Ahlfors-regularity to hold. Hence,
\begin{equation} \label{e:36}
C_{A}^{-1}(q(y) - t_{1})^{n-1} \le |\mu_{E}|(B(y,q(y) - t_{1})).
\end{equation}
Since $x \in B(y, q(y) - t_{1})$ implies $q(y) - q(x) \le |x-y| < q(y) - t_{1}$ and since $y \in C_{1/2}$ with $q(y) - t_{1} < \frac{1}{2}$,
\begin{equation} \label{e:37}
B(y, q(y) - t_{1}) \subset \{ x \in C_{1} \mid q(x) > t_{1} \}.
\end{equation} 
Thus $B(y, q(y) - t_{1})\cap \partial^* E \subset M\cap \{q > t_1\}$. So, \eqref{e:36} and \eqref{e:37} imply
\begin{equation} \label{e:38}
C_{A}^{-1}(q(y) - t_{1})^{n-1} \le |\mu_{E}| ( C_{1} \cap \{ q(x) > t_{1} \} ) = \cH^{n-1}(M \cap \{ q(x) > t_{1} \}) = f(t_{1}).
\end{equation}
By the choice of $t_{1}$ in \eqref{e:35}, it follows that under the standing assumption $q(y) - t_{1} > 0$ we have 
\begin{equation} \label{e:39}
q(y) - t_{1} \le c(n,C_{A}) e_{n}(2)^{\frac{1}{2(n-1)}},
\end{equation}
as desired. Note, \eqref{e:39} is trivially true when $q(y) \le t_{1}$.

Next we show that $t_{1} - t_{0} \le c_{n}e_{n}(2)^{\frac{1}{2(n-1)}}$, which verifies the Claim 1. We will use a slicing result, see \cite[Theorem 18.11]{maggi2012sets} which ensures that for almost every $t \in (-1,1)$,
\begin{equation} \label{e:40}
\cH^{n-2}( (\partial^{*} E_{t}) ~ \Delta ~ (\partial^{*} E)_{t} ) = 0,
\end{equation}
where $(\partial^{*}E)_{t} = \{ z \in \R^{n-1} : (z,t) \in \partial^{*} E\} \subset \R^{n-1}$ and $E_{t} = \{ z \in \R^{n-1} \mid (z,t) \in E \} \subset \R^{n-1}$. Furthermore, the co-area formula ensures that for any $g : \R^{n} \to [0, \infty]$ a non-negative Borel function,
\begin{equation} \label{e:caf}
\int_{\partial^{*}E} g \sqrt{1 - \langle \nu_{E}, e_{n} \rangle^{2}} d \cH^{n-1}= \int_{\R} \left( \int_{(\partial^{*}E)_{t}} g \dif \cH^{n-2} \right) \dif t.
\end{equation}
In particular, realizing the square-root term on the left is just the Jacobian of the projection $p$, and choosing the function $g = \chi_{C_{1}}$, recalling that $C_{1} \cap \partial^{*}E \supset M$ is Ahlfors regular up to scale $2$,

\begin{align*}
\int_{-1}^{1} \cH^{n-2} \left( (\partial^{*}E)_{t} \cap D \right) dt 
& = \int_{M} \sqrt{1 - \langle \nu_{E}, e_{n} \rangle^{2}} \dif \cH^{n-1} \\
& \le \left( 2 \cH^{n-1}(M) \right)^{\frac{1}{2}} \left( \int_{M} (1 - \langle \nu_{E}, e_{n} \rangle ) \dif \cH^{n-1} \right)^{\frac{1}{2}} \\
& \le c(n,C_{A}) \sqrt{e_{n}(2)},
\end{align*}
 We extract from the above that
\begin{equation} \label{e:41}
\int_{t_{0}}^{1} \cH^{n-2}(\partial^{*}E_{t} \cap D) \dif t \le \int_{-1}^{1} \cH^{n-2}(\partial^{*}E_{t} \cap D) dt \le c(n) \sqrt{e_{n}(2)}.
\end{equation}

For almost all $t \in [t_{0}, 1)$ it follows from $\cH^{n-1}(E_{t} \cap D) \le \cH^{n-1}(M \cap \{q(x) > t \})$,  \eqref{e:30}, \eqref{e:31}, and \eqref{e:34} that
\begin{equation*}
\cH^{n-1}(E_{t} \cap D) \le \frac{\cH^{n-1}(M)}{2} \le \frac{\cH^{n-1}(D)}{2} + 2^{n-2} e_{n}(2) \le \frac{3}{4} \cH^{n-1}(D)
\end{equation*}
where we used that $e_n(2) \le 2^{-n} \cH^{n-1}(D)$.

Applying the relative isoperimetric inequality (see \cite[(12.45)]{maggi2012sets}) in $\R^{n-1}$ to the set $E_{t} \cap D$ we have
\begin{equation} \label{e:42}
\cH^{n-2}(D \cap \partial^{*} E_{t} ) \ge c(n) \cH^{n-1}(E_{t} \cap D)^{\frac{n-2}{n-1}} \quad \text{for a.e. } t \in [ t_{0},1 )
\end{equation}
 \eqref{e:41} and \eqref{e:42} together imply (where the constant $c(n)$ can change in every instance, but only depends on $n$)
\begin{equation} \label{e:43}
\int_{t_{0}}^{t_{1}} \cH^{n-1} (E_{t} \cap D)^{\frac{n-2}{n-1}} dt \le c(n) \int_{t_{0}}^{1} \cH^{n-1}(E_{t} \cap D)^{\frac{n-2}{n-1}} dt \le c(n) \sqrt{e_{n}(2)}.
\end{equation}

Finally, \eqref{e:31} and \eqref{e:35} yield for $t < t_{1}$,
\begin{align*}
\cH^{n-1}(E_{t} \cap D) &\ge \cH^{n-1}(M \cap \{q(x) > t \}) - 2^{n-1} e_{n}(2)  \\
& \ge \sqrt{e_{n}(2)} - 2^{n-1}e_{n}(2) \ge c(n) \sqrt{e_{n}(2)},
\end{align*}
which combined with \eqref{e:43} ensures
$$
(t_{1} - t_{0}) e_{n}(2)^{ \frac{n-1}{2(n-1)} - \frac{1}{2(n-1)}} = (t_{1} - t_{0}) \sqrt{e_{n}(2)}^{\frac{n-2}{n-1}} \le c(n) \sqrt{e_{n}(2)},
$$
so that $t_{1} - t_{0} \le c(n) e_{n}(2)^{\frac{1}{2(n-1)}}$ as desired.
\end{proof}

We are now ready to prove Theorem \ref{t:SL} which first appears in Section \ref{s:excessdecay} above. We restate it here for convenience: 

 \begin{theorem*}
Fix $C_{A} \ge 1, r_{0} > 0$, and $n \ge 2$. Let $\epsilon_1 = \epsilon(C_A, n) > 0$ be as in Theorem \ref{t:HB}. If $E \in \cA(C_{A}, 4r_{0})$ and $x_{0} \in \partial E$ satisfies 
\begin{equation} \label{e:CSL-1}
e(E,x_{0}, 2r, \nu) \le \epsilon_{1} 
\end{equation}
for some $\nu \in \St^{n}$ and $0 < r < 2r_{0}$ then 
\begin{equation} \label{e:ssep1b}
\{ x \in C(x_{0}, r, \nu) \cap E \mid \langle x - x_{0}, \nu \rangle > r C_{1} e(E,x_{0}, 2r, \nu)^{\frac{1}{2(n-1)}} \} = \emptyset
\end{equation}
and
\begin{equation} \label{e:ssep2b}
\{ x \in C(x_{0}, r, \nu) \cap E^{c} \mid \langle x - x_{0}, \nu \rangle < - r C_{1} e(E,x_{0}, 2r, \nu)^{\frac{1}{2(n-1)}} \} = \emptyset.
\end{equation}
\end{theorem*}

\begin{proof}[Proof of Theorem \ref{t:SL}]
We will verify \eqref{e:ssep1b}, and \eqref{e:ssep2b} follows similarly. By translation and rotation, without loss of generality we suppose $x_{0} = 0$ and $\nu = e_{n}$. 

Suppose \eqref{e:ssep1b} fails. Then, there exists $x \in C_{r} \cap E$ with $q(x) > r C_{1} e_{n}(2r)^{\frac{1}{2(n-1)}}$.  However, $\epsilon_{1} \le \omega(n, 1/4, C_{A})$ guarantees that \eqref{e:16} holds with $t_{0} = \frac{1}{4}$. However, \eqref{e:16} guarantees that there exists some $y \in C_{r} \cap E^{c}$ with $q(x) < q(y) < r$. But then, there exists $z \in \partial E $ which lies on the line segment connecting $x$ and $y$. In particular, $q(z) > q(x) > r C_{1} e_{n}(2r)^{\frac{1}{2(n-1)}}$ contradicting Theorem \ref{t:HB}.
\end{proof}

The following theorem is another consequence of the height bound, Theorem \ref{t:HB}. Hereafter, $\nabla^{\prime}$ denotes the gradient in $\R^{n-1}$.

\begin{theorem}[Lipschitz function approximation: {compare with \cite[Theorem 23.7]{maggi2012sets}}] \label{t:LA}
There exist positive $C_{3} = C(n,C_{A}),$ $\epsilon_{3} = \epsilon(n,C_{A})$, $\delta_{0} = \delta(n,C_{A})$, and $L = L(n,C_{A}) < 1$ with the following properties. If $E \in \cA(C_{A}, 13r)$ and $e_{n}(x_{0}, 13r) \le \epsilon_{3}$ with $x_{0} \in \partial E$, then for
$M = C(x_{0}, r) \cap \partial E$ and $M_{0} = \{ y \in M \mid \sup_{0 < s < 8r} e_{n}(y,s) < \delta_{0} \}$ there exists $u : \R^{n-1} \to \R$ with $\lip(u) \le L$ and
\begin{equation} \label{e:47}
\sup_{\R^{n-1}} \frac{|u|}{r} \le C_{3} e_{n}(x_{0}, 13r)^{ \frac{1}{2(n-1)} }
\end{equation}
such that $M_{0} \subset M \cap \Gamma$ where 
\begin{equation} \label{e:48}
\Gamma = x_{0} + \{(z,u(z)) \mid z \in D_{r}\}.
\end{equation}
Furthermore,
\begin{equation} \label{e:49}
\frac{ \cH^{n-1}(M \Delta \Gamma)}{r^{n-1}} \le C_{3} e_{n}(x_{0},13r),
\end{equation}

\begin{equation} \label{e:50}
\frac{1}{r^{n-1}} \int_{D_{r}} | \nabla^{\prime} u|^{2} \le C_{3} e_{n}(x_{0}, 13r),
\end{equation}

and 
\begin{equation} \label{e:49.1}
\dist(x, (p(x),u(p(x)))) = |q(x) - u(p(x))| \le 2 L \dist(p(x), p(M_{0})) \quad \forall x \in M.
\end{equation}

In fact, \eqref{e:49.1} ensures there exist Lipschitz functions $u_{\pm}$ defined by
\begin{equation} \label{e:up}
u_{+}(x) = 
\begin{cases}
u(x) & x \in p(M_{0}) \\
\inf_{y \in p(M_{0})} u(y) + L |x - y| & x \in D \setminus p(M_{0})
\end{cases} 
\end{equation}

\begin{equation} \label{e:um}
u_{-}(x) = \begin{cases}
u(x) & x \in p(M_{0}) \\
\sup_{y \in p(M_{0})} u(y) - L | x - y|  & x \in D \setminus p(M_{0})
\end{cases}
\end{equation}
with the property that

\begin{equation} \label{e:LipEnv}
u_{-}(p(x)) \le q(x) \le u_{+}(p(x)) \qquad \forall x \in M.
\end{equation}
\end{theorem}

\begin{proof} 

Step 1: Up to replacing $E$ with $E_{x_{0},r}$ and correspondingly replacing $u$ with $u_{r}(z) = r^{-1} u(rz)$, we can reduce to proving that if $E \in \cA(C_{A}, 13)$ with $0 \in \partial E$, if 
\begin{equation} \label{e:51}
M = C \cap \partial E, \qquad M_{0} = \{ y \in M \mid \sup_{0 < s < 8} e_{n}(y,s) < \delta_{0}(n,C_{A}) \},
\end{equation}
and if $e_{n}(0,13) \le \epsilon_{3}$ then there exists a Lipschitz function $u : \R^{n-1} \to \R$ with $\lip(u) \le L < 1$ such that
\begin{equation} \label{e:52}
\sup_{\R^{n-1}} |u| \le C_{3} e_{n}(0,13)^{{ \frac{1}{2(n-1)} }}
\end{equation}
such that $M_{0} \subset M \cap \Gamma$ where
\begin{equation} \label{e:53}
\Gamma = \{(z, u(z)) \mid z \in D \}.
\end{equation}
Furthermore,
\begin{equation} \label{e:54}
\cH^{n-1}(M \Delta \Gamma) \le C_{3} e_{n}(0,13)
\end{equation}
and
\begin{equation} \label{e:55}
\int_{D} | \nabla^{\prime}u|^{2} \le C_{3} e_{n}(0,13).
\end{equation}

By Theorem \ref{t:HB} it follows that
\begin{equation} \label{e:56}
\sup \left\{ |q(x)| \mid x \in C_{2} \cap \partial E \right\} \le C_{1} e_{n}(0,13)^{\frac{1}{2(n-1)}}.
\end{equation}
By choosing $\epsilon_{3} \le \epsilon_{1} \le \omega(n,\frac{1}{4}, C_{A})$, $E$ satisfies the hypotheses of Lemma \ref{l:SL}. Consequently, Lemma \ref{l:EM} and \eqref{e:compscales} imply,
\begin{equation} \label{e:57}
0 \le \cH^{n-1}(M \cap p^{-1}(G)) - \cH^{n-1}(G) \le e_{n}(0,1) \le 13^{n-1}e_{n}(0,13),
\end{equation}
for every Borel set $G \subset D$. Meanwhile, Theorem \ref{t:SL} ensures
\begin{equation} \label{e:58}
\left\{ x \in C_{2} \mid q(x) < - \frac{1}{4} \right\} \subset C_{2} \cap E \subset \left\{ x \in C_{2} \mid q(x) < \frac{1}{4} \right\}.
\end{equation}

Step 2: We show that $M_{0}$ is contained in the graph of a Lipschitz function $u$, satisfying \eqref{e:52}  and \eqref{e:54}. In order to create the Lipschitz function, we first need to know $M_{0}$ is non-empty. This follows from a covering argument done later in more detail in \eqref{e:nonempty}.

Define $\| \cdot \| = \max \{ |p( \cdot)|, |q(\cdot)| \}$. Then, $C(y,s) = \{ z \in \R^{n} \mid \|z - y \| < s \}$. For fixed $y \in M_{0}$ and $x \in M$ and consider $F = E_{y, \|x-y\|}$.  Notably, $\|x-y\| < 2$. Since $y \in M_{0}$ and $4 \| x - y \| < 8$ it follows from \eqref{e:scalex} and \eqref{e:51} that
$$
e_{n}(F,0,4) = e_{n}(E,y, 4 \|x-y\|) \le \delta_{0}.
$$
So, choosing $\delta_{0} \le \epsilon_{1}$ allows us to apply Theorem \ref{t:HB} to $F \in \cA(C_{A},4)$ and conclude that
\begin{equation} \label{e:59}
\sup \{ |q(w)| \mid w \in C \cap \partial F \} \le  C_{1} e_{n}(F,0,4)^{\frac{1}{2(n-1)}} \le C_{1} \delta_{0}^{\frac{1}{2(n-1)}}.
\end{equation}
Applying this height-bound to the specific point $w = \frac{x-y}{\|x-y\|}$ we find
\begin{equation} \label{e:60}
|q(x) - q(y)| \le C_{0}(n) \delta_{0}^{\frac{1}{2(n-1)}}\|y-x\|.
\end{equation}
If we now define $L = C_{1} \delta_{0}^{\frac{1}{2(n-1)}}$ and choose $\delta_{0}$ so small that $L< 1$ it follows from \eqref{e:60} that $|q(x) - q(y)| < \|x-y\|$ which ensures $\|x-y\| = |p(x) - p(y)|$, and hence \eqref{e:60} can be written
\begin{equation} \label{e:61}
|q(x) - q(y)| \le L  |p(x) - p(y)|, \qquad \forall y \in M_{0}, x \in M,
\end{equation}
which implies that $p|_{M_{0}}$ is invertible. Define $u : p(M_{0}) \to \R$ such that $u(p(x)) = q(x)$ for every $x \in M_{0}$. Evidently, \eqref{e:61} ensures $u$ satisfies
\begin{equation} \label{e:62}
|u(p(x))- u(p(y))| \le L |p(x) - p(y)|, \qquad \forall x,y \in M_{0}.
\end{equation}

Since $M_{0} \subset M$, it follows from \eqref{e:56} that
\begin{equation} \label{e:63}
|u(p(x))| = |q(x)| \le C_{1} e_{n}(0,13)^{\frac{1}{2(n-1)}}, \qquad \forall x \in M_{0}.
\end{equation}

Via Kirzbraun's theorem and trunction we extend $u$ from $p(M_{0})$ to $\R^{n-1}$ with Lipschitz constant $L < 1$ such that the $L^{\infty}$-bound from \eqref{e:63} holds on all of $\R^{n-1}$, which verifies \eqref{e:52}. The definition of $u$ on $p(M_{0})$ guarantees $M_{0} \subset M \cap \Gamma$ where $\Gamma$ is as in \eqref{e:53}.

Next we show \eqref{e:54}. By definition of $M_{0}$, for every $y \in M \setminus M_{0}$ there exists $s_{y} \in (0,8)$ with
\begin{equation} \label{e:64}
\delta_{0} s_{y}^{n-1} < \int_{C(y,s_{y}) \cap \partial E} \frac{ |\nu_{E} - e_{n}|^{2}}{2} \dif \cH^{n-1}.
\end{equation}

Let $\cF$ be the set of all balls $B(y_{k},\sqrt{2} s_{k})$ centered on $M \setminus M_{0}$ satisfying \eqref{e:64} of radius at most $8 \sqrt{2}$. Each ball is contained in $C_{1 + 8 \sqrt{2}} \subset C_{13}$.  By Besicovitch's covering theorem (see \cite[Theorem 2, Seciton 1.5.2]{evans1992measure}) we partition $\cF$ into $N_{n}$ disjoint families of balls $\cG_{j}$. Then, there exists $j$ such that
\begin{align*}
\cH^{n-1}(M \setminus M_{0}) &\le N_{n} \sum_{B(y_{k}, s_{k}) \in \cG_{j}} \cH^{n-1} \left( (M \setminus M_{0} ) \cap B(y_{k}, \sqrt{2} s_{k})\right) \\
& \le N_{n} \sum_{k \in \N} \cH^{n-1} \left( M \cap B(y_{k}, \sqrt{2}s_{k}) \right) \\
& \le N_{n} C_{A} 2^{\frac{n-1}{2}} \sum_{k \in \N} s_{k}^{n-1}.
\end{align*}
Since $C(y_{k}, s_{k}, e_{n}) \subset B(y_{k}, \sqrt{2} s_{k})$ the family of cylinders are also mutually disjoint. So, \eqref{e:64} combined with the preceding computation yields
\begin{align} \nonumber 
\cH^{n-1}(M \setminus M_{0}) &\le C \sum_{k \in \N} s_{k}^{n-1} \\
\nonumber & \le \frac{C}{\delta_{0}} \sum_{k} \int_{C(y_{k}, s_{k})} \frac{|\nu_{E} - e_{n}|^{2}}{2} \dif \cH^{n-1} \\
\label{e:nonempty} &  \le \frac{C}{\delta_{0}} e_{n}(0,13).
\end{align}
Keeping in mind that $\delta_{0} < \min\{C_{1}^{-2(n-1)}, \epsilon_{1}\}$, if $\epsilon_{3}$ is small enough that $\delta_{0} \ge \frac{C \epsilon_{3}}{\cH^{n-1}(D)}$ it follows that $M_{0}$ is non-empty.  This also adds an additional constraint on $\epsilon_{3}$.
A consequence of \eqref{e:nonempty} and $M \setminus \Gamma \subset M \setminus M_{0}$ is 
\begin{equation} \label{e:65}
\cH^{n-1}(M \setminus \Gamma) \le C e_{n}(0,13).
\end{equation} 
To finish verifying \eqref{e:54} it remains to bound $\cH^{n-1}(\Gamma \setminus M)$.

Indeed, $\lip(u) \le 1$ and $M_{0} \subset \Gamma$ together ensure

$$
\cH^{n-1}(\Gamma \setminus M) \le \sqrt{1 + |\nabla^{\prime}u|^{2}} \cH^{n-1}(p(\Gamma \setminus M)) \le \sqrt{2} \cH^{n-1} \left( M \cap p^{-1} \left(p(\Gamma \setminus M) \right) \right).
$$
But, $M \cap p^{-1} \left(p(\Gamma \setminus M) \right) \subset M \setminus \Gamma$, so by the bound in \eqref{e:65}, we have the necessary bound on $\cH^{n-1}(\Gamma \setminus M)$, verifying \eqref{e:54} with a constant we denote as $C_{3}$.

Step 3: We verify \eqref{e:55}.

The first necessary observation is to note that for almost every $x \in M \cap \Gamma$,
\begin{equation} \label{e:66}
\nu_{E}(x) = \lambda(x) \frac{ ( - \nabla^{\prime}u(p(x)), 1 )}{\sqrt{1 + | \nabla^{\prime}u(p(x))|^{2}}}
\end{equation}
where $\lambda(x) \in \{-1, 1\}$. Since $| \nu_{E} - e_{n}|^{2} = |p(\nu_{E})|^{2}$, \eqref{e:66} implies
\begin{align*}
e_{n}(0,1) & \ge \frac{1}{2} \int_{M \cap \Gamma} |p(\nu_{E})|^{2} \dif \cH^{n-1} \\
& = \frac{1}{2} \int_{M \cap \Gamma} \frac{ | \nabla^{\prime}u(p(x))|^{2}}{1 + | \nabla^{\prime}u(p(x))|^{2}} \dif \cH^{n-1}(x) \\
& = \frac{1}{2} \int_{p(M \cap \Gamma)} \frac{ |\nabla^{\prime}u(z)|^{2}}{\sqrt{1 + | \nabla^{\prime}u(z)|^{2}}} \dif \cH^{n-1}(z).
\end{align*}
Since $\lip(u) < 1$ it follows that
\begin{equation} \label{e:67}
\int_{p(M \cap \Gamma)} | \nabla^{\prime}u(z)|^{2} \le 2^{\frac{3}{2}} e_{n}(0,1).
\end{equation}
On the other hand, $\lip(u) < 1$ and \eqref{e:54} imply
\begin{equation} \label{e:68}
\int_{p(M \Delta \Gamma)} | \nabla^{\prime}u|^{2} \le \cH^{n-1}(p(M \Delta \Gamma)) \le \cH^{n-1}(M \Delta \Gamma) \le C_{3} e_{n}(0,13).
\end{equation}

Since $e_{n}(0,1) \le 13^{n-1} e_{n}(0,13)$, \eqref{e:67} and \eqref{e:68} together guarantee \eqref{e:55}.

Step 4: Note that \eqref{e:61} and the definition of $u_{\pm}$ in \eqref{e:up} and \eqref{e:um} ensure \eqref{e:LipEnv} holds. So we conclude by showing \eqref{e:49.1}. In fact, if $M_{0}$ were closed, then \eqref{e:61} would immediately verify \eqref{e:49.1}. 

In case $M_{0}$ is is not closed, fix $\epsilon > 0$ small. For $x \in M \setminus M_{0}$ choose $y \in M_{0}$ such that $\dist(p(x), p(y)) \le \dist(x, p(M_{0})) + \epsilon$. Then,
\begin{align*}
|q(x) - u(p(x))| & \le u_{+}(p(x)) - u_{-}(p(x)) \\
& \le \left( u(p(y)) + L |p(x) - p(y)|\right) - \left(u(p(y)) - L|p(x) - p(y)| \right) \\
& \le 2 L |p(x) - p(y)| \\
& \le 2 L \dist(x, p(M_{0})) + 2 L \epsilon.
\end{align*}
Taking $\epsilon \to 0$ verifies \eqref{e:49.1}.
\end{proof}

\section{Approximation of UR domains with doubly local two-sided corkscrews}\label{URapprox.sect}

In this appendix we will build UR domains\footnote{Recall that in  \cite{hofmann2010singular} they use the word domain to mean an open set and we have adopted this convention only in the context of ``UR domains".} which (locally) approximate 
open sets satisfying a (doubly) local two-sided corkscrew (DLTSCS) condition with Ahlfors regular boundary.
This will allow us to directly use the work of \cite{hofmann2010singular} on singular integrals on UR domains.

\begin{definition}[Doubly local two-sided corkscrew condition]\label{DLTSCS.def}
Let $R_0 \in (0,\infty)$, $M_0 \ge 2$ and $x_0 \in \rn$. We say an open set $\om \subset \rn$, with $x_0 \in \pom$ satisfies the $(x_0,M_0,R_0)$-doubly local two-sided corkscrew condition or {\bf $(x_0,M_0,R_0)$-DLTSCS condition}, if for every $x \in B(x_0, R_0) \cap \pom$ and $r \in (0, R_0)$ there exist two points $X_1,X_2$ such that $B(X_1,r/M_0) \subset B(x,r) \cap \om$ and $B(X_2,r/M_0) \subset B(x,r) \setminus \overline{\om}$.
\end{definition}

The first step in the construction is to introduce the appropriate notion of boundary ``cubes" for sets with $(n-1)$-dimensional Ahlfors regular boundary. These constructions were introduced in the work of David \cite{david1988morceaux} and were refined by Christ \cite{Christ90}. The dyadic ``families" built later by Hyt\"onen and Kairema in \cite{Hyt2012}
are better adapted to our needs, thus we describe them below.

\begin{lemma}[Dyadic cubes \cite{david1988morceaux, Christ90, Hyt2012}]\label{Dcubes.lem}
Suppose $E\subset \rn$ is an $(n-1)$-dimensional, closed Ahlfors regular set. Then there exist $N, a_0, \gamma, C_2$ and $C_3$ depending on $n$ and the Ahlfors regularity constant such that the following holds. For each $t \in \{1,\dots, N\}$ there exists a  collection of Borel sets (``cubes'')
$$
\mathbb{D}^{t}_{k}(E) \defeq \mathbb{D}^t_k :=\{Q_{j}^k\subset E: j\in \mathfrak{I}_k\},$$ where
$\mathfrak{I}_k$ denotes some (possibly finite) index set depending on $k$, satisfying

\begin{list}{$(\theenumi)$}{\usecounter{enumi}\leftmargin=.8cm
\labelwidth=.8cm\itemsep=0.2cm\topsep=.1cm
\renewcommand{\theenumi}{\roman{enumi}}}

\item $E=\cup_{j}Q_{j}^k\,\,$ for each
$k\in{\mathbb Z}$.

\item If $m\geq k$ then either $Q_{i}^{m}\subset Q_{j}^{k}$ or
$Q_{i}^{m}\cap Q_{j}^{k}=\emptyset$.

\item For each $(j,k)$ and each $m<k$, there is a unique
$i$ such that $Q_{j}^k\subset Q_{i}^m$.

\item $\diam\big(Q_{j}^k\big)\leq C_2 2^{-k}$.

\item Each $Q_{j}^k$ contains some ``surface ball'' $\Delta \big(x^k_{j},a_02^{-k}\big):=
B\big(x^k_{j},a_02^{-k}\big)\cap E$.

\item $\cH^{n-1}\big(\big\{x\in Q^k_j:{\rm dist}(x,E\setminus Q^k_j)\leq \varrho \,2^{-k}\big\}\big)\leq
C_2\,\varrho^\gamma\,\cH^{n-1}\big(Q^k_j\big),$ for all $k,j$ and for all $\varrho\in (0,a_0)$.

\item For every surface ball $\Delta(x,r) = B(x,r) \cap E$, $x \in E$ and $r \in (0,\diam E)$ there exists $t$ and $Q \in \mathbb{D}^t \defeq  \cup_{k} \mathbb{D}^{t}_{k}$ with $B \subset Q$ and $\diam(Q) \le C_3r$.
\end{list}

If $Q \in \mathbb{D}^t_k$ for some some $t \in \{1,\dots, N\}$ and $k \in \mathbb{Z}$ we set $\ell(Q) = 2^{-k}$. Evidently, $\diam(Q) \approx \ell(Q)$, provided $2^{-k} \lesssim \diam(E)$
\footnote{We ignore the cubes for which, $2^{-k} \gg \diam(E)$, because $(v)$ implies that eventually $\mathbb{D}_k^t$ consists of a single cube if $\diam(E) < \infty$ and $k$ is sufficiently large.}, and we refer to $\ell(Q)$ as the ``side length" of $Q$.

\end{lemma}

\begin{remark}\label{goodchoicecube.rmk} When we use these dyadic cubes we always start knowing that the DLTSCS condition holds on some ball $B(x_0,R_0)$. The flexibility of the families (the index $t$ above) allows us to use property ($vii$) to find a cube $Q$ such that $B(x_0, C_3^{-1}R_0) \cap \pom\subset Q \subset B(x_0,R_0)\cap\pom$. 
\end{remark}

From this point onward, we work with $E \subset \rn$, an $(n-1)$-dimensional Ahlfors regular set ($E$ will eventually be the boundary of an open set) and a particular dyadic grid $\mathbb{D}:= \mathbb{D}^t$ for some $t$ to be chosen when needed to ensure the existence of a cube as in Remark \ref{goodchoicecube.rmk}. There will be no constants that depend on $t$.

For $E\subset \rn$ an $(n-1)$-dimensional Ahlfors regular set, we denote by $\W = \W(E^c)$ the collection of (closed) $n$-dimensional dyadic Whitney cubes of $\rn \setminus E$, that is the collection $\W = \{I\}$ form a pairwise non-overlapping (their boundaries may intersect) covering of $\rn \setminus E$ with the property that
\[4\diam(I) \le \dist(4I, E) \le \dist(I,E) \le 40 \diam(I),\]
(see \cite[Chapter VI]{Stein70}). Moreover, whenever $I_1,I_2 \in \W$ with $I_1 \cap I_2 \neq \emptyset$ 
\[\diam(I_1) \approx \diam(I_2).\]
 For $I \in \W$ we let $\ell(I)$ denote the side length of $I$. 

Now we relate these two notions of cubes, to form Carleson and Whitney-type regions associated to each boundary cube $Q$. These are almost exactly as in \cite{hofmann2014}\footnote{The difference here is that the regions are not `augmented' by exploiting connectivity which was present in \cite{hofmann2014}.}.

We let $K \gg 1$ be a large parameter and for $Q \in \mathbb{D}(E)$ we define
\[\W_Q := \W_Q(K) := \{I \in \W(E^c): K^{-1} \ell(Q) \le \ell(I) \le K\ell(Q), \dist(I, Q) \le K\ell(Q)\}.\]
Since $E$ is Ahlfors regular, one can show that $\W_Q$ is non-empty provided $K$ is chosen large enough. We do not fix $K$ at this point because we will eventually set $E = \pom$ and want to choose $K$ to take advantage of the existence of the (local) corkscrew points afforded by the DLTSCS condition. 

Next we fix $\tau$ a small parameter depending on dimension so that the $(1 + \tau)$-dilates of $I \in \W$, $I^*:= I^*(\tau) = (1+\tau)I$ maintain the Whitney property 
\[\ell(I) \approx \ell(I^*) \approx \dist(I^*, E) \approx \dist(I,E) \]
and $I^*$ meets $J^*$ if and only if $I \cap J \neq \emptyset$. We also may ensure (by choice of $\tau$ small) that if $I \cap J \neq \emptyset$ and $I \neq J$ then $I^* \cap (\tfrac{3}{4}J) = \emptyset$.  

Finally, we define the {\bf Whitney regions relative to $Q$}
\begin{equation}\label{tt-100}
U_Q(K) := \bigcup_{I \in \W_Q(K)} I^*
\end{equation}
and the {\bf Carleson boxes relative to Q}
\begin{equation}\label{tt-101}
T_Q(K): = \interior\left(\bigcup_{Q' \in \mathbb{D}_Q} U_{Q'}(K) \right),
\end{equation}
where $\mathbb{D_{Q}}: = \{Q' \in \mathbb{D} : Q' \subseteq Q\}$.

Now we are ready to state our approximation lemma.
\begin{lemma}\label{URdomApprox.lem}
Let $M_0 \ge 2$ and $R_0 > 0$. 
If $\om \subset \rn$ is an open set with $(n-1)$-dimensional Ahlfors regular boundary $\pom$ satisfying $\partial_*\om = \pom$ with $x_0 \in \pom$ such that $\om$ satisfies the $(x_0,M_0, 2R_0)$-DLTSCS condition, then there exist $K \gg 1$ and $M_0' \ge M_0$ depending on $n, R_0, M_0$ and the Ahlfors regularity constant such that the following holds.

Let $E = \pom$, $\mathbb{D}(E)$, $\W = \W(E^c)$, etc. be as above. Suppose $Q \in \mathbb{D}^t$ for some $t$ such that $ B(x_{0}, C_{3}^{-1} R_{0}) \cap \partial \Omega \subseteq Q \subseteq B(x_0,R_0)$\footnote{See Remark \ref{goodchoicecube.rmk}.}, then the sets 
\[T_Q^+ := T_Q^+(K):= T_Q(K) \cap \om\]
and 
\[T_Q^- := T_Q^-(K):= T_Q(K) \cap (\overline{\om})^c\]
are non-empty. They satisfy the $(M_0', \ell(Q))$-two sided corkscrew condition (see Definition \ref{tscs.def}) and $\partial T_Q^\pm$ are $(n-1)$-Ahlfors regular with constant depending on $M_0, R_0$ and the Ahlfors regularity constant for $\pom$. In particular, $T_Q^\pm$ are UR domains with constants depending on $n, R_0, M_0$ and the Ahlfors regularity constant for $\pom$\footnote{See the discussion following Definition \ref{URdom} and note that since $\diam(T_Q) \approx_K \ell(Q)$, $T_Q$ satisfies the two-sided corkscrew condition.}, and 
\[\partial T_Q^\pm \cap Q = Q. \]
 Moreover, for $\cH^{n-1}$-a.e. $x \in Q$ the measure theoretic outer normals to $T_Q^\pm$, denoted by $\nu_{T_Q^\pm}(x)$, exist and satisfy
\[\nu_{T_Q^\pm}(x) = \pm \nu_\om(x).\]
\end{lemma}

\begin{proof} Fix $Q \subseteq B(x_0, R_0)$. We choose $K$ big enough to ensure that for $Q' \in \mathbb{D}_{Q}$ with $Q' \subseteq B(x_0,R_0)$ the sets $U_{Q'}^+ := U_{Q'}^+(K):= U_{Q'}(K) \cap \om$ and $U_{Q'}^- := U_{Q'}^-(K):= U_{Q'}(K) \cap (\overline{\om})^c$ 
are non-empty. To see that such a choice (depending on $M_0, R_0$ and Ahlfors regularity constant for $\pom$) exists, we note that if $x \in Q' \subseteq B(x_0,R_0)$ then necessarily $\ell(Q') \le CR_0$ and the ball $B(x_{Q^{\prime}}, \tfrac{1}{C}\ell(Q'))$ contains two corkscrew points, one for $\om$ and one for $(\overline \om)^c$. Choosing $K^{-1} \ll 1/(CM_0)$ ensures that these points are contained in $U_Q(K)$. 

We also have that $\partial T_Q^\pm$ are both Ahlfors regular by the work of \cite{hofmann2014} (see the Appendix therein). It is also easy to see that $\partial T_Q^\pm \cap Q = Q$, since for every $x \in Q$, $x \in Q_j \in \mathbb{D}_Q$ with $\ell(Q_j) \to 0$ as $j \to \infty$. Using that $U_{Q_j}^\pm$ are non-empty we see that there exist $X_j \in U_{Q_j} \to x$ as $j \to \infty$ and hence $x \in \partial T_Q^\pm$ (see \eqref{tt-100} and \eqref{tt-101}).

Next, we show that $T_Q^\pm$ both satisfy the $(M_0', \ell(Q))$-two sided corkscrew condition. The hypotheses are symmetric so we may just show $T_Q^+$ satisfies the $(M_0', \ell(Q))$-two sided corkscrew condition. To this end, let $x \in \partial T_Q^+$ and $r \in (0, \ell(Q))$ and fix $A_0$ to be chosen\footnote{Note that the choice of $A_0$ depends on $K$, which is now fixed.}. We break into cases, following closely \cite{hofmann2014,HofMarMay}.

{\bf Case 1:} $r < A_0 \delta(x)$, where $\delta(x) := \dist(x, \pom)$. 
In this case, $\delta(x) > 0$ and $x$ is `far' from $\pom$. Necessarily (since $\delta(x) > 0$), $x \in \partial I^*$ for some `fat' Whitney cube $I^*$ with $\interior(I^*)  \subset T_Q^+$ and also $x \in J$ for some $J \in \W \setminus (\cup_{Q' \in \mathbb{D}_Q} \W_{Q'})$. The Whitney property of $I^*$ and $J$ yields $\ell(I^*) \approx \ell(J) \approx \delta(x) \gtrsim r/A_0$. It follows (from our choice of $\tau$) that $J$ contains an exterior corkscrew point and $I^*$ contains an interior corkscrew point for $T_Q^+$ at $x$ at scale $r$, with constants depending on $A_0$, for now. 

{\bf Case 2:} $r \ge A_0 \delta(x)$. In this case, we are close enough to the boundary so that we may exploit the $(M_0, R_0)$-DLTSCS condition for $\om$. We break into further cases.

{\bf Case 2a:} $\delta(x) > 0$. In this case $x \in \partial I^*$ for some $I$ as in Case 1. Let $\hat{x} \in \overline{Q}$ be such that $\delta(x) \approx |x - \hat{x}|$, where the implicit constants depend on $K$ (which we have fixed). Note that the existence of $\hat{x}$ is afforded by the Whitney property of $I^*$. Moreover, $I \in \W_{Q'}$ for some $Q' \subset Q$. Since
\[|x - \hat{x}| \le C_K \delta(x) \le C_K r/A_0 < C_K\ell(Q)/A_0,\]
choosing $A_0$ large enough we may find $Q^*$ whose closure contains $\hat{x}$,  $Q^* \subset Q$
and
\[\ell(Q^*) \approx r/A_0,\]
where the implicit constants depend on $n$, the Ahlfors regularity constant and $K$. Note that by the $(x_0, M_0,2R_0)$-DLTSCS condition of $\Omega$, and choice of $K$,  $U_{Q^*}^\pm$ are both non-empty, we may find two points $X_{Q^*}^\pm \in U_{Q^*}^\pm$ with 
\[\dist(X_{Q^*}^\pm, \partial T_Q^+) \ge C_K \ell(Q^*) \approx r/A_0.\]
Here one may take each $X_{Q^*}^\pm$ to be the center of a Whitney cube in $\W_{Q^*}$. We then choose $A_0 \gg 2$ such that
\[|x - X_{Q^*}^\pm| \le |x - \hat{x}| + |\hat{x} - X_{Q^*}^\pm| \lesssim r/A_0 < r/2.\]
Having fixed such a $A_0$, depending on the allowable parameters, we have
\[\dist(X_{Q^*}^\pm, \partial T_Q^+) \ge C_K \ell(Q^*) \gtrsim r\] 
so that $X_{Q^*}^\pm$ may serve as interior and exterior corkscrews (resp.) for $T_Q^+$ at $x$ at scale $r$.

{\bf Case 2b:} $\delta(x) = 0$. In this case, things are easier than Case 2a, provided we can show $x \in \overline{Q}$. Indeed, we may forgo the step of finding $\hat{x}$ above, by setting $\hat{x} = x$ and repeating the above argument verbatim. To show $x \in \overline{Q}$, we use that $\delta(x) = 0$ and $x \in \partial T_Q^+$ so there exists a sequence of points $X_i \in U_{Q_i}^+$ with $Q_i \subset Q$ and $\ell(Q_i)\to 0$, $|X_i - x| \to 0$ as $i \to \infty$. Here we used $\delta(X_i) \approx \ell(Q_i)$ by the Whitney property of cubes in $\W_{Q_i}$ and that $\delta(\cdot)$ is continuous. Moreover, for each $i$ there exists $\widehat{X_i} \in Q_i$ with $|\widehat{X_i} - X_i| \lesssim \ell(Q_i)$ so that
\[|x - \widehat{X_i}| \le |x - X_i| + |X_i - \widehat{X_i}| \to 0 \text{ as } i \to \infty.\]
Since $\widehat{X_i} \in Q$ this shows $x \in \overline{Q}$ and we can proceed as in Case 2a.

Again by \cite[Theorem 1]{david1990lipschitz}, an open set with Ahlfors regular boundary that satisfies a two-sided corkscrew condition on scales up to its diameter is a UR domain.
Thus, the only thing left to do is show that the measure theoretic unit normals for $T_Q^\pm$ agree with the unit normal of $\om$ up to a sign. Again, the symmetry of the hypotheses in the theorem and the fact that $\partial_*\om = \pom$ allow us only consider $T_Q^+$.

Since $T_Q^+$ have $(n-1)$-Ahlfors regular boundary and satisfy the two-sided \\
 corkscrew condition, Federer's criteria ensures that $T_Q^+$ is a set of locally finite perimeter
\cite[Theorem 1,  Section 5.11]{evans1992measure}. The structure theorem for sets of locally finite perimeter ensures that the measure 
theoretic unit normal to $\partial T_Q^+$ exists $\cH^{n-1}$-a.e. \cite[Theorem 2,  Section 5.7.3]{evans1992measure}. {  Since $Q \subset \pom$ and $\partial T_Q^+\cap Q=Q$ the measure theoretic tangents to $\partial T_Q^+$ and $\pom$ must agree $\mathcal H^{n-1}$-a.e in $Q$.  Thus the measure theoretic outer unit normal for $T_Q^+$ and $\om$ must agree up to a sign for $\cH^{n-1}$-a.e. $x \in Q$.}

To show that $\nu_{T_Q^+}(x) = \nu_\om(x)$ for $\cH^{n-1}$ a.e. in $Q$, assume that $x \in \partial^*T_Q \cap Q$ then $\nu_{T_Q^+}(x) = \pm \nu_\om(x)$. Suppose, for the sake of obtaining a contradiction, that 
$\nu_{T_Q^+}(x) = - \nu_\om(x)$ and set 
\[H^+:=\{y \in \rn: (y-x) \cdot \nu_{\om}(x) \ge 0\}.\] 
This is a half-space through $x$, perpendicular to $\nu_{\om}(x)$. The blow-up of the reduced boundary \cite[Section 5.7, Corollary 1]{evans1992measure}
gives
\[\lim_{r \to 0^+} \frac{\mathcal{L}^n(B(x,r) \cap \om \cap H^+)}{\mathcal{L}^n(B(x,r))} = 0, \]
which of course implies 
\[\lim_{r \to 0^+} \frac{\mathcal{L}^n(B(x,r) \cap T_Q^+ \cap H^+)}{\mathcal{L}^n(B(x,r))} = 0. \]
On the other hand, using $\nu_{T_Q^+}(x) = - \nu_\om(x)$, and applying \cite[Section 5.7, Corollary 1]{evans1992measure} to the set $T_Q^+$ gives
\[\lim_{r \to 0^+} \frac{\mathcal{L}^n(B(x,r) \cap T_Q^+ \cap H^+)}{\mathcal{L}^n(B(x,r))} = 1/2, \]
which is impossible. Therefore $\nu_{T_Q^+}(x) =  \nu_\om(x)$ and we have proved the lemma.
\end{proof}

\bibliographystyle{alpha}
\bibdata{references}
\bibliography{references}

\begin{thebibliography}{{Sem}91b}

\bibitem[ABKY11]{ABKY}
Daniel Aalto, Lauri Berkovits, Outi~Elina Kansanen, and Hong Yue.
\newblock John-{N}irenberg lemmas for a doubling measure.
\newblock {\em Studia Math.}, 204(1):21--37, 2011.

\bibitem[Bad12]{badger2012null}
Matthew Badger.
\newblock Null sets of harmonic measure on {NTA} domains: {L}ipschitz
  approximation revisited.
\newblock {\em Math. Z.}, 270(1-2):241--262, 2012.

\bibitem[BE17]{bortzengelstein}
Simon Bortz and Max Engelstein.
\newblock Reifenberg flatness and oscillation of the unit normal vector.
\newblock ArXiv Preprint, 08 2017.

\bibitem[BH16]{bortz2017singular}
Simon Bortz and Steve Hofmann.
\newblock A singular integral approach to a two phase free boundary problem.
\newblock {\em Proc. Amer. Math. Soc.}, 144(9):3959--3973, 2016.

\bibitem[Chr90]{Christ90}
Michael Christ.
\newblock A {$T(b)$} theorem with remarks on analytic capacity and the {C}auchy
  integral.
\newblock {\em Colloq. Math.}, 60/61(2):601--628, 1990.

\bibitem[Dav88]{david1988morceaux}
Guy David.
\newblock Morceaux de graphes lipschitziens et int\'{e}grales singuli\`eres sur
  une surface.
\newblock {\em Rev. Mat. Iberoamericana}, 4(1):73--114, 1988.

\bibitem[Dav91]{david1991singular}
Guy David.
\newblock {\em Wavelets and singular integrals on curves and surfaces}, volume
  1465 of {\em Lecture Notes in Mathematics}.
\newblock Springer-Verlag, Berlin, 1991.

\bibitem[DG61]{DG61}
Ennio De~Giorgi.
\newblock {\em Frontiere orientate di misura minima}.
\newblock Seminario di Matematica della Scuola Normale Superiore di Pisa,
  1960-61. Editrice Tecnico Scientifica, Pisa, 1961.

\bibitem[DJ90]{david1990lipschitz}
G.~David and D.~Jerison.
\newblock Lipschitz approximation to hypersurfaces, harmonic measure, and
  singular integrals.
\newblock {\em Indiana Univ. Math. J.}, 39(3):831--845, 1990.

\bibitem[DS98]{david1998quasiminimal}
Guy David and Stephen Semmes.
\newblock Quasiminimal surfaces of codimension {$1$} and {J}ohn domains.
\newblock {\em Pacific J. Math.}, 183(2):213--277, 1998.

\bibitem[EG92]{evans1992measure}
Lawrence~C. Evans and Ronald~F. Gariepy.
\newblock {\em Measure theory and fine properties of functions}.
\newblock Studies in Advanced Mathematics. CRC Press, Boca Raton, FL, 1992.

\bibitem[HK12]{Hyt2012}
Tuomas Hyt\"{o}nen and Anna Kairema.
\newblock Systems of dyadic cubes in a doubling metric space.
\newblock {\em Colloq. Math.}, 126(1):1--33, 2012.

\bibitem[HLMN17]{hofmann2017weak}
Steve Hofmann, Phi Le, Jos\'{e}~Mar\'{i}a Martell, and Kaj Nystr\"{o}m.
\newblock The weak-{$A_\infty$} property of harmonic and {$p$}-harmonic
  measures implies uniform rectifiability.
\newblock {\em Anal. PDE}, 10(3):513--558, 2017.

\bibitem[HM14]{hofmann2014}
Steve Hofmann and Jos\'{e}~Mar\'{i}a Martell.
\newblock Uniform rectifiability and harmonic measure {I}: {U}niform
  rectifiability implies {P}oisson kernels in {$L^p$}.
\newblock {\em Ann. Sci. \'{E}c. Norm. Sup\'{e}r. (4)}, 47(3):577--654, 2014.

\bibitem[HM15]{hofmann2015uniform}
Steve Hofmann and JM~Martell.
\newblock Uniform rectifiability and harmonic measure iv: Ahlfors regularity
  plus poisson kernels in $ l^{p}$ implies uniform rectifiability.
\newblock {\em arXiv preprint arXiv:1505.06499}, 2015.

\bibitem[HMM16]{HofMarMay}
Steve Hofmann, Jos\'{e}~Mar\'{i}a Martell, and Svitlana Mayboroda.
\newblock Uniform rectifiability, {C}arleson measure estimates, and
  approximation of harmonic functions.
\newblock {\em Duke Math. J.}, 165(12):2331--2389, 2016.

\bibitem[HMT10]{hofmann2010singular}
Steve {Hofmann}, Marius {Mitrea}, and Michael {Taylor}.
\newblock {Singular integrals and elliptic boundary problems on regular
  Semmes-Kenig-Toro.}
\newblock {\em {Int. Math. Res. Not.}}, 2010(14):2567--2865, 2010.

\bibitem[JK82]{jerison1982boundary}
David~S. Jerison and Carlos~E. Kenig.
\newblock Boundary behavior of harmonic functions in nontangentially accessible
  domains.
\newblock {\em Adv. in Math.}, 46(1):80--147, 1982.

\bibitem[KT97]{kenig1997harmonic}
Carlos~E. Kenig and Tatiana Toro.
\newblock Harmonic measure on locally flat domains.
\newblock {\em Duke Math. J.}, 87(3):509--551, 1997.

\bibitem[KT99]{kenig1999free}
Carlos~E. Kenig and Tatiana Toro.
\newblock Free boundary regularity for harmonic measures and {P}oisson kernels.
\newblock {\em Ann. of Math. (2)}, 150(2):369--454, 1999.

\bibitem[KT03]{kenig2003poisson}
Carlos~E. Kenig and Tatiana Toro.
\newblock Poisson kernel characterization of {R}eifenberg flat chord arc
  domains.
\newblock {\em Ann. Sci. \'{E}cole Norm. Sup. (4)}, 36(3):323--401, 2003.

\bibitem[KT06]{kenigtorotwophase}
Carlos~E. Kenig and Tatiana Toro.
\newblock Free boundary regularity below the continuous threshold: 2-phase
  problems.
\newblock {\em J. Reine Angew. Math.}, 596:1--44, 2006.

\bibitem[Mag12]{maggi2012sets}
Francesco Maggi.
\newblock {\em Sets of finite perimeter and geometric variational problems},
  volume 135 of {\em Cambridge Studies in Advanced Mathematics}.
\newblock Cambridge University Press, Cambridge, 2012.
\newblock An introduction to geometric measure theory.

\bibitem[Mat95]{Mattila}
Pertti Mattila.
\newblock {\em Geometry of sets and measures in {E}uclidean spaces}, volume~44
  of {\em Cambridge Studies in Advanced Mathematics}.
\newblock Cambridge University Press, Cambridge, 1995.
\newblock Fractals and rectifiability.

\bibitem[Mer16a]{jessica1}
Jessica Merhej.
\newblock On the geometry of rectifiable sets with carleson and poincar\'e-type
  conditions.
\newblock preprint, \text{arXiv:1510.05056}. To appear in Indiana Univ. Math.
  J., 2016.

\bibitem[Mer16b]{jessica2}
Jessica Merhej.
\newblock Poincar\`e-type inequalities and finding good parameterizations.
\newblock preprint, \text{arXiv:1605.07655}., 2016.

\bibitem[MMV96]{mattila1996cauchy}
Pertti {Mattila}, Mark~S. {Melnikov}, and Joan {Verdera}.
\newblock {The Cauchy integral, analytic capacity, and uniform rectifiability.}
\newblock {\em {Ann. Math. (2)}}, 144(1):127--136, 1996.

\bibitem[NTV14]{nazarov2014uniform}
Fedor {Nazarov}, Xavier {Tolsa}, and Alexander {Volberg}.
\newblock {On the uniform rectifiability of AD-regular measures with bounded
  Riesz transform operator: the case of codimension 1.}
\newblock {\em {Acta Math.}}, 213(2):237--321, 2014.

\bibitem[PT19]{Prats-Tolsa}
Mart{\'\i} {Prats} and Xavier {Tolsa}.
\newblock {The two-phase problem for harmonic measure in VMO}.
\newblock preprint, \text{arXiv:1904.00751}., 2019.

\bibitem[Rei60]{reifenberg1960solution}
E.~R. Reifenberg.
\newblock Solution of the {P}lateau {P}roblem for {$m$}-dimensional surfaces of
  varying topological type.
\newblock {\em Acta Math.}, 104:1--92, 1960.

\bibitem[{Sem}91a]{semmes1}
Stephen {Semmes}.
\newblock {Chord-arc surfaces with small constant. I.}
\newblock {\em {Adv. Math.}}, 85(2):198--223, 1991.

\bibitem[{Sem}91b]{semmes2}
Stephen {Semmes}.
\newblock {Chord-arc surfaces with small constant. II: Good parametrizations.}
\newblock {\em {Adv. Math.}}, 88(2):170--199, 1991.

\bibitem[Sem91c]{semmes3}
Stephen Semmes.
\newblock Hypersurfaces in {${\bf R}^n$} whose unit normal has small {BMO}
  norm.
\newblock {\em Proc. Amer. Math. Soc.}, 112(2):403--412, 1991.

\bibitem[ST89]{ST89}
Jan-Olov Str\"{o}mberg and Alberto Torchinsky.
\newblock {\em Weighted {H}ardy spaces}, volume 1381 of {\em Lecture Notes in
  Mathematics}.
\newblock Springer-Verlag, Berlin, 1989.

\bibitem[Ste70]{Stein70}
Elias~M. Stein.
\newblock {\em Singular integrals and differentiability properties of
  functions}.
\newblock Princeton Mathematical Series, No. 30. Princeton University Press,
  Princeton, N.J., 1970.

\bibitem[Tor97]{toronotices}
Tatiana Toro.
\newblock Doubling and flatness: geometry of measures.
\newblock {\em Notices Amer. Math. Soc.}, 44(9):1087--1094, 1997.

\bibitem[Tor19]{toronotices-2}
Tatiana Toro.
\newblock Geometric measure theory---recent applications.
\newblock {\em Notices Amer. Math. Soc.}, 66(4):474--481, 2019.

\end{thebibliography}

\end{document}